\documentclass{amsart} 
\usepackage{amssymb, amscd, epsfig, supertabular, times}
 
\def\BB{{\mathbb B}} 
\def\CC{{\mathbb C}}

\def\HH{{\mathbb H}}

\def\LL{{\mathbb L}} 
 
\def\PP{{\mathbb P}} 
\def\QQ{{\mathbb Q}} 
\def\RR{{\mathbb R}}
\def\SS{{\mathbb S}}

\def\ZZ{{\mathbb Z}}

\def\G{{\Gamma}} 
\def\g{{\gamma}}

\def\mon{{\rm mon}}
\def\eps{\varepsilon}
\def\even{{\rm even}}
\def\Flat{{\rm flat}}
\def\spec{{\rm Spec}} 
 
\def\irr{{\rm irr}}

\def\refl{{\rm refl}}
 
\def\ss{{\rm ss}} 
\def\stein{{\mathfrak{St}}}

\def\hyp{{\rm hyp}}

\def\bs{\backslash} 
 
\def\pt{{\bullet}}
\def\one{{\mathbf 1}}

\def\Fcal{{\mathcal F}}
\def\Hcal{{\mathcal H}}

\def\Lcal{{\mathcal L}}
 
\def\Mcal{{\mathcal M}}
 
\def\Ocal{{\mathcal O}}
 
\def\Scal{{\mathcal S}}

\def\Vcal{{\mathcal V}}

\def\la{\langle} 
\def\ra{\rangle}
\def\half{\textstyle{\frac{1}{2}}}

\newcommand\Aff{\operatorname{Aff}}
\newcommand\ar{\operatorname{Ar}}
\newcommand\aut{\operatorname{Aut}}
\newcommand\bl{\operatorname{Bl}}
\newcommand\codim{\operatorname{codim}}

\newcommand\End{\operatorname{End}}
\newcommand\ev{\operatorname{ev}}

\newcommand\Hom{\operatorname{Hom}}

\newcommand\im{\operatorname{Im}}

\newcommand\proj{\operatorname{Proj}}
\newcommand\re{\operatorname{Re}}
\newcommand\res{\operatorname{Res}}

\newcommand\GL{\operatorname{GL}}

\newcommand\SL{\operatorname{SL}}

\newcommand\PSL{\operatorname{PSL}} 
\newcommand\U{\operatorname{U}}

\newcommand\PU{\operatorname{PU}}

\newtheorem{theorem}{Theorem}[section]
\newtheorem{lemma}[theorem]{Lemma}
\newtheorem{propdef}[theorem]{Proposition-definition}
\newtheorem{thmdef}[theorem]{Theorem-definition}

\newtheorem{proposition}[theorem]{Proposition}
\newtheorem{corollary}[theorem]{Corollary}

\newtheorem{assumptions}[theorem]{Assumptions}

\theoremstyle{definition}
\newtheorem{definition}[theorem]{Definition}
\newtheorem{example}[theorem]{Example}

\newtheorem{tconvention}[theorem]{Terminological convention}

\newtheorem{discussion}[theorem]{Discussion}

\theoremstyle{remark}
\newtheorem{remark}[theorem]{Remark}
\newtheorem{remarks}[theorem]{Remarks}
\newtheorem{question}[theorem]{Question}

\begin{document}
\title{Geometric structures on the  complement of a
projective arrangement}

\author{Wim Couwenberg}
\address{Korteweg-de Vries Instituut voor Wiskunde, 
Plantage Muidergracht 24, NL-1018 TV Amsterdam,
Nederland}\email{wcouwenb@science.uva.nl}
\author{Gert Heckman}
\address{Mathematical Institute, University of
Nijmegen, P.O. Box 9010, NL-6500 GL Nijmegen, Nederland}
\email{heckman@math.kun.nl}
\author{Eduard Looijenga}
\address{Faculteit Wiskunde en Informatica, Universiteit 
Utrecht, P.O. Box 80.010, NL-3508 TA Utrecht, Nederland} 
\email{looijeng@math.uu.nl}

\subjclass{Primary: 33C67; 33D80, secondary: 33C65; 20F36}
\keywords{Dunkl connection, hypergeometric function, ball quotient}

\begin{abstract}
Consider a complex projective space with its Fubini-Study metric.
We study certain one parameter deformations of this metric on the complement
of an arrangement ($=$ finite union of hyperplanes) whose Levi-Civita
connection is of Dunkl type. Interesting 
examples are obtained from the arrangements defined by finite 
complex reflection groups. We determine a parameter interval for which the 
metric is locally of Fubini-Study type, flat, or complex-hyperbolic. 
We find a finite subset of this interval for which 
we get a complete orbifold or at least a Zariski open subset thereof, 
and we analyze these cases in 
some detail (e.g., we determine their orbifold fundamental group). 

In this set-up, the principal results of Deligne-Mostow 
on the Lauricella hypergeometric differential equation and work of
Barthel-Hirzebruch-H\"ofer on arrangements in a projective plane appear 
as special cases. Along the way we produce in a geometric manner 
all the pairs of complex reflection groups with isomorphic discriminants, 
thus providing a uniform approach to work of Orlik-Solomon. 
\end{abstract}
\maketitle

\hfill{\emph{In memory of Peter Slodowy (1948--2002)}}

\section*{Introduction}
This article wants to be the child of two publications which saw the light of day 
in almost the same year. One of them is the book by Barthel-Hirzebruch-H\"ofer (1987)
\cite{barthel}, which, among 
other things, investigates Galois coverings of $\PP^2$ that ramify in a specified manner 
over a given configuration of lines and characterizes the ones  
for which a universal such cover is a complex ball (and thus make $\PP^2$ appear 
as a---perhaps compactified---ball quotient).  
The other is a long paper by Deligne and Mostow (1986) \cite{delmost1}, which 
completes the work of Picard and Terada on the Lauricella functions and which leads to a ball 
quotient structure on $\PP^n$ relative to a hyperplane configuration of type $A_{n+1}$. 
Our reason for claiming such a descendence is that we develop a higher 
dimensional  generalization of the work by Hirzebruch et al.\  in such a manner that it 
contains the cited work of Deligne-Mostow as a special case. In other words, 
this paper's subject matter is projective arrangements which can be  understood 
as discriminants of geometric orbifold structures. 
Our approach yields new, and we believe, interesting, examples of ball quotients (which was 
the original goal) and offers at the same time a novel perspective 
on the material of the two parent papers.

It starts out quite simply with the data of a finite dimensional 
complex  inner product space $V$ in which is given a \emph{hyperplane arrangement},
that is, a finite collection of (linear) hyperplanes. We write $V^\circ$ for the complement 
of the union of these hyperplanes and $\PP (V^\circ)\subset \PP(V)$  for its 
projectivization. 
The inner product determines a (Fubini-Study) metric on $\PP(V)$ and the idea
is to deform continuously (in a rather specific manner) the restriction 
of this metric to $\PP (V^\circ)$ as to obtain a complex hyperbolic metric, i.e., 
a metric that makes $\PP (V^\circ)$ locally isometric to a complex ball. We do this in 
two stages.

We first attempt to produce a one-parameter deformation 
$\nabla^t$, $t\ge 0$, of the standard translation 
invariant connection $\nabla^0$  on (the tangent bundle) of  $V$ restricted to 
$V^\circ$ \emph{as a flat torsion free connection} on $V^\circ$.
For the reflection hyperplane arrangement of a finite
Coxeter group such a deformation is given by Dunkl's construction and we
try to imitate this.  Although this is not always possible---the existence of such
a deformation imposes strong conditions on  the arrangement---plenty of examples do exist. 
For instance, this is always possible for the 
reflection hyperplane arrangement of a complex reflection group. Besides,
it is a property that is inherited by the  arrangement that  is naturally defined
on a given intersection of members of the arrangement. 

The inner product defines a translation invariant metric on $V$. Its restriction 
$h^0$ to $V^\circ$ is obviously flat for  $\nabla^0$ and the next step is to show that  
we can deform $h^0$ as a nonzero flat hermitian form $h^t$ \emph{which is flat for} $\nabla^t$ 
(so that $\nabla^t$ becomes a Riemannian connection as long as  $h^t$ 
is nondegenerate). This is done in such a manner that scalar multiplication 
in $V$ acts locally like homothety and as a consequence, $\PP (V^\circ)$ 
inherits from $V^\circ$ a hermitian form $g^t$. For $t=0$ this gives us the Fubini-Study
metric.  We only allow $t$  to move in an interval for which $g^t$ stays positive definite. 
This still makes it possible for $h^t$ to become degenerate or of hyperbolic signature as long
as for every $p\in V^\circ$, the restriction of $h^t$ to a hyperplane 
supplementary and perpendicular to $T_p(\CC p)$ is positive definite.
If $T_p(\CC p)$ is the kernel of $h^t$ (we refer to this situation 
as the \emph{parabolic} case), then
$g^t$ is a flat metric, whereas when $h^t$ is negative on $T_p(\CC p)$
(the \emph{hyperbolic} case), $g^t$ is locally 
the metric of a complex ball. It is necessary to impose additional conditions
of a simple geometric nature in order to have a neat global picture, that is,
to have $\PP (V^\circ)$ of finite volume and realizable as a quotient of a dense open 
subset of a flat space resp.\ a  ball by a discrete group of isometries. We call these
the \emph{Schwarz conditions}, because they are reminiscent of the ones found by 
H.A. Schwarz which ensure that the Gau\ss\ hypergeometric function is algebraic.

Deligne and Mostow gave a modular interpretation of their ball quotients. Some of them
are in fact Shimura varieties and indeed, particular cases were already studied by 
Shimura and Casselman (who was then Shimura's student) in the sixties. The natural 
question is whether such an interpretation also exists for the ball quotients 
introduced here.  We know this to be the case for some of them, but
we do not address this issue in the present paper. 

We mention some related work, without however any pretension of attempted completeness. 
A higher dimensional generalization of Hirzebruch's original approach with 
Fermat covers and fixed weights along all hyperplanes and emphasizing 
the three dimensional case was developped by Hunt \cite{hunt}. His paper 
with Weintraub \cite{hw} fits naturally in our framework; their Janus-like algebraic 
varieties are exactly related to the various ramification orders $q$ allowed in the
tables of our final Section \ref{section:examples}. The articles
by  Holzapfel \cite{holz1}, \cite{holz2} and Cohen-W\"ustholz \cite{cw} 
contain applications to transcendency theory.

\medskip
We now briefly review the contents of the separate sections of this paper.
In the first section we develop a bit of the general theory of affine structures
on complex manifolds, where we pay special attention to a simple kind of 
degeneration of such a structure along a normal crossing divisor. Although it 
is for us the occasion to introduce
some terminology and notation, the reader is perhaps well-advised  
to skip this section during a first reading and use it for consultation only.

Section two  focuses on a notion which is central to this paper,  that
of a Dunkl system.  We prove various hereditary properties and we give a 
number of examples. We show in particular that the
Lauricella functions fit in this setting. In fact, in the last subsection  
we classify all the Dunkl systems whose underlying arrangement is a Coxeter
arrangement and show that the Lauricella examples exhaust the cases of type
$A$. For the other Coxeter arrangements of rank at least three  the Dunkl system
has automatically the symmetry of the corresponding Coxeter group, except
for those of type $B$, for which we essentially reproduce the Lauricella series.

The next section discusses the existence of a nontrivial hermitian form which is 
flat relative to the Dunkl connection. We prove  among other things that such a
form always exists in the case of a complex reflection arrangement and in the 
Lauricella case and we determine when this form 
is positive definite, parabolic or hyperbolic.

Section four is devoted to the Schwarz conditions. We show that when these 
conditions are satisfied, the holonomy cover extends as a ramified cover 
over an open subset of $V$ of codimension at least two, that the developing
map extends to this ramified cover, and that the latter extension becomes a 
local isomorphism if we pass to the quotient by a finite group $G$
(which acts as a complex reflection group on $V$, but lifts to the ramified cover). 
This might explain why we find it  reasonable to impose such a condition. 
From this point onward we assume such conditions satisfied  and
concentrate on the situations that really matter to  us. 

Section five deals with the elliptic and the parabolic cases. 
The elliptic case can be characterized as having finite holonomy.
It is in fact treated in two somewhat different situations: at first 
we deal with a situation where we find that $\PP (G\bs V)$ is 
the metric completion of $\PP (G\bs V^\circ)$ and acquires the structure
of an elliptic orbifold.
What makes this interesting is that this is not the natural $G$-orbifold structure that 
$\PP(G\bs V)$  has a priori: it is the structure of the quotient of a  projective space 
by the holonomy group. This is also a complex reflection group,
but usually differs from $G$. Still the two reflection groups are related  
by the fact that their discriminants satisfy a simple inclusion relation.
We prove that all pairs of complex reflection groups with isomorphic discriminants
are produced in this fashion. The other elliptic case we discuss is when the metric 
completion of $\PP (G\bs V^\circ)$ differs from $\PP(G\bs V)$ but is gotten from
the latter by means of an explicit blowup followed by an explicit blowdown. We have to deal
with such a situation, because it is one which we encounter when we treat the hyperbolic 
case. The parabolic case presents little trouble and is dealt with in a straightforward 
manner.

Our main interest however concerns the hyperbolic situation and that is saved for last. 
We first treat the case where we get a compact hyperbolic orbifold, because it is 
relatively easy and takes less than half a page. The general case is rather delicate, 
because the metric completion of $\PP (G\bs V^\circ)$ (which should be a ball quotient of 
finite volume) may differ from $\PP(G\bs V)$. Deligne and Mostow used at this point
geometric invariant theory for effective divisors on $\PP^1$, but in the present
situation this tool is not available to us and we use an argument based
on Stein factorization instead. As it is rather difficult to briefly  summarize the 
contents of our main theorem, we merely refer to \ref{thm:main} for its statement. 
It suffices to say here that it produces new examples of discrete  complex hyperbolic groups
of cofinite volume. We also discuss the implications for the allied
algebra of automorphic forms. We close this section with a presentation for the
holonomy group, which is also valid for elliptic and the parabolic cases.  

The final section tabulates the elliptic, parabolic and hyperbolic examples 
of finite volume with the property that the associated arrangement is that of a 
finite reflection group of rank at least three (without requiring it to have the 
symmetry of that group). In the hyperbolic case we mention whether the holonomy 
group is cocompact.

 \medskip
This work has its origin in the thesis by the first author \cite{couw}  at the University
of Nijmegen (1994) written under the supervision of the second author. Although that 
project went quite far in carrying out the program described above, the results 
were never formally 
published, in part, because both felt that it should be completed first. This remained
the state of affairs until around 2001, when the idea emerged that work of the
third author \cite{looij:ball} might be relevant here. After we had joined forces in 2002,
the program was not only completed as originally envisaged, but we were even able to
go well beyond that, including the adoption of a more general point of view and a 
change in perspective.

We dedicate this paper to the memory of our good friend and colleague Peter
Slodowy.

\medskip
\emph{Acknowledgements.}
Three letters written in 1994 by P.~Deligne to Couwenberg 
(dated Nov.\ 12 and Nov.\ 16) and to M.~Yoshida (dated Nov.\ 12) were quite helpful to us. 
Couwenberg thanks Masaaki Yoshida for his encouragements
and support during his 1998 visit to Kyushu University, Heckman
expresses his thanks to Dan Mostow for several stimulating discussions 
and Looijenga is grateful to the MSRI at Berkeley where he stayed the first 
three months of 2002 and where part of his contribution to this work was done (and 
via which he received NSF-support through grant DMS-9810361).

\section*{Terminological index and list of notation}
The terminological index is alphabetical, but the list of notation is by 
order of introduction.

{\small
\subsection*{Terminological index}\hfill

\smallskip\noindent
\emph{admissible} hermitian form: Definition \ref{def:admissible}\\
\emph{affine quotient}: Remark\ref{rem:simpledeg1}\\
\emph{affine structure}: Subsection \ref{subsect:affinestr}\\
\emph{apex curvature}: Subsection \ref{subsect:thurston}\\
\emph{arrangement complement}: Subsection \ref{subsect:review}\\
\emph{Artin group}:  Subsection \ref{subsect:hecke}\\
\emph{Borel-Serre extension}: Subsection \ref{subsect:borelserre}\\
\emph{co-exponent}: Subsection \ref{subsect:hypexp}\\
\emph{cone manifold}: Subsection \ref{subsect:thurston}\\
\emph{Coxeter matrix}: Subsection \ref{subsect:hecke}\\
\emph{degenerate hyperbolic form}: Subsection \ref{subsect:degenhyp}\\
\emph{developing map}: Definition \ref{def:devmap}\\
\emph{dilatation  field}: Definition \ref{def:dilatation}\\
\emph{discriminant} of a complex reflection group: Subsection \ref{subsect:hypexp}\\
\emph{Dunkl}, connection of $\sim$ type, $\sim$ form, $\sim$ system: 
Definition \ref{def:dunkl}\\
\emph{elliptic structure}: Definition \ref{def:admissible}\\
\emph{Euler field}: Corollary \ref{cor:dunklchar}\\
\emph{exponent} of a complex reflection group: Subsection \ref{subsect:hypexp}\\
\emph{fractional divisor}: Remark \ref{rem:delmost}\\
\emph{germ}: See \emph{Some notational conventions}\\
\emph{Hecke algebra}:  Subsection \ref{subsect:hecke}\\
\emph{holonomy group}: Terminological convention \ref{tconv}\\
\emph{hyperbolic exponent}: Theorem-definition \ref{thm:indefinite}\\
\emph{hyperbolic structure}: Definition \ref{def:admissible}\\
\emph{index} of a hermitian form: Lemma \ref{lemma:lauricellaform}\\
\emph{infinitesimally simple degeneration of an affine structure along  a divisor}: 
Definition \ref{def:infisimpledeg}\\
\emph{irreducible arrangement, stratum of an $\sim$, component of an $\sim$}: 
Subsection \ref{subsect:review}\\
\emph{Lauricella} connection, $\sim$ function: 
Proposition-definition \ref{def:lauricellaconn}\\
\emph{longitudinal} Dunkl connection: Definition \ref{def:transvdunkl}\\
\emph{mildly singular} function, $\sim$ differential: 
discussion preceding Lemma \ref{lemma:mild}\\
\emph{monodromy group}: Terminological convention \ref{tconv}\\
\emph{normal linearization of a hypersurface}: Definition  \ref{def:normallin}\\
\emph{nullity}: Lemma \ref{lemma:lauricellaform}\\
\emph{parabolic structure}: Definition \ref{def:admissible}\\
\emph{projective quotient}: Remark \ref{rem:simpledeg1}\\
\emph{pure degeneration}: Definition \ref{def:simpledeg}\\
\emph{pure quotient}: Remark \ref{rem:simpledeg1}\\
\emph{reflection representation}:  Subsection \ref{subsect:hecke}\\
\emph{residue} of a connection: Subsection \ref{subsect: logdeg}\\
\emph{semisimple holonomy} around a stratum: paragraph preceding
Corollary \ref{cor:affineproj}\\
\emph{simple degeneration} of an affine structure along  a divisor: 
Definition \ref{def:simpledeg}\\
\emph{Schwarz} condition, $\sim$ rotation group, $\sim$ symmetry group,
$\sim$ in codimension one: Definition \ref{def:schwarz}\\
\emph{special} subball, $\sim$ subspace: Subsection \ref{subsect:main}\\
\emph{splitting} of an arrangement: Subsection \ref{subsect:review}\\
\emph{stratum} of an arrangement: Subsection \ref{subsect:review}\\
\emph{topological Stein factorization}: paragraph preceding Lemma \ref{lemma:stein}\\
\emph{transversal} Dunkl connection: Definition \ref{def:transvdunkl}\\

 \subsection*{List of notation}
\begin{itemize}
\item[$\Aff_M$]  Subsection \ref{subsect:affinestr}: 
the local system of locally affine-linear functions on an affine manifold.
\item[$\Aff (M)$] Subsection \ref{subsect:affinestr}: the space of global
sections of $\Aff_M$.
\item[$\G$] Subsection \ref{subsect:affinestr}: the holonomy group.
\item[$A$] Subsection \ref{subsect:affinestr}: the affine space which receives
the developing map.
\item[$\res_D(\nabla)$] Subsection \ref{subsect: logdeg}: 
Residue of a connection along $D$.
\item[$\nu_{D/W}$] Lemma \ref{lemma:affineresidue}: normal bundle of $D$ in $W$. 
\item[$D_{p,0}$] Remark \ref{rem:simpledeg1}: the affine quotient of $D_p$.
\item[$D_{p,\lambda}$] Remark \ref{rem:simpledeg1}: the projective quotient of $D_p$.
\item[$W_{p,\lambda}$] Remark \ref{rem:simpledeg1}: the pure quotient of $W_p$.
\item[$V^\circ$]  Subsection \ref{subsect:review}: the complement of
an arrangement in $V$.
\item[$\Lcal (\Hcal)$]  Subsection \ref{subsect:review}: the intersection  lattice of the
arrangement $\Hcal$.
\item[$\Hcal_L$] Subsection \ref{subsect:review}: the members of $\Hcal$ containing $L$.
\item[$\Hcal^L$] Subsection \ref{subsect:review}: the intersections of the
members of $\Hcal-\Hcal_L$ with $L$.
\item[$\Lcal_\irr(\Hcal)$]  Subsection \ref{subsect:review}: the irreducible members
of $\Lcal (\Hcal)$.
\item[M(L)] Lemma \ref{lemma:dich}.
\item[$\phi_H$] Subsection \ref{subsect:affinearc}: a linear from 
which defines the hyperplane $H$.
\item[$\omega_H$]  Subsection \ref{subsect:affinearc}: the logarithmic form 
defined by the hyperplane $H$.
\item[$\nabla^0$]  Subsection \ref{subsect:affinearc}: the translation invariant
connection on an affine space.
\item[$E_V$]  Corollary \ref{cor:dunklchar}: the Euler vector field on a vector space $V$.
\item[$\pi_L$]  Subsection \ref{subs:hered}: the orthogonal projection 
in an inner product space with kernel $L$.
\item[$\kappa_L$]  Lemma \ref{lemma:hereditary}.
\item[$\nabla^\kappa$]  paragraph preceding Corollary \ref{cor:monotony}.
\item[$\Omega^\kappa$]  paragraph preceding Corollary \ref{cor:monotony}.
\item[$\CC^{\Hcal,\Flat}$] paragraph preceding Corollary \ref{cor:monotony}:
the set of exponents $\kappa$ for which $\nabla^\kappa$ is flat.
\item[$\Hcal^L_\perp$]  Discussion preceding Lemma \ref{lemma:condecomp}.
\item[$\pi_I^L$] Lemma \ref{lemma:condecomp}.
\item[$\omega_I^L$] Lemma \ref{lemma:condecomp}.
\item[$\bl_L V$] Subsection \ref{subsect:loctriv}: blowup of $V$ in $L$.
\item[$\tilde\pi_L^*$] Subsection \ref{subsect:loctriv}.
\item[$h^0$] Subsection \ref{subs:admrange}: the hermitian form defined by the 
inner product.
\item[$m_\hyp$] Theorem-definition \ref{thm:indefinite}: the hyperbolic exponent.
\item[$\Hcal (\Fcal)$] Lemma \ref{lemma:flatsubbundle}.
\item[$d_i$] Subsection \ref{subsect:hypexp}: the $i$th degree of a reflection group.
\item[$m_i$] Subsection \ref{subsect:hypexp}: the $i$th exponent of a reflection group.
\item[$d_i^*$] Subsection \ref{subsect:hypexp}: the $i$th codegree of a reflection group.
\item[$m_i^*$] Subsection \ref{subsect:hypexp}: the $i$th co-exponent of a reflection 
group.
\item[$\ar (M)$] Subsection \ref{subsect:hecke}: the Artin group attached to the 
Coxeter matrix $M$.
\item[$\Hcal (M,t)$] Subsection \ref{subsect:hecke}: the universal Hecke algebra.
\item[$R$] Subsection \ref{subsect:hecke}: a domain associated to a 
Hecke algebra.
\item[$\Hcal (M)$] Subsection \ref{subsect:hecke}: the Hecke algebra over $R$.
\item[$\rho^\mon$] Subsection \ref{subsect:hecke}: the monodromy representation
of $\Hcal (M)$.
\item[$\rho^\refl$] Subsection \ref{subsect:hecke}: the reflection representation
of $\Hcal (M)$.
\item[$N$] Subsection \ref{subsect:lauricellaflat}.
\item[$w_k$] Subsection \ref{subsect:lauricellaflat}: a complex number of norm one
attached to $\mu$.
\item[$A_w$] paragraph preceding Lemma \ref{lemma:lauricellaform}: a hyperplane
of $\RR^{n+1}$.
\item[$Q_w$] paragraph preceding Lemma \ref{lemma:lauricellaform}: a quadratic form
on $A_w$.
\item[$p_L,q_L$] Definition \ref{def:schwarz}: numerator resp.\ denominator of
$1-\kappa_L$.
\item[$G_L$] Definition \ref{def:schwarz}: Schwarz rotation group.
\item[$G$] Definition \ref{def:schwarz}: Schwarz symmetry group.
\item[$V^f$] Subsection \ref{subsect:extension}: locus of finite holonomy  in $V$.
\item[$\ev_G$] Theorem \ref{thm:devmap}: a factor of an extension of the developing map.
\item[$\Lcal^-$] Subsection \ref{subs:finiteholctd}.
\item[$\Lcal^0$] Subsection \ref{subs:finiteholctd}.
\item[$\Lcal^+$] Subsection \ref{subs:finiteholctd} and Discussion \ref{disc:devmap2}.
\item[$\BB^-$] Theorem \ref{thm:main}.
\item[$V^+$] Discussion \ref{disc:devmap}.
\item[$V^-$] Discussion \ref{disc:devmap}.
\item[$E(L)$] Discussion \ref{disc:devmap}.
\item[$D(L)$] Discussion \ref{disc:devmap}.
\item[$\SS$] Discussion \ref{disc:devmap}: the sphere of rays.
\item[$E(L_\pt)$] Discussion \ref{disc:devmap}.
\item[$S(L_\pt)$] Discussion \ref{disc:devmap} and Discussion \ref{disc:devmap2}.
\item[$\stein$] (as a subscript) paragraph preceding Lemma \ref{lemma:stein}:
formation of a Stein quotient.
\item[$\BB^+$] Subsection \ref{subsect:borelserre}: the Borel-Serre extension of $\BB$.
\item[$B^+$] Discussion \ref{disc:devmap2}.
\end{itemize}
}

\subsection*{Some notational conventions.}
If $\CC^\times$ acts on a variety $X$, then we often write $\PP (X)$ for 
the orbit space of the subspace of $X$ where $\CC^\times$ acts with finite 
isotropy groups. This notation is of course suggested by the case when
$\CC^\times$ acts by scalar multiplication on a complex vector space $V$, for
$\PP (V)$ is then the associated projective space. This example also shows that
a $\CC^\times$-equivariant map  $f:X\to Y$ may or may not induce
a morphism $\PP (f) :\PP(X)\to \PP (Y)$.

If $X$ is a space with subspaces $A$ and $Y$, then the \emph{germ} of $Y$
at $A$ is the filter of neighborhoods of $A$ in $X$ restricted to $Y$; we
denote it by $Y_A$. Informally, $Y_A$  may be thought of as an unspecified 
neighborhood of $A$ intersected with $Y$. For
instance, a \emph{map germ} $Y_A\to Z$ is given by a pair 
$(U, f:U\cap Y\to Z)$, where $U$ is some neighborhood of $A$, and another 
such pair $(U', f':U'\cap Y\to Z)$ defines the same map-germ if $f$ and $f'$ 
coincide on $U''\cap Y$ for some neighborhood $U''$ of $A$ in $U\cap U'$.

\tableofcontents

\section{Affine structures with logarithmic singularities}
We first recall a few basic properties regarding the notion of an affine
structure.

\subsection{Affine structures}\label{subsect:affinestr}
Let be given a connected complex manifold $M$ of complex dimension $n$.
An \emph{affine structure} on $M$ is an atlas (of 
complex-analytic charts) for which the transitions maps are
complex affine-linear and which is maximal for that property.
Given such an atlas, then the complex valued functions that are locally complex-affine linear
make up a local system $\Aff_M$ of $\CC$-vector
spaces in the structure sheaf $\Ocal_M$. This local system 
is of rank $n+1$ and  contains the constants $\CC_M$. The quotient
$\Aff_M/\CC_M$ is a local system whose underlying vector bundle
is the complex cotangent bundle of $M$, hence is given by a flat connection
$\nabla : \Omega_M\to \Omega_M\otimes \Omega_M$. This 
connection is torsion free, for
it sends closed forms to symmetric tensors. (This is indeed equivalent
to the more conventional definition which says that the associated 
connection on the tangent bundle is symmetric:
for any pair of local vector fields $X,Y$ on $M$, we have 
$\nabla_XY-\nabla_YX=[X,Y]$.)

Conversely, any flat, torsion free connection
$\nabla$ on the complex cotangent bundle of $M$ defines an affine structure:
the subsheaf  $\Aff_M\subset \Ocal_M$ of holomorphic functions  whose total 
differential is flat for $\nabla$ is then a local system of rank $n+1$ containing
the constants and the atlas in question consists of the charts whose components
lie in $\Aff_M$.

\begin{tconvention}\label{tconv}
With regard to a flat, torsion free connection $\nabla$ on the complex cotangent 
bundle of a connected complex manifold $M$, we reserve the term 
\emph{monodromy group} 
as the monodromy of that connection on the cotangent bundle of $M$, whereas the 
\emph{holonomy group} shall be the monodromy group of the local system 
$\Aff_M$. 
\end{tconvention}

So the holonomy group is an extension of the monodromy group by a group of translations.
In this situation one defines a developing map as follows. 
If $\G$ denotes the holonomy group, let   $\widetilde{M}\to M$ be an associated 
$\G$-covering. It is  unique up to isomorphism and it has the property that
the pull-back of $\Aff_M$ to this covering is generated by its sections. 
Then the space of affine-linear functions on $\widetilde{M}$,
$\Aff (\widetilde{M}) :=H^{0}(\widetilde{M},\Aff_{\widetilde{M}})$, is a 
$\G$-invariant vector space of holomorphic functions 
on $\widetilde{M}$. This vector space  contains the constant functions
and the quotient $\Aff (\widetilde{M})/\CC$ can be identified with the space 
of flat holomorphic differentials on $\widetilde{M}$; it has the same dimension as $M$. 
The set $A$ of linear forms $\Aff (\widetilde{M})\to \CC$
which are the identity on $\CC$ is an affine $\G$-invariant
hyperplane in $\Aff (\widetilde{M})^{*}$. 

\begin{definition}\label{def:devmap}
The \emph{developing map} of the affine structure is the evaluation mapping
$\ev :\widetilde{M}\to A$ which assigns to $\tilde z$ the linear form 
$\ev_{\tilde z}: \tilde f\in  \Aff(\widetilde{M}) \mapsto \tilde f(\tilde z)\in\CC$.
\end{definition}

Notice that this map is $\G$-equivariant and a local affine isomorphism.
In fact, it determines a natural affine atlas on $M$ whose charts take values in $A$ 
and whose transition maps lie in $\G$.

\begin{definition}\label{def:dilatation}
We call a nowhere zero holomorphic vector field $E$ on $M$ a 
\emph{dilatation  field with factor $\lambda\in\CC$}  when for every
local vector field $X$ on $M$, $\nabla_X(E)=\lambda X$. 
\end{definition}

Let us have a closer look at this property. If $X$ is flat, then
the torsion freeness yields: $[E,X]=\nabla_E(X)-\nabla_X(E)=-\lambda X$. In other words, 
Lie derivation with respect to $E$ acts on flat vector fields simply as multiplication by 
$-\lambda$. Hence it acts on flat differentials as multiplication by $\lambda$.
So $E$ acts on $\Aff_M$ with eigenvalues $0$ (on $\CC$) and $\lambda$ 
(on $\Aff_M /\CC_M$).

Suppose first that $\lambda\not=0$. Then the $f\in \Aff_M$ for which 
$E(f)=\lambda f$ make up a flat supplement of $\CC_M$ in $\Aff_M$.
This singles out a fixed point  $O\in A$ of $\G$ so that the affine-linear structure 
is in fact a linear structure and the developing map takes the lift of $E$ on
$\widetilde{M}$ to $\lambda$ times the Euler vector field on $A$ relative to $O$.
This implies that locally the leaf space of the foliation defined by 
$E$ is identified with an open set of the projective space of $(A,O)$ (which is naturally 
identified with the projective space of the space of flat vector fields on 
$\widetilde{M}$). Hence this leaf space acquires a complex projective structure.

Suppose now that $\lambda=0$. Then $\CC$ need not be a direct summand of $\Aff_M$.
All we can say is that $E$ is a flat vector field so that its lift to $\widetilde{M}$ 
maps a constant nonzero vector field on $A$. So locally the leaf space of 
the foliation defined by $E$ has an affine-linear structure defined by an atlas 
which takes values in the quotient of $A$ by the translation group generated by 
a constant vector field.

\subsection{Logarithmic degeneration}\label{subsect: logdeg}
In this subsection $W$ is a complex manifold with a given
affine structure $\nabla$ on the complement  $W-D$ of a hypersurface $D$. 
At first, $D$ is smooth connected, later we allow $D$ to have normal
crossings.

We recall that if we are given a holomorphic vector bundle $\Vcal$ on $W$,
then a flat connection  $\nabla$ on $\Vcal$ with a logarithmic pole along $D$ is a map
$\Vcal\to\Omega_W(\log D)\otimes \Vcal$ satisfying the usual properties of a
flat connection. Then the  residue map $\Omega_W(\log D)\to \Ocal_D$ 
induces an $\Ocal_D$-endomorphism $\res_D(\nabla)$ of $\Vcal\otimes\Ocal_D$, 
called the  \emph{residue of the connection}. It is well-known that the 
conjugacy class of this endomorphism is constant along $D$.  In particular, 
$\Vcal\otimes\Ocal_D$ decomposes according to the generalized
eigen spaces of $\res_D(\nabla)$. This becomes clear if
we choose at $p\in D$ a chart $(t,u_1,\dots ,u_n)$ such that $D_p$ is given 
by $t=0$: then 
$R:=t\frac{\partial}{\partial t}, U_1:=\frac{\partial}{\partial u_1},\dots , 
U_n:=\frac{\partial}{\partial u_n}$ is a set of commuting vector fields,
covariant derivation with respect to these fields preserves $\Vcal_p$
(and since $\nabla$ is flat, the resulting endomorphisms of $\Vcal_p$ pairwise 
commute) and $R$ induces in $\Vcal_p\otimes\Ocal_{D,p}$ the residue
endomorphism. In particular, the kernel of $R$ is preserved by $U_i$. The
action of $U_i$ on this kernel restricted to $D_p$ only depends on the restriction 
of $U_i$ to $D_p$. This shows that $\nabla$ induces on the kernel of the residue 
endomorphism a flat connection. (A similar argument shows that the 
projectivization of the subbundle of $\Vcal\otimes\Ocal_D$ associated to 
an eigen value of $\res_D(\nabla)$ comes 
with a projectively flat connection.)

\begin{lemma}\label{lemma:affineresidue}
Suppose that the affine structure $\nabla$ on ${W-D}$ extends 
to $\Omega_W$ with a genuine logarithmic pole. Letting $\nu_{D/W}$ stand for the 
normal bundle of $D$ in $W$, then:
\begin{enumerate}
\item[(i)] the residue of $\nabla$ on $\Omega_W$  respects the natural exact sequence
\[
0\to \nu^*_{D/W}\to \Omega_W\otimes \Ocal_D\to \Omega_D\to 0
\]
and induces the zero map in $\Omega_D$,
\item[(ii)] the connection induces in $D$ an affine structure,
\item[(iii)] the connection has a logarithmic pole on $\Omega_W(\log D)$ as well,
its residue on this sheaf  respects the exact sequence
\[
0\to \Omega_D\to \Omega_W(\log D)\otimes \Ocal_D\to \Ocal_D\to 0.
\]
and is zero on $\Omega_D$. The scalar operator in $\Ocal_D$ is one less
than the one in  $\nu^*_{D/W}$.
\end{enumerate}
\end{lemma}
\begin{proof} By assumption, $\nabla$ defines a map
$\Omega_W\to \Omega_W(\log D)\otimes\Omega_W$. Since $\nabla$ is
torsion free, this extension then takes values in 
\[
\Big(\Omega_W(\log D)\otimes\Omega_W\Big)\cap 
\Big(\Omega_W\otimes\Omega_W(\log D)\Big)\subset
\Omega_W(\log D)\otimes\Omega_W(\log D).
\]
If $t$ be a local equation of $D$, then this intersection 
is spanned by $t^{-1}dt\otimes dt$ and $\Omega_W\otimes\Omega_W$.
Hence the residue of $\nabla$ on $\Omega_W$ maps
$\Omega_W\otimes \Ocal_D$ to the span of $dt$, that is, to $\nu^*_{D/W}$.
So (i) follows. It is also clear that $\nabla$ drops to map
$\Omega_D\to\Omega_D\otimes\Omega_D$ and so (ii) follows as well.
Finally, let $R$ be a local 
vector field  with $R(t)=t$. Then $\nabla_R$ induces
the residue map and so $\nabla_R (dt)$ is of the form  $cdt+t\omega$
for some constant $c\in\CC$ and some $\omega\in \Omega_W$. It follows that 
$\nabla_R (t^{-1}dt)=(c-1)t^{-1}dt +\omega\in \Omega_W(\log D)$. 
This proves assertion (iii) .
\end{proof}

The converse is not true: if the affine 
structure extends with a logarithmic pole to $\Omega_W(\log D)$, then
it need not have that property on $\Omega_W$.  The advantage of this logarithmic 
extension (over $\Omega_W$) is that has better stability properties
with respect to blowing up. 

\begin{definition}\label{def:simpledeg}
Let $D$ be a smooth connected hypersurface in an analytic manifold $W$.
We say that an affine structure on $W-D$ has 
\emph{simple degeneration along $D$ of logarithmic exponent $\lambda\in\CC$} if 
at any $p\in D$ there exist a local equation $t$ for $D_p$ in $W_p$, a 
morphism $F_0 :W_p\to T_0$ to an affine space $T_0$ such that
\begin{enumerate}
\item[($\lambda=0$)] $(F_0,t):W_p\to T_0\times \CC$ is a local isomorphism and
there exists an affine-linear function $u:T_0\to\CC$ such that $uF_0(p)\not=0$ and
the developing map near $p$ is affine equivalent to $(F_0, \log t .(uF_0))$,
\item[($\lambda\not=0$)] there exists  a morphism
$F_1: W_p\to T_1$ to  a linear space $T_1$, such that $(F_0,t,F_1): 
W_p\to T_0\times\CC\times T_1$ is a
local isomorphism and the developing map near $p$ is 
affine equivalent to $(F_0,t^{-\lambda}, t^{-\lambda}F_1)$.
\end{enumerate}
If in the last case ($\lambda\not=0$), $T_0=0$, we say that 
the degeneration is \emph{pure}. 
\end{definition}

Before we analyze the structural implications of this property it is 
useful to have the following notion at our disposal.

\begin{definition}\label{def:normallin}
If $D$ is a smooth analytic hypersurface in an analytic manifold $W$, then
a \emph{normal linearization} of $D$ is a vector field on $W_D$ which 
is tangent to the fibers of some retraction $W_D\to D$ and has a simple
zero at $D$ with residue $1$. If we are also given an affine structure 
$\nabla$ on $W_D-D$, then we say that the normal linearization is \emph{flat}
if the vector field is an infinitesimal affine-linear transformation.
\end{definition}

It is clear that this retraction is then unique.
Note that such a vector field generates a $\CC^\times$-action on
$W_D$ with $D$ as fixed point set which preserves each fiber of the
retraction. Thus the germ $W_D$ gets identified with the germ of $D$ in its normal
bundle (in other words, an analytic version of the tubular neighborhood theorem
holds); this explains the chosen terminology. If it is flat 
with respect to a given affine structure, then the 
$\CC^\times$ action lifts to the holonomy cover as a one-parameter group
of affine-linear transformations.

\begin{remarks}[The case $\lambda=0$]\label{rem:simpledeg0}
Let us begin with noting that $\nabla$ extends to $\Omega_W(\log D)$ 
with a logarithmic singularity along $D$:  We get
\[
\nabla (\frac{dt}{t})=-\Big( \frac{dt}{t}\otimes \frac{d(uF_0)}{uF_0}+
\frac{d(uF_0)}{uF_0}\otimes\frac{dt}{t}\Big),\quad \nabla (F_0^*\alpha)=0,
\]
where $\alpha$ is any translation invariant differential on $T_0$.
We also see that the residue endomorphism of 
$\Omega_W(\log D)\otimes\Ocal_D$ preserves $\Omega_D$ and is either
trivial ($u$ is constant) or has  image a rank one subbundle of $\Omega_D$
($u$ nonconstant).
An element of $\Ocal_{D,p}$ that is the restriction of an
element of $\Ocal_{W,p}$ which is affine-linear outside $D$ is in fact the
composite of the local  isomorphism $F_0|D_p$ and an affine-linear function on $T_0$.
So $D$ has a natural affine structure and $F_0$
determines a retraction of $W_p\to D_p$
whose restriction to $W_p-D_p$ is affine. Notice that 
$t\frac{\partial}{\partial t}$ is a flat vector field which is tangent to the fibers
of this affine retraction.  It is easy 
to see that both this vector field and the retraction are canonical (independent
of our choice of coordinates).  Hence they are globally defined
and determine a flat normal linearization of $D\subset W$.
The total space  of the normal bundle deprived from its zero section comes with an 
affine structure. The holonomy  respects that structure, hence the 
holonomy group of $W_D-D$ is a central extension of the
holonomy group of the affine structure of $D$. Notice also that if we let 
$t\to 0$ in a fixed sector  (on which $\log t$ is continuous), then the  projectivization 
of the developing map tends to a singleton.  
\end{remarks}

\begin{remarks}[The case $\lambda\not=0$]\label{rem:simpledeg1}
The affine  structure is given in terms of our chart by
\[
\nabla (\frac{dt}{t})=\lambda \frac{dt}{t}\otimes \frac{dt}{t},\quad
\nabla (\alpha_1)=\lambda\Big(\frac{dt}{t}\otimes \alpha_1 +
\alpha_1\otimes\frac{dt}{t}\Big),\quad \nabla (\alpha_0)=0 
\]
(here $\alpha_0$ resp.\ $\alpha_1$ is a translation invariant form on 
$T_0$ resp.\ $T_1$)  and so has a logarithmic singularity on $\Omega_W(\log D)$.
The residue endomorphism is semisimple with eigen values $0$ and $\lambda$,
respects $\Omega_{D,p}\subset \Omega_{W,p}(\log D)\otimes\Ocal_{D,p}$
and acts on the quotient with eigenvalue $\lambda$.
The eigen space decomposition of $\Omega_{D}$ is integrable in the
sense that it underlies the decomposition defined by the local 
isomorphism $(F_0,F_1)|D_p: D_p\to T_0\times T_1$.
In  particular, this decomposition of $D_p$ is natural;
we denote this $D_p=D_{p,0}\times D_{p,\lambda}$, where
the factors are understood as quotients of $D_p$ (the leaf
spaces of foliations).

For the same reason as  in the case $\lambda=0$, 
$D_{p,0}$ has a natural affine structure; we call it therefore
the \emph{affine quotient} of $D_p$. The elements of  $\Ocal_{D,p}$ 
that are quotients of affine functions that have order $\lambda$ at $D_p$
factor through $F_1|D_p$. So $D_{p,\lambda}$ has a natural projective structure;
we call it therefore the \emph{projective quotient} of $D_p$. 
So this makes $D_p$  look like the exceptional divisor of the blowup
of a copy of $D_{p,0}$ in some smooth space whose dimension is that of $W$.
 
Although the triple $(F_0,t,F_1)$ is not unique, there is not a great deal of choice: 
for any other system $(F'_0,t',F'_1)$,
$(F'_0, t'{}^{-\lambda}, t'{}^{-\lambda} F'_1)$ must be obtained from 
$(F_0, t^{-\lambda},t^{-\lambda} F_1)$ by an affine-linear transformation. 
If $\lambda$ is not a negative integer, then
$F'_0$ is clearly the composite of $F_0$ and an affine-linear 
isomorphism $T_0\to T'_0$. This means that the foliation defined by $F_0$
naturally extends to a morphism $W_p\to D_{p,0}$. 
A similar argument shows that if $\lambda$ is not a positive integer, the morphism
$(t,F_1)$ defines a natural quotient $W_p\to W_{p,\lambda}$. 
We call this the \emph{pure quotient} of $W_p$ since the latter is a
pure degeneration.

So if $\lambda\notin\ZZ$, then, just as in the case $\lambda=0$, we have
a natural retraction $r: W_D\to D$,  the vector field
$t\frac{\partial}{\partial t}$ is naturally defined on $W_D$ (as a dilatation
field with factor $-\lambda$) so that we have a flat normal linearization.
Furthermore, the degeneration is locally canonically the product of a pure degeneration
and an affine space and the holonomy along $D$ is a central extension the product
of a projective linear group acting on $D_{p,\lambda}$
and an affine-linear group acting on $D_{p,0}$. 

If we let $t\to 0$ in a fixed sector (on which $\log t$ is continuous), then for
$\re(\lambda)<0$ the developing map has a limit affine equivalent to 
the projection onto $D_{p,0}$ and if 
$\re (\lambda)>0$, then the  projectivization of the developing map has a limit
projectively equivalent to the projection onto  $D_{p,\lambda}$.
\end{remarks}

\begin{definition}\label{def:infisimpledeg}
Let $D$ be a smooth connected hypersurface in an analytic manifold $W$ and let
be given an affine structure on $W-D$.
We say that the affine structure on $W-D$ has 
\emph{infinitesimally simple degeneration along $D$ of logarithmic exponent 
$\lambda\in\CC$} if 
\begin{enumerate}
\item[(i)] $\nabla$ extends to $\Omega_{W}(\log D)$ with a logarithmic pole along $D$,  
\item[(ii)] the residue of this extension along $D$ preserves the subsheaf 
$\Omega_D\subset \Omega_{W}(\log D)\otimes\Ocal_D$ and its eigenvalue on 
the quotient sheaf $\Ocal_D$ is $\lambda$ and
\item[(iii)] the residue endomorphism restricted to $\Omega_D$ is
 semisimple and all of its eigenvalues are $\lambda$ or $0$. 
\end{enumerate}
\end{definition}

It is clear from the preceding that our insertion of the adjective \emph{infinitesimally}
a priori  weakens the property in question. We show that this is often
only apparently so.

\begin{proposition}\label{prop:localmodel}
Let $D$ be a smooth connected hypersurface in an analytic manifold $W$ and let
be given an affine structure on $W-D$ which along $D$ is an infinitesimally 
simple degeneration of logarithmic exponent  $\lambda\in\CC$. 
If $\lambda\notin \ZZ -\{ 0\}$, then this is true without the adjective 
\emph{infinitesimally}, so that all the properties discussed in Remarks 
\ref{rem:simpledeg0} and \ref{rem:simpledeg1}
hold; in particular, we have a flat normal linearization.

If $\lambda$ is a nonzero integer, then at any $p\in D$ there exist a 
local equation $t$ for $D_p$ in $W_p$ and a morphism
$F=(F_0,F_1): W_p\to T_0\times T_1$ to the product of an affine space $T_0$ and 
a linear space $T_1$ such that $(F_0,t,F_1)$ is a chart for $W_p$ 
and the developing map near $p$ is affine equivalent to
\begin{align*}
(F_0, t^{-n}+\log t .c^0F_0, t^{-n}F_1+\log t. C^0_1F_0)\quad & 
\text{ when $\lambda=n$ is a positive integer,}\\
(F_0+t^n\log t. C_0^1F_1, t^n,t^nF_1 )\quad & \text{ when $\lambda=-n$ 
is a negative integer.}
\end{align*}
Here $c^0: T_0\to \CC$, $C^0_1: T_0\to T_1$ and $C^1_0: T_1\to T_0$ are affine-linear maps.
\end{proposition}

\begin{corollary}\label{cor:loctrivial}
Suppose we are in the situation of Proposition \ref{prop:localmodel}. 
If $\lambda$ is not an integer $\le 0$, then $F_0$ defines a natural the affine quotient 
$W_p\to D_{p,0}$; if $\lambda$ is not an integer $\ge 0$, then $(t,F_1)$ defines
a natural pure quotient $W_p\to W_{p,\lambda}$.
If the monodromy around $D$ is 
semisimple, then the affine structure degenerates simply along $D$.
\end{corollary}
\begin{proof}
The first two assertions are clear. As for the last, 
if $\lambda$ is a positive integer $n$, then according to \ref{prop:localmodel} the
monodromy is given the unipotent transformation in $T_0\times\CC\times T_1$ with matrix
\[
\begin{pmatrix} 
1 & 0 & 0\\
2\pi\sqrt{-1}c^0 & 1 & 0\\
2\pi\sqrt{-1}C^0_1 & 0 & 1
\end{pmatrix}.
\]
This matrix is semisimple if and only if $c^0$ and $C^0_1$ are both zero, in which case
we  a simple degeneration, indeed. The proof for  the case when $\lambda$ is a 
negative integer is similar.
\end{proof}

For the proof of Proposition \ref{prop:localmodel}
we need the following well-known fact \cite{deligne}.

\begin{lemma}\label{lemma:flatsections}
Let $\Vcal$ be a holomorphic vector bundle over the germ $W_p$
endowed with a flat connection with a logarithmic pole along $D_p$. 
Then $\Vcal$  (with its flat connection) naturally decomposes naturally
according to the images of the eigenvalues of the residue map
in $\CC/\ZZ$: $\Vcal=\oplus_{\zeta\in\CC^\times}\Vcal [\zeta]$, 
where  $\Vcal [\zeta]$ has a residue endomorphism whose 
eigenvalues $\lambda$ have the property that $\exp (2\pi\sqrt{-1}\lambda)=\zeta$.

Assume now that the residue map is semisimple and
that a local equation  $t$ for $D_p$ is given.
If the residue map has a single eigenvalue $\lambda$, then there exists a  unique
$\CC$-linear section  $s: \Vcal \otimes \CC_p\to  \Vcal$ of the reduction map 
such that $t^{-\lambda}s(u)$ is a multivalued flat section and any 
multivalued flat section is thus obtained.
If the residue has two eigenvalues $\lambda$ and $\lambda+n$ with $n$ a positive 
integer, and $\Vcal \otimes \CC_p=V_\lambda\oplus V_{\lambda +n}$ is
the eigenspace decomposition, then there exist a
$\CC$-linear section  $s: \Vcal \otimes \CC_p\to  \Vcal$ of the reduction map 
and a $C\in\Hom (V_{\lambda +n},V_\lambda)$ such that the image of
\begin{align*}
u\in V_\lambda &\mapsto  t^{-\lambda}s(u); \\ 
u\in V_{\lambda+n}&\mapsto  t^{-\lambda-n}s(u)-\log t. t^{-\lambda}sC(u).
\end{align*}
spans the space of flat multivalued sections. 
\end{lemma}

We also need a Poincar\'e lemma, the proof of which is left as an exercise.

\begin{lemma}\label{lemma:poincare}
Let $\lambda\in\CC$ and $\omega\in\Omega_{W,p}(\log D)$ be such that
$t^{-\lambda}\omega$ is closed. Then $t^{-\lambda}\omega =d(t^{-\lambda}f)$
for some $f\in\Ocal_{W,p}$ unless $\lambda$ is a nonnegative integer: then 
$t^{-\lambda}\omega =d(t^{-\lambda}f)+c\log t$ for some $f\in\Ocal_{W,p}$
and  some  $c\in\CC$.
\end{lemma}

\begin{proof}[Proof of  Proposition \ref{prop:localmodel}]
The case $\lambda=0$, although somewhat special, is relatively easy; we leave
it to the reader. We therefore assume that $\lambda\not=0$. 
Choose a local equation $t$ for $D_p$.
Put $V:=\Omega_{W,p}(\log D)\otimes \CC_p$ and let $V=V_0\oplus V_\lambda$
be the eigenspace decomposition. If $\lambda\notin\ZZ$, then according to 
Lemma \ref{lemma:flatsections}
there is a section $s=s_0+s_\lambda: V_0\oplus V_\lambda\to \Omega_{W,p}(\log D)$ of the 
reduction map
such that $s_0$ resp.\ $t^{-\lambda}s_\lambda $ map to flat
sections. Any flat section is closed, because the  connection is symmetric.
Since the residue has eigenvalue $\lambda$ on the logarithmic differentials
modulo the regular differentials, $s_0$ will take its values in the regular
differentials. So by our Poincar\'e lemma \ref{lemma:poincare} both
$s_0$ and $t^{-\lambda}s_\lambda$ take values in the exact forms: there exists
a linear $\tilde s=\tilde s_0+\tilde s_\lambda: V_0\oplus V_\lambda\to \Ocal_{W,p}$
such that $d\tilde s_0 =s_0$ and $d(t^{-\lambda}\tilde s_\lambda) =
t^{-\lambda}s_\lambda $. We put $T_0:=V_0^*$ and take for
$F_0: W_p\to T_0$ the morphism defined by $s_0$. Choose $v\in V_\lambda$
not in the cotangent space $T^*_pD$ so that $V_\lambda$ splits as
the direct sum of $\CC v\oplus (T^*_pD)_\lambda $. Then $\tilde s_\lambda(v)$ is a unit
and so $t^{-\lambda}\tilde s_\lambda(v)$ is of the form $\tilde t^{-\lambda}$ for
another defining equation $\tilde t$ of $D_p$. So upon replacing  $t$ by $\tilde t$
we can assume that $\tilde s_\lambda(v)=1$. Then we take 
$T_1=(T_pD)_\lambda$, and let $F_1:W_p\to T_1$ be defined by the
set of elements in the image of $s_\lambda$ which vanish in $p$.
The proposition then follows in this case.

Suppose now that $\lambda$ is a positive integer $n$.
Then Lemma \ref{lemma:flatsections} gives us a section
$s_0+s_n : V_0\oplus V_n\to \Omega_{W,p}(\log D)$ and a 
linear map $C:V_n\to V_0$ such that the images of $s_0$
and $t^{-n}s_n- \log t .s_0C$  are flat. 
The image of $s_0$ consists of exact forms for the same reason as before
so that we can still define $\tilde s_0: V_0\to \Ocal_{W,p}$ and 
a  flat morphism $F_0: W_p\to T_0=V_0^*$. 
If $u\in V_\lambda$, then $t^{-n}s_n(u)- \log t .s_0C(u)$ is flat and hence closed. Since
$s_0C(u)=d\tilde s_0 C(u)$ we have that $t^{-n}s_n(u)+  \tilde s_0C(u)t^{-1}dt$
is also closed. Invoking our Poincar\'e lemma yields that this
must have the form $d(\tilde s_n (u)+ c(u)\log t)$ for some $\tilde s_n (u)\in \Ocal_{W,p}$
and $c(u)\in\CC$. So $\tilde s_n (u)+\log t .(c(u)-\tilde s_0C(u))$ is a multivalued
affine function. the argument is then finished as in the previous case.

The remaining case: $\lambda$ a negative integer is done similarly.
\end{proof}

We shall need to understand what happens in 
the case of a normal crossing divisor $D\subset W$ with smooth
irreducible components $D_i$ so that we have a simple degeneration along
each irreducible component. Fortunately, we do not have to deal with the most 
general case.

Suppose for a moment that we are in the simple situation where $D$ has
only two smooth irreducible components $D_1$ and $D_2$, with 
nonzero logarithmic exponents $\lambda_1$, 
$\lambda_2$. Put $S:=D_1\cap D_2$ and let $p\in S$. We have two 
residue operators acting in $\Omega_W(\log D)\otimes\CC_p$. They 
mutually commute and  respect the exact residue sequence
\[
0\to \Omega_{S,p}\to \Omega_W(\log D)\otimes\Ocal_{S,p}\to 
\Ocal_{S_p}\oplus\Ocal_{S_p}\to 0.
\]
The affine-linear functions near $p$ will have along $D_i$ 
order zero or $-\lambda_i$. The formation of the affine quotient of 
$D$ as a quotient of its ambient germ persists as a submersion 
$W_p\to (D_1)_{p,0}$ precisely when there are no affine-linear functions
which have order zero on $D_1$ and order $-\lambda_2$ on $D_2$. So we see
that we have a local  equation $t_i$ for $D_i$ and
a morphism $F=(F_0,F_1,F_2): W_p\to T_0\times T_1\times T_2$ 
to a product of which the first factor is an affine space and the other
two are linear, which makes up with
with $t_1, t_2$ a chart and has the property that the developing map is
affine-equivalent to
\[
(F_0,t_1^{-\lambda_1}(1,F_1),
t_1^{-\lambda_1}t_2^{-\lambda_2}(1,F_2)):
W_p\to T_0\times(\CC\times T_1)\times (\CC\times T_2). 
\]
Notice that the decomposition of $S_p$ defined by $F|S_p$
underlies the eigenspace decomposition defined by the two
residue operators; the factors $T_0$, $T_1$, $T_2$ correspond to
the eigenvalue pairs $(0,0)$, $(\lambda_1,0)$ and $(\lambda_1,\lambda_2)$
respectively.
 
If $\lambda_2=0$ (but $\lambda_1\not= 0$), then only a small modification 
is needed:$T_2$ is a singleton, so that
we only have a morphism $F=(F_0,F_1): W_p\to T_0\times T_1$, and
the developing map is affine-equivalent to
\[
(F_0,t_1^{-\lambda_1}(1,F_1,\log t_2)):
W_p\to T_0\times (\CC\times T_1\times \CC). 
\]
So in this case $S_p$ is decomposed into two factors.

This immediately generalizes to

\begin{proposition}\label{prop:localmodel2}
Let $W$ be an analytic manifold, $D$ a normal crossing divisor on $W$
which is the union of smooth irreducible components $D_1,\dots ,D_k$,
and $\nabla$ an affine structure on $W-D$ which is simple
of logarithmic exponent $\lambda_i$ along $D_i$.
Assume that $\lambda_i\not=0$ for $i<k$. Suppose that 
for any pair $1\le i<j\le l$ the formation of the  affine quotient 
of the generic point of  $D_j$ 
extends across the generic point of $D_i\cap D_j$. Then at $p\in \cap_i D_i$,
we have a local equation $t_i$ for $D_i$ and a morphism
\[
F=:
\begin{cases} 
(F_0,F_1,\dots ,F_k): W_p\to T_0\times T_1\times\cdots\times T_k & 
\text{ if $\lambda_k\not=0$,}\\
(F_0,F_1,\dots ,F_{k-1}): W_p\to T_0\times T_1\times\cdots\times T_{k-1}& 
\text{ if $\lambda_k=0$,}
\end{cases}
\]
to a product of an affine space $T_0$ and linear spaces $T_1,\dots ,T_k$
which together with $(t_1,\dots ,t_k)$ define a chart for $W_p$ such that
the developing map is affine equivalent to the multivalued map 
\[
\begin{cases} 
(F_0,(t_1^{-\lambda_1}\cdots t_i^{-\lambda_i}(1,F_i))_{i=1}^k) 
&\text{ if $\lambda_k\not=0$,}\\
(F_0,(t_1^{-\lambda_1}\cdots t_i^{-\lambda_i}(1,F_i))_{i=1}^{k-2},
t_1^{-\lambda_1}\cdots t_{k-1}^{-\lambda_{k-1}}(1,F_{k-1},\log t_k))
&\text{ if $\lambda_k=0$.} 
\end{cases}
\] 
\end{proposition}

\subsection{Admissible metrics}
If $M$ is a connected complex manifold with an affine structure and $p\in M$, then a flat
hermitian form on (the tangent bundle of) $M$ restricts to a hermitian form
on $T_pM$ which is invariant under the monodromy. Conversely, a monodromy
invariant hermitian form on $T_pM$ extends to flat hermitian form on $M$.
This also shows that the kernel of such a hermitian form is integrable to a foliation 
in $M$ whose local leaf space comes with an affine structure endowed with a flat
\emph{nondegenerate} hermitian form.

\begin{remark}
Consider the situation of definition \ref{def:simpledeg}, where 
$M=W-D$ and the affine structure has simple degeneration along $D$ with exponent 
$\lambda$. A flat hermitian form $h$ on $M$ must be compatible with the structure
that we have near $D$. So when $\lambda=0$, then this gives rise to
flat hermitian structure $h_D$ on $D$. When the degeneration is pure 
(so that $D$ has a projective structure), then this determines a hermitian 
form $h_D$ on $D$ which is flat for the projective connection on $D$, so that
if $h_D$ is nondegenerate, the connection on $D$ is just the Levi-Civita
connection for $h_D$. We will be mostly concerned with the case when
$h_D$ is positive definite. Of particular interest are the cases
when $h$ is positive definite (then $h_D$ is isomorphic to a 
Fubini-Study metric) and when $h$ has hyperbolic signature $(k,1)$ and is 
negative on the normal dilatation field (then $h_D$ is isometric to a complex
hyperbolic metric).

In general we have locally on $D$ a metric product of these two cases.  
\end{remark}

\begin{definition}\label{def:admissible}
Let be given an affine analytic manifold $M$ and a dilatation field $E$ on $M$.
We say that a flat hermitian form $h$ on the tangent
bundle of $M$ is \emph{admissible relative to $E$} if we are in 
one of the following three cases:
\begin{enumerate}
\item[\emph{(ell)}] $h$ is positive definite. 
\item[\emph{(par)}] $h$ is positive semidefinite with kernel spanned by $E$. 
\item[\emph{(hyp)}] $h$ has a hyperbolic signature and $h(E,E)$ is negative everywhere.
\end{enumerate}
They define on the leaf space a Fubini-Study metric,
a flat metric  and a complex hyperbolic metric respectively, to which we shall  
simply refer  as a \emph{elliptic}, \emph{parabolic}, \emph{hyperbolic structure}.
\end{definition}

\section{Linear arrangements with a Dunkl connection}

\subsection{Review of the terminology concerning linear arrangements}
\label{subsect:review}
We adhere mostly to the notation used in the book by Orlik and 
Terao \cite{orlikterao}. 

Let $(V,\Hcal)$ be a  \emph{linear arrangement}, that is, a finite 
dimensional complex vector space $V$ and a finite collection 
$\Hcal$ of (linear) hyperplanes of $V$. We shall suppose that
$\Hcal$ is nonempty so that $\dim (V)\ge 1$.
The \emph{arrangement complement}, that is, the 
complement in $V$ of the union of the members of $\Hcal$, will be 
denoted by $V^\circ$. We will also use the superscript ${}^\circ$
to denote such a complement in analogous situations (such as the case 
of a projective setting), assuming that the arrangement is understood.

The collection of hyperplane intersections in $V$ taken from  
subsets of $\Hcal$ is denoted $\Lcal(\Hcal)$ (this includes $V$ itself
as the intersection over the empty subset of $\Hcal$).
We consider it as a poset for the reverse inclusion relation: $L\le 
M$ means $L\supseteq M$. (This is in fact a lattice with \emph{join} 
$L\vee M=L\cap M$ and with \emph{meet} $L\wedge M$ the intersection of 
the $H\in\Hcal$ containing $L\cup M$.)
The  members of $\Hcal$ are the minimal elements (the \emph{atoms}) 
of $\Lcal(\Hcal)-\{ V\}$ and $\cap_{H\in\Hcal} H$ is the unique maximal element.  
For $L\in \Lcal(\Hcal)$ we denote by $\Hcal_L$ the collection of 
$H\in\Hcal$ which contain $L$. We often think of $\Hcal_L$ as 
defining a linear arrangement on $V/L$. Clearly, $\Lcal(\Hcal_{L})$ 
is the lower link of $L$ in $\Lcal(\Hcal)$, that is, the  set of  
$M\in \Lcal(\Hcal)$ 
with $M<L$. The assignment $L\mapsto\Hcal_{L}$ identifies 
$\Lcal(\Hcal)$ with a subposet of the lattice of subsets of $\Hcal$ and
we will often tacitly use that identification in our notation.

Given an $L\in \Lcal(\Hcal)$, then each $H\in\Hcal-\Hcal_L$ meets $L$ 
in a hyperplane of $L$. The 
collection of these hyperplanes of $L$ is denoted $\Hcal^L$.
We call the arrangement complement $L^\circ\subset L$  defined 
by $\Hcal^L$ an \emph{$\Hcal$-stratum}; these partition $V$.

A \emph{splitting} of $\Hcal$ is a nontrivial decomposition of
$\Hcal$ of the form $\Hcal=\Hcal_L\sqcup \Hcal_{L'}$ with $L, L'\in\Lcal(\Hcal)$ 
and $L+L'=V$. If no splitting exists, then we say that $\Hcal$ is \emph{irreducible}. 
A member $L\in \Lcal(\Hcal)$ is called  \emph{irreducible} if 
$\Hcal_L$ is. This amounts to the property that there exist 
$(\codim (L)+1)$ hyperplanes from $\Hcal_L$ such that 
$L$ is the intersection of any $\codim (L)$-tuple out of them. Or 
equivalently, that the identity component of $\aut(V/L, \Hcal_L)$ 
is the group of scalars $\CC^\times$. 
It is clear that a member of $\Hcal$ is irreducible.
We denote by $\Lcal_\irr (\Hcal)\subset \Lcal(\Hcal)$ the subposet of 
irreducible members.

Given $L\in \Lcal(\Hcal)$, then an \emph{irreducible component} of $L$
is a maximal irreducible member of $\Lcal(\Hcal_L)$.  
If $\{L_i\}_i$ are the distinct irreducible components of $L$, then
$L$ is the transversal intersection of these in the sense that the
map $V\to \oplus_i V/L_i$ is onto and has kernel $L$.

\begin{lemma}\label{lemma:dich}
Given $L,M\in \Lcal (\Hcal)$ with $M\subset L$,
denote by $M(L)\in\Lcal (\Hcal)$ the common intersection of the members of 
$\Hcal_M-\Hcal_L$. If $M\in\Lcal_\irr (\Hcal^L)$, then
$M(L)$ is the unique irreducible component of $M$ in $\Lcal (\Hcal)$ which is not an
irreducible component of $L$.
In particular, if $L\in\Lcal_\irr (\Hcal)$ and $M\in\Lcal_\irr (\Hcal^L)$,
then either $M=M(L)\in\Lcal_\irr(\Hcal)$ or  $\{ L, M(L)\}$ are the distinct irreducible 
components of  $M$ in $\Lcal (\Hcal)$. 
\end{lemma}
\begin{proof}
Left as an exercise.
\end{proof}

 \subsection{Affine structures on arrangement complements}\label{subsect:affinearc}
Let $\Hcal$ be a linear arrangement in the complex vector space $V$.
For $H\in\Hcal$, we denote by $\omega_H$  (or $\omega_H^{V}$, if a reference 
to the ambient space is appropriate)  the unique meromorphic 
differential on $V$ with divisor $-H$ and residue $1$ along  $H$. 
So  $\omega_H=\phi_{H}^{-1}d \phi_H$, where $\phi_H$ is a linear 
equation for $H$.

Suppose $\nabla$ is a torsion free flat connection on the complement $V^\circ$ 
of the union of the members of $\Hcal$. We regard it in the first place as a 
connection on the tangent bundle and then write it as 
$\nabla:=\nabla^0 -\Omega$, where $\nabla^0$ is the standard 
(translation invariant) flat connection on the tangent bundle of $V$ and 
$\Omega$ is a $\End (V)$-valued holomorphic differential on $V^\circ$:  
$\Omega\in H^0(V^\circ ,\Omega_V)\otimes_\CC \End (V)$,
the \emph{connection form} of $\nabla$.  
The associated (dual) connection on the cotangent bundle of $V^{\circ}$ 
(also denoted by $\nabla$) is characterized by the property 
that the pairing between vector fields and differentials is flat. So 
its connection form is $-\Omega^{*}$.

\begin{corollary}\label{cor:dunklchar}
Suppose that  $\nabla$ is invariant under scalar multiplication 
(as a connection on the tangent bundle of $V$) and has a 
logarithmic singularity along the generic point of 
every member of $\Hcal$. Then for every $H\in\Hcal$, $\res_H(\nabla)$ is 
a constant endomorphism $\rho_H\in \End (V)$ whose kernel contains $H$ and $\Omega$ 
has the form 
\[
\Omega :=\sum_{H\in\Hcal} \omega_H\otimes \rho_H.
\]
If $E_V$ denotes the Euler vector field on $V$, then 
the covariant derivative of $E_V$ with respect to the constant vector field
parallel to a vector $v\in V$ is the constant vector field parallel to 
$v-\sum_{H\in\Hcal}\rho_H(v)$.

If $\rho_H\not= 0$, then $\nabla$ induces on $H\in\Hcal$ a connection 
of the same type.  
\end{corollary}
\begin{proof}
The assumption that $\nabla$ is invariant under scalar multiplication means
that the coefficient forms of $\Omega$ in $H^0(V^\circ ,\Omega_V)$ are 
$\CC^\times$-invariant. This implies that these forms are $\CC$-linear 
combinations
of the logarithmic differentials $\omega_H$ and so $\Omega$ has indeed the form
$\sum_{H\in\Hcal} \omega_H\otimes \rho_H$ with $\rho_H\in \End (V)$.
Following Lemma  \ref{lemma:affineresidue}, $\rho_H$ is zero or has 
has kernel $H$. This lemma also yields the last assertion.

Finally, let $\phi_H$ be a defining linear form
for $H$ so that we can write $\omega_H=\phi_H^{-1}d\phi_H$
and $\phi_H (u)=\phi_H(u)v_H$ for some $v_H\in V$. Then
\[
\omega_H(\partial_v) \rho_H (E_V)=\frac{\phi_H (v)}{\phi_H(z)}\phi_H(z)\partial_{v_H}=
\partial_{\rho_H(v)}.
\]
 Since $\nabla^0_{\partial_v}(e)=\partial_v$, it follows
that $\nabla_{\partial_v}(E_V)=\partial_v-\sum_{H\in\Hcal}\partial_{\rho_H (v)}$.
\end{proof}

We denote by $\overline{V}$ the projective compactification of $V$ 
obtained by adding the hyperplane at infinity $\PP (V)$.

\begin{proposition}\label{prop:dunkl}
Suppose that for every $H\in\Hcal$ we are given $\rho_H\in\End (V)$ 
with kernel $H$ and let $\Omega :=\sum_{H\in\Hcal} \omega_H\otimes \rho_H$.
Then the connection on the tangent bundle of $V^\circ$ defined 
by $\nabla :=\nabla^\circ-\Omega$
is  $\CC^{\times}$-invariant and torsion free. As a connection on
the cotangent bundle it extends to
$\Omega_{\overline{V}}(\log (\PP (V))$ with logarithmic 
singularities so that $\nabla$ is regular-singular.  Moreover, the
following properties are equivalent 
\begin{enumerate}
\item[(i)] $\nabla$ is flat,
\item[(ii)] $\Omega\wedge \Omega=0$,
\item[(iii)]  for every pair  $L,M\in \Lcal(\Hcal)$ with $L\subset M$, 
$\sum_{H\in\Hcal_L}  \rho_H$ and
$\sum_{H\in\Hcal_{M}}  \rho_H$ commute,
\item[(iv)] for every $L\in \Lcal(\Hcal)$ of codimension $2$, the sum
$\sum_{H\in\Hcal_L}  \rho_H$ commutes with each of its terms.
\end{enumerate}
\end{proposition}
\begin{proof}
The $\CC^\times$-invariance of $\nabla$ is clear. Let $\phi_H\in V^*$
have zero set $H$. Then there exist $e_H\in V$ such that 
\[
\Omega=\sum_{H\in\Hcal} \phi_H^{-1} d\phi_H\otimes d\phi_H\otimes \partial_{e_H}
\]
which plainly shows that $\Omega$ is symmetric in the first two factors. So
$\nabla$ is symmetric. 
The connection $\nabla$ has on $\Omega_{\overline{V}}(\log (\PP (V))$ 
visibly a logarithmic  singularity along each member of $\Hcal$ and so 
it remains to verify that this is also the case along $\PP (V)$. 
It is clear that $\PP(V)$ is pointwise fixed under 
the $\CC^{\times}$-action. 
The generic point $w$ of $\PP(V)$ has a local defining equation $u$ in 
$\overline{V}$ that is homogeneous of degree $-1$. The
$\CC^{\times}$-invariance of $\nabla$ implies
that its matrix  has the form
\[
\frac{du}{u}\otimes A(w)+\Omega'(w),
\]
where $A$ is a matrix and  $\Omega'$ a matrix valued 
differential in the generic point of $\PP(V)$. 

The proof that the four properties are indeed equivalent can be found in
\cite{looijenga}.
\end{proof}

\begin{example}[The case of dimension two]\label{example:dimtwo}
Examples abound in dimension two: suppose $\dim V=2$ and let $\{ \rho_i\in\End 
(V)\}_{i\in I}$ a finite collection of rank one endomorphisms with
$\ker ( \rho_i)\not=\ker (\rho_j)$ if $i\not=j$ and which has more 
than one member. So if $\omega_i$ is the logarithmic differential defined by 
$\ker(\rho_i)$, then the connection defined by 
$\Omega =\sum_i \omega_i\otimes\rho_i$ is flat, precisely when 
$\sum_i \rho_i$ is a scalar operator.

Notice that in that case $I$ has just two elements $\rho_1, \rho_2$, then
both must be semisimple. This is because the centralizer of  
$\rho_i$ in $\End (V)$ is spanned by $\rho_i$ and the  identity. 
\end{example} 

\begin{example}[Complex reflection groups]\label{ex:complexreflection}
Irreducible examples in dimension $\ge 2$ can be obtained 
from finite complex reflection groups. Let 
$G\subset\GL (V)$ be a finite irreducible subgroup generated by complex 
reflections and let $\Hcal$ be the collection of fixed point hyperplanes 
of the complex reflections in $G$. Choose a $G$-invariant positive definite
inner product on $V$ and let for $H\in\Hcal$, $\pi_{H}$
be the orthogonal projection along $H$ onto $H^\perp$. 
If $\kappa\in\CC^{\Hcal}$ is $G$-invariant, then the connection 
defined by the 
form $\sum_{H\in\Hcal} \omega_{H}\otimes \kappa_{H}\pi_{H}$ is flat \cite{looijenga}.
\end{example} 

The next subsection describes a classical example.

\subsection{The Lauricella local system}\label{subsect:lauricella}
Let $V$ be the quotient of $\CC^{n+1}$ by its main diagonal.
Label  the standard basis of $\CC^{n+1}$ as $e_0,\dots ,e_{n+1}$ and let
for $0\le i<j\le n$,  $H_{ij}$ be the hyperplane $z_i=z_j$ (either 
in $\CC^{n+1}$ 
or in $V$) and $\omega_{ij}:= (z_i-z_j)^{-1}d(z_i-z_j)$ the associated 
logarithmic form.
We let $\Hcal$ be the collection of these hyperplanes so that
we can think of $V^\circ$ as the configuration space of $n+1$ distinct
points in $\CC$ given up  to translation. 

Let be given  positive real numbers $\mu_0,\dots ,\mu_n$ and define an
inner product $\la\, ,\, \ra$ on $\CC^{n+1}$ by $\la e_i,e_j\ra=\mu_i\delta_{i,j}$.  
We may identify $V$ with the orthogonal complement of the main diagonal, that
is, with the hyperplane defined by $\sum_i\mu_iz_i=0$.
The line orthogonal to the hyperplane $z_i-z_j=0$ is spanned by the
vector $\mu_je_i-\mu_ie_j$. (For this reason it is often convenient to use
the basis $(e_i':=\mu_i^{-1}e_i)_i$ instead, for then the hyperplane in 
question is the orthogonal complement of $e'_i-e'_j$; notice that 
$\la e'_i,e'_j\ra=\mu_i^{-1}\delta_{i,j}$.) So the
endomorphism $\tilde\rho_{ij}$ of $\CC^{n+1}$ which sends
$z$ to $(z_i-z_j)(\mu_je_i-\mu_ie_j)$
is selfadjoint, has $H_{ij}$ in its kernel and 
has $\mu_je_i-\mu_ie_j$ as eigenvector with eigenvalue $\mu_i+\mu_j$.
In particular, $\tilde\rho_{ij}$ induces an endomorphism $\rho_{ij}$ in $V$.

\begin{propdef}\label{def:lauricellaconn}
The connection 
\[
\nabla:=\nabla^0-\sum_{i<j} \omega_{ij}\otimes \rho_{ij}
\]
is flat (we call it the  \emph{Lauricella connection}) and 
has the Euler vector field on $V$ as a dilatation field with factor 
$1-\sum_i\mu_i$.

Let $\gamma$ be a path in $\CC$ which 
connects $z_i$ with $z_j$ but otherwise avoids
$\{ z_0,\dots ,z_n\}$ in $\CC$. If both $\mu_i<1$ and 
$\mu_j<1$ and a determination of the integrand in 
\[
\int_\gamma (z_0-\zeta )^{-\mu_0}\cdots (z_n-\zeta )^{-\mu_n} d\zeta 
\]
is chosen, then this integral converges. It is translation invariant and thus defines a multivalued 
holomorphic (so-called \emph{Lauricella}) function on $V^\circ$. This 
function is homogeneous of degree $1-\sum_i\mu_i$ and its differential  
is flat for the Lauricella connection.
\end{propdef}
\begin{proof}
The first assertion follows from a straightforward computation based on Proposition 
\ref{prop:dunkl}: one verifies that for $0\le i<j<k\le n$ the transformation
$\tilde\rho_{ij}+\tilde\rho_{ik}+\tilde\rho_{jk}$ acts on the orthogonal complement
of $e_i+e_j+e_k$ in the span of  $e_i, e_j, e_k$ as multiplication by 
$\mu_i+\mu_j+\mu_k$ so that this sum commutes with each of its terms.

The convergence and the translation invariance and the homogeneity
property of the integral are clear. If $F$ denotes 
the associated multivalued function, then the flatness of $dF$  comes down to
\[
\sum_{i,j} \frac{\partial ^2 F}{\partial z_i\partial z_j} dz_i\otimes dz_j
=-\sum_{i<j} \frac{1}{z_i-z_j}\Big( \mu_j \frac{\partial F}{\partial z_i}-
\mu_i \frac{\partial F}{\partial z_j}\Big)  (dz_i-dz_j)\otimes (dz_i-dz_j).
\]
For $i<j$, we have
\begin{multline*}
\frac{1}{z_i-z_j}\Big( \mu_j \frac{\partial F}{\partial z_i}-
\mu_i \frac{\partial F}{\partial z_j}\Big)=
\frac{-\mu_i\mu_j}{z_i-z_j}\int_\gamma 
\Big(\frac{1}{z_i-\zeta }-\frac{1}{z_j-\zeta }\Big) \prod_{\nu =0}^n
(z_\nu-\zeta )^{-\mu_\nu}d\zeta \\
= \mu_i\mu_j\int_\gamma (z_i-\zeta )^{-1}(z_j-\zeta ) ^{-1}\prod_{\nu =0}^n 
(z_\nu-\zeta )^{-\mu_\nu} d\zeta  =\frac{\partial ^2 F}{\partial z_i\partial z_j} .
\end{multline*}
If we combine this with the observation that $\sum_i \frac{\partial F}{\partial z_i}=0$,
we find the desired identity.

That $E_V$ is a dilatation field with factor $1-\sum_i\mu_i$ is left to the reader.
\end{proof}

This implies that locally, the Lauricella functions span a vector space of dimension $\le n+1$
($\le n$ in case $\sum_i\mu_i\not= 1$). We can be more precise:

\begin{proposition}\label{prop:lauricellaspan}
If $\mu_i<1$ for all $i$, then the Lauricella functions span a 
vector space of dimension $\ge n$. So if $\sum_i\mu_i\not= 1$, then  their differentials 
span the local 
system of Lauricella-flat $1$-forms.
\end{proposition}
\begin{proof}
For $i=1,\dots ,n$, we choose a path $\gamma_i$ from
$z_0$ to $z_i$ such that these paths have disjoint interior. We prove
that the corresponding Lauricella functions $F_1,\dots ,F_n$ are linearly independent.
For this it is enough to show that $F_n$ is not a linear combination of $F_1,\dots ,F_{n-1}$.
Let $T\subset \CC$  be the union of the images  of $\gamma_1,\dots ,
\gamma_n$ minus $z_n$.  We fix $z_1,\dots ,z_{n-1}$, but let 
let $z_n$ move along a path $z_n(s)$ in $\CC-T$ that eventually follows
a ray to infinity. Then $F_i(z_0,\dots ,z_{n-1},z_n (s))$ is for $s\to\infty$
approximately a constant times
$z_n(s)^{-\mu_n}$ in case $i\not= n$, and a nonzero constant times 
$z_n(s)^{1-\mu_n}$ when $i=n$. The assertion follows.
\end{proof} 

\subsection{Connections of  Dunkl type}\label{subs:hered}
The examples coming from complex reflection groups and the Lauricella 
examples suggest:

\begin{definition}\label{def:dunkl}
We say that a flat connection on $V^\circ$ whose connection form has the shape
$\Omega:=\sum_{H\in\Hcal}\omega_H\otimes \rho_H$ with $\rho_H\in\GL (V)$ 
is of \emph{Dunkl type} if there exists a positive definite inner
product on $V$ for which each $\rho_H$ is selfadjoint, in other words, 
if $\pi_H$ denotes the orthogonal projection onto $H^\perp$, then 
$\rho_H=\kappa_H\pi_H$ for some $\kappa_H\in\CC$. We call $\Omega$ a
\emph{Dunkl form} and the pair $(V,\nabla^\circ-\Omega)$ a
\emph{Dunkl system}.
\end{definition}

So in the complex reflection example we have a connection of Dunkl type and 
the same is true for the Lauricella  example. This last class  shows that it is
possible that not just the exponent function $\kappa$, but also the 
hermitian inner product (and hence the orthogonal projections $\pi_H$) 
that can deform continuously in an essential manner while retaining 
the Dunkl property. We shall see in Subsection \ref{subsect:reflarr} that for the
arrangement of type $A_n$, any connection of Dunkl type is essentially
a  Lauricella connection: its connection form is proportional to a Lauricella form.

\begin{example}\label{ex:dunkl2}
There are still many examples in dimension two. 
In order to understand the situation here,
let be given a complex vector space $V$ of dimension two 
and a finite set $\Hcal$ of lines in $V$ which comprises
at least three elements. 

Suppose that is given an inner product $\la\, ,\,\ra$ on $V$.
Choose a defining linear form $\phi_H\in V^*$ for $H$ 
of unit length relative the dual inner product and
let $e_H\in V$ be the unique vector perpendicular to $H$ on which 
$\phi_H$ takes the value $1$.  So $e_H$ is also of unit length.  
By Proposition \ref{prop:dunkl}-iv, $\kappa\in(\CC^\times)^\Hcal$
defines a Dunkl form relative to this inner product if and only if the linear map
\[
v\in V\mapsto \sum_{H\in\Hcal} \kappa_H\phi_H (v)e_H\in V
\]
commutes with each orthogonal projection $\pi_H$. This means that the map
is multiplication by a scalar $\kappa_0$. Since $\la v, e_H\ra=
\phi_H(v)$, we can also write this as
\[
\sum_{H\in\Hcal} \kappa_H\phi_H (v)\overline{\phi_H(v')} =\kappa_0\la v,v'\ra.
\]
This equality remains valid if we replace each coefficient by its real resp.\
imaginary part.  Notice, that if every $\kappa_H$ is real and positive, then 
$\kappa_H\phi_H \otimes \overline{\phi_H}$ can be thought of as 
an inner product on the line $V/H$. 

Conversely, if we are given for every $H\in\Hcal$ an inner product 
$\la \, ,\, \ra_H$ on $V/H$, and $a_H\in\RR$ is such that
 $\la\, ,\,\ra:=\sum_{H\in\Hcal} a_H\la \, ,\, \ra_H$
is an inner product on $V$, then we get a Dunkl system relative the latter with
$\kappa_H =a_H\la v,v\ra_H /\la\,v ,v\ra$ for a generator $v$ of $H^\perp$.
\end{example}

\begin{assumptions}\label{ass:positive} 
Throughout the rest of this paper we assume that $\Hcal$ is
irreducible, that the common intersection of the members of $\Hcal$ 
is reduced to $\{ 0\}$ (these are  rather innocent) and that the
residues $\rho_H$ are selfadjoint with respect to some inner 
product $\la \, ,\, \ra$ on $V$ (this is more substantial).
\end{assumptions}

Then there exist complete flags of irreducible
intersections:

\begin{lemma}\label{lemma:irrflag}
Every $L\in\Lcal_\irr (\Hcal)$ of positive dimension contains 
member of $\Lcal_\irr (\Hcal)$ of codimension one in $L$. In particular
there exists a complete flag $V>L_1>L_2>\cdots >L_n=\{ 0\}$ of irreducible
intersections from $\Hcal$. 
\end{lemma}
\begin{proof} 
If all  members of $\Hcal-\Hcal_{L}$ would contain $L^\perp$, then
$\Hcal$ would be reducible, so there exists a $H\in\Hcal-\Hcal_{L}$ 
which does contain $L^\perp$. It is clear that $L\cap H$ is 
then irreducible.
\end{proof}
  
For each linear subspace $L\subset V$
we denote by $\pi_L$ the orthogonal projection with kernel
$L$ and image $L^\perp$. So each residue $\rho_H$ is written
as $\kappa_H\pi_H$ for some $\kappa_H\in\CC$. The following
lemma shows that $\pi_L$ is independent of the inner product.

\begin{lemma}\label{lemma:uniqueform}
Suppose that none of the residues $\rho_H$ is zero.
Then any inner product on $V$ for which each of the $\rho_H$ is selfadjoint
is a positive multiple of $\la\; ,\;\ra$. (So the Dunkl form
$\Omega:=\sum_H\omega_H\otimes \kappa_H\pi_H$ then determines 
both $\Hcal$ and the inner product up to scalar.)
\end{lemma}
\begin{proof}
Suppose $\la\; ,\;\ra'$ is another hermitian
form on $V$ for which the  residues $\rho_H$ are selfadjoint. 
Then $\la\; ,\;\ra' -c\la\; ,\;\ra$
will be degenerate for some real $c\in\RR$. We prove that
this form is identically zero, in other words
that its $K\subset V$ is all of $V$. Since
$\rho_H$ is selfadjoint for this form,
we either have $K^\perp\subset H$ or $K\subset H$. So if
$\Hcal'\subset \Hcal$ resp.\ $\Hcal''\subset \Hcal$ denote 
the corresponding subsets, then for every pair 
$(H',H'')\in\Hcal'\times \Hcal''$, $H'{}^\perp\perp{H''}^\perp$.
Since $\Hcal$ is irreducible, this implies that either
$\Hcal'=\emptyset$ or $\Hcal'=\Hcal$. In the first case 
$K$ lies in the common intersection of the $H\in\Hcal$ and hence 
is reduced to $\{ 0\}$, contrary to our assumption. So we are
in the second case: $K^\perp=\{ 0\}$, that is, $K=V$. 
\end{proof}

\begin{lemma}\label{lemma:hereditary}
Let $\nabla$ be a Dunkl connection with residues $\kappa_H\pi_H$
and let $L\in \Lcal_{\irr}(\Hcal)$. 
Then the  transformation $\sum_{H\in\Hcal_L} \kappa_H \pi_H$ is of the 
form $\kappa_L\pi_L$, where 
\[
\kappa_L=\frac{1}{\codim (L)}\sum_{H\in\Hcal_L}  \kappa_H.
\]
In particular, the Euler vector field $E_V$ is a dilatation field
for $\nabla$ with factor $1-\kappa_0$.
\end{lemma}
\begin{proof}
It is clear that $\sum_{H\in\Hcal_L}  \kappa_H \pi_H$ is zero on $L$
and preserves $L^\perp$.  
Since this sum commutes with each of its terms, it will preserve 
$H$ and $H^\perp$, for each $H\in\Hcal_L$.
Since $\Hcal_L$ contains $\codim (L)+1$ members of which each 
$\codim (L)$-element subset is in general position, the induced 
transformation in $L^\perp$ will be scalar.
This scalar operator must have the same trace as $\sum_{H\in\Hcal_L} 
\kappa_H \pi_H$, and so the scalar equals the number $\kappa_L$ 
above. Since $L^{\perp}$ is the span of the lines $H^\perp$, 
$H\in \Hcal_L$, the first part of the lemma follows. The last assertion
follows from Corollary \ref{cor:dunklchar}.
\end{proof}

\begin{example}
In the Lauricella case a member $L$ of $\Lcal_{\irr}(\Hcal)$ is simply
given by a subset $I\subset\{0,\dots ,n\}$ which is not a singleton: 
it is then the set of $z\in V$ for which $z_i-z_j=0$ when 
$i,j\in I$. It is straightforward 
to verify that  $\kappa_{L}=\sum_{i\in I}\mu_i$.
\end{example}

For $\kappa\in \CC^\Hcal$, put
\[
\nabla ^\kappa:=\nabla^0-\Omega^\kappa,\quad 
\Omega^\kappa:= \sum_{H\in\Hcal} \omega_H\otimes\kappa_H\pi_H.
\]
Notice that the set of $\kappa\in (\CC^\times)^\Hcal$ for which
$\nabla^\kappa$ is flat is the intersection of a linear subspace of 
$\CC^{\Hcal}$ with $(\CC^\times)^\Hcal$. 
We shall denote that subspace by $\CC^{\Hcal, \Flat}$.

\begin{corollary}\label{cor:monotony}
Choose for every $H\in\Hcal$ a unit vector $e_H\in V$ 
spanning $H^\perp$. Then the connection $\nabla ^\kappa$
is flat if and only if for every $L\in \Lcal_{\irr}(\Hcal)$ of codimension 
two we have
\[
\sum_{H\in \Hcal_L}\kappa_H \la v,e_H\ra \la e_H,v'\ra =
\kappa_L\la \pi_L(v),\pi_L(v')\ra 
\]
for some $\kappa_L\in\CC$. In particular, $\CC^{\Hcal, \Flat}$ is 
defined over $\RR$. Moreover, any $\kappa\in (0,\infty)^{\Hcal,\Flat}$ 
is monotonic in the sense that if $L,M\in\Lcal_\irr (\Hcal)$ and $M$ 
strictly contains $L$, then  $\kappa_M<\kappa_L$.
\end{corollary}
\begin{proof}
Lemma \ref{lemma:hereditary} and condition 
(iv) of Proposition \ref{prop:dunkl} show that the flatness of 
$\nabla^\kappa$ is equivalent to the condition that
for every $L\in \Lcal_{\irr}(\Hcal)$,  
$\sum_{H\in \Hcal_L}\kappa_H\pi_H$ is proportional to $\pi_L$, in other words
that $\sum_{H\in \Hcal_L}\kappa_H\la v,e_H\ra e_H=\kappa_L\pi_H (v)$ for some
$\kappa_L\in\CC$. If we take the inner product with $v'\in V$, we see that
this comes down to the stated equality.
Since  the terms $\la v,e_H\ra \la e_H,v'\ra$ and $\la\pi_H (v),\pi_H (v')$
are hermitian, this equality still holds if we replace  
the coefficients by their complex conjugates.

Finally, if $\kappa\in (0,\infty)^{\Hcal,\Flat}$ and $L\in\Lcal_\irr(\Hcal)$ then
\[
\kappa_L \la v,v \ra=\sum_{H\in \Hcal_L}\la \kappa_H\pi_H(v),v\ra=
\sum_{H\in \Hcal_L} \kappa_H|\la v,e_H\ra |^2.
\]
If $M\in\Lcal (\Hcal)$ strictly contains $L$, 
then $\Hcal_L$ strictly contains $\Hcal_M$, and from 
\[
(\kappa_L-\kappa_M )\la v,v \ra =
\sum_{H\in \Hcal_L-\Hcal_M} \kappa_H|\la v,e_H \ra |^2
\]
it follows (upon taking $v\in L^\perp$) that $\kappa_M<\kappa_L$.
\end{proof}

Proposition \ref{prop:dunkl} shows that for every $L\in \Lcal(\Hcal)$, 
\[
\Omega_L:=\sum_{H\in\Hcal_L} \omega_H\otimes \kappa_H\pi_H
\]
defines a Dunkl-connection $\nabla_L$ in $(V/L)^\circ$. We shall see
that $L^\circ$ also inherits such a connection.

Denote by $i_L: L\subset V$ the inclusion. Notice that if $H\in\Hcal 
-\Hcal_{L}$, then $i_{L}^{*}(\omega_{H})$ is the logarithmic differential
$\omega^L_{L\cap H}$ on $L$ defined by $L\cap H$. 

The set $\Hcal^L$ of hyperplanes in $L$ injects into $\Lcal_\irr(\Hcal)$ by
sending $I$ to $I(L)$, the unique irreducible intersection such that $L \cap
I(L) = I$ as in Lemma \ref{lemma:dich}.  The set of $I \in \Hcal^L$
for which $I(L) \not \in \Hcal$ will be denoted $\Hcal^L_{\irr}$ so that
$\Hcal^L-\Hcal^L_{\irr}$ injects into  $\Hcal-\Hcal_{L}$.
We denote the image of the latter by $\Hcal^{L}_{\perp}$.

\begin{lemma}\label{lemma:condecomp}
Given $L\in\Lcal (\Hcal)$, then the connection on the tangent bundle 
of $V$ restricted to $L^\circ$ 
defined by
\[
i_{L}^{*}(\Omega -\Omega _{L})=
\sum_{H\in\Hcal -\Hcal_L} i_L^*\omega_H\otimes \kappa_H\pi_H.
\]
is flat. Moreover, the decomposition $V=L^\perp\oplus L$ defines a 
flat splitting of this bundle; on the normal bundle (corresponding to 
the first summand) the connection is given by the scalar valued 
$1$-form
$\sum_{I\in \Hcal^L_{\irr}} (\kappa_{I}-\kappa_{L})\omega_{I}^{L}$,
whereas on the tangent bundle of $L$ (corresponding to the second 
summand) it is given by the $\End (L)$-valued $1$-form
\[
\Omega^{L}:=\sum_{I\in \Hcal^L} \omega _{I}^{L}\otimes \kappa_{I(L)}\pi^L_I;
\]
here $\pi^L_I$ denotes the restriction of $\pi_I$ to $L$. We thus have a natural affine structure on 
$L^\circ$ defined by a Dunkl connection $\nabla^L$ whose form 
is defined by restriction of the inner product to $L$ and the function 
$\kappa ^L :I\in \Hcal^{L}\mapsto \kappa_{I(L)}$. The extension of that
function to $\Lcal_\irr(\Hcal^L)$ (as defined by Lemma 
\ref{lemma:hereditary}) is given by $M\in \Lcal_\irr(\Hcal^L)\mapsto\kappa_{M(L)}$.
\end{lemma}
\begin{proof}   
Let $M\in  \Lcal_{\irr }(\Hcal^{L})$. We verify that 
$\sum_{H\in\Hcal_{M}-\Hcal_{L}}\kappa_H\pi_H$ commutes with $\pi_L$
and that its restriction to $L$ equals $\kappa_{M(L)}\pi_{M(L)}$.
If $M$ is irreducible relative to $\Hcal$ (so that $M(L)=M$), then
\[
\sum_{H\in\Hcal_{M}-\Hcal_{L}} 
\kappa_H\pi_H=\kappa_{M}\pi_{M}-\kappa_L\pi_L.
\]
It is clear that the right-hand side  commutes with $\pi_L$ and 
that its restriction to $L$ is $\kappa_M\pi_M^L$.
If $M$ is  reducible relative to $\Hcal$, then $M(L)$ is the unique 
irreducible component of $M$ distinct from $L$ so that 
\[
\sum_{H\in\Hcal_{M}-\Hcal_{L}} 
\kappa_H\pi_H=\kappa_{M(L)}\pi_{M(L)}.
\]
Since $M(L)$ and $L$ are perpendicular, the right-hand side 
commutes with $\kappa_L\pi_L$  and its restriction to
 $L$ is $\kappa_{M(L)}\pi_M^L$. 
 
 The very last assertion of the proposition now follows:
by grouping the members of $\Hcal_M-\Hcal_L$ according to their
 intersection with $L$, we see that 
 \[
 \sum_{H\in\Hcal_{M}-\Hcal_{L}} 
\kappa_H\pi_H=\sum_{I\in\Hcal^L_{M}}\sum_{H\in\Hcal_I-\Hcal_L}\kappa_H\pi_H
 \]
 and according to the discussion above, the left-hand side equals
 $\kappa_{M(L)}\pi_{M(L)}$, whereas the internal sum of the 
right-hand side equals
 $\kappa_{I(L)}\pi_{I(L)}$. 
For the flatness of $\nabla^L$ we invoke criterion (iv) of 
Proposition \ref{prop:dunkl}:
if $M, N\in \Lcal_{\irr }(\Hcal^{L})$ satisfy an inclusion relation, 
then it follows from the above, that the sums $\sum_{H\in\Hcal_{M}-\Hcal_{L}} 
\kappa_H\pi_H$ and $\sum_{H\in\Hcal_{N}-\Hcal_{L}} \kappa_H\pi_H$ 
commute and the flatness follows from this. 

If we let $M$ run over the members of $\Hcal^L$ we get 
\[
\Omega -\Omega _{L}
=\sum_{I\in \Hcal^L_{\irr}} \omega _{I}^{L}\otimes \kappa_{I}\pi_{I}+ 
\sum_{H\in \Hcal^{L}_{\perp}} \omega _{H}\otimes\kappa_{H}\pi_{H}.
\]
Since all the terms commute with $\pi_{L}$ it  follows that  
$\pi_{L}$ is flat,  when viewed as an endomorphism  of the tangent bundle of $V$ 
restricted  to $L$. It also follows that the components of the connection are 
as asserted. 
\end{proof}

\begin{remark}\label{rem:kappacondition}
The last property of Lemma \ref{lemma:condecomp} imposes a very strong 
condition on $\kappa$ when viewed as a function on the poset 
$\Lcal_\irr(\Hcal)$: it implies that for any pair 
$L<M$ in this poset we have the equality  
$\kappa_M - \kappa_L = \sum_I(\kappa_I - \kappa_L)$, where the sum is taken
over all $I\in \Lcal_\irr(\Hcal)$ which satisfy $L < I \le M$ and are minimal 
for that property. In fact, it turns out that this condition yields all 
the possible weights for Coxeter arrangements of rank at least three. 
We we will not pursue this here, since 
we will obtain this classification by a different method in 
Subsection \ref{subsect:reflarr}.
\end{remark}

\begin{definition}\label{def:transvdunkl}
The Dunkl connection on $(V/L)^\circ$ resp.\ $L^\circ$ defined by 
$\Omega_L$ resp.\ $\Omega^L$ is called the \emph{$L$-transversal} resp.\
\emph{$L$-longitudinal} Dunkl connection.
\end{definition}

\subsection{Local triviality}\label{subsect:loctriv} 
Let $L\in \Lcal_{\irr}(\Hcal)$, $f:\bl_L V\to V$ be the blow-up of $L$ in $V$
and denote by $D$ the exceptional divisor. The inner product identifies
$V$ with $L\times V/L$ and this identifies
$\bl_L V$ with  $L\times\bl_0 (V/L)$, $D$ with  $L\times \PP (V/L)$ and
$\Omega_{\bl_L V}(\log D)$ with 
$pr_L^*\Omega_L\oplus pr_{V/L}^*\Omega_{\bl_0 (V/L)} (\log \PP (V/L))$.
The projection on the second factor defines  a  natural projector
$\Omega_{\bl_L V}(\log D)\to \Omega_{\bl_L V}(\log D)$, which we shall denote
by $\tilde\pi_L^*$.

\begin{lemma}\label{lemma:infsimple}
The  affine structure on $V^\circ$ is of infinitesimal simple type along $D$ with 
logarithmic exponent $\kappa_L-1$: its  residue  is $(\kappa_{L}-1)\tilde\pi^{*}_{L}$.
When $\kappa_L\not=1$, the first factor of the  product decomposition  
$D=L\times \PP (V/L)$ is the affine quotient and the second factor the projective
quotient of $D$ (in the sense of Remark \ref{rem:simpledeg1}).
\end{lemma}
\begin{proof}
The last assertion is clearly a consequence of the first.
Let $p$  be a generic point of $D$, precisely,  suppose that $p\in D$ and $p$ 
not in the strict transform of any $H\in\Hcal-\Hcal_L$.
We identify $V^*$ with $(V/L)^*\oplus (V/L^\perp)^*$. We must show
that for $y\in (V/L^\perp)^*$ and $x\in (V/L)^*$, $f^*\nabla (dy)$ and 
$f^*\Big( \nabla (x^{-1}dx)-(\kappa_L-1)x^{-1}dx\otimes x^{-1}dx\Big)$
both lie in $\Omega_{\bl_L V,p}\otimes\Omega_{\bl_L V,p}(\log D)$.
 The pull-back of $\omega_H$
to $\bl_L V$ is a regular differential at $p$ unless $H\in\Hcal_L$, in which case
it is logarithmic differential with residue one. We have that 
\[
\nabla (dy)=\sum_{H\in\Hcal} \omega_H\otimes\kappa_H\pi_H^*(dy)
\]
and since $\pi_H^*(dy)=0$ in case $H\in\Hcal_L$, we see right away that
$f^*\nabla (dy)\in \Omega_{\bl_L V,p}\otimes  \Omega_{\bl_L V,p}$.
Now consider 
\[
\nabla (x^{-1}dx)=-\frac{dx}{x}\otimes \frac{dx}{x}
+\sum_{H\in\Hcal} \omega_H\otimes\kappa_H\pi_H^*(\frac{dx}{x}).
\]
Let us first concentrate on the subsum over $\Hcal_L$. Fix a local defining
equation $t$ of $D$ at $p$. Then $(f^*x^{-1}dx)_p-t^{-1}dt$ is regular
and  so is $(f^*\omega_H)_p-t^{-1}dt$ when $H\in\Hcal_L$.
So if we calculate modulo $\Omega_{\bl_L V,p}\otimes \Omega_{\bl_L V,p}(\log D)$, 
then we find
\begin{multline*}
f^*\Big(\sum_{H\in\Hcal_L} \omega_H\otimes\kappa_H\pi_H^*(\frac{dx}{x})\Big)
\equiv \frac{dt}{t}\otimes f^*\Big(\sum_{H\in\Hcal_L} 
\kappa_H\pi_H^*(\frac{dx}{x})\Big)\equiv\\
\equiv \frac{dt}{t}\otimes f^* (\kappa_L\pi_L^*\frac{dx}{x})
\equiv \kappa_L f^*\Big( \frac{dx}{x}\otimes  \frac{dx}{x}\Big)
\end{multline*}
So it remains to show that
\[
 f^*\Big(\sum_{H\in\Hcal-\Hcal_L} 
 \omega_H\otimes\kappa_H\pi_H^*(\frac{dx}{x})\Big)\in 
 \Omega_{\bl_L V,p}\otimes \Omega_{\bl_L V,p}(\log D).
 \]
 Here all the $f^*\omega_H$ are regular at $p$, so it is rather the denominator 
 of $x^{-1}dx$ that is cause for concern.
 For this we group the $H\in \Hcal-\Hcal_L$ according to their intersection with $L$.
Let $I\in \Hcal^L$. Then for $H\in \Hcal_I-\Hcal_L$, the restriction of $\omega_H$
to $L$ as a form is $\omega_I^L$, hence independent of $H$. The same is true 
for $f^*\omega_H$: its restriction to $D$ as a form is the pull-back of 
$\omega_I^L$ and hence independent of $H$.
This means that if $H,H'\in \Hcal_I-\Hcal_L$, then the image of 
$f^*\omega_H-f^*\omega_{H'}$ in $\Omega_{\bl_L V,p}(\log D)$ can we written
as a form divisible by $dt$ plus a form divisible by $t$. In other words, it lies in
$t\Omega_{\bl_L V,p}(\log D)$. Since $f^*(x^{-1})\in t^{-1}\Ocal_{\bl_LV,p}$,
it follows that if we fix some $H_0\in \Hcal_I-\Hcal_L$, then
\[
f^*\Big(\sum_{H\in\Hcal_I-\Hcal_L} 
 \omega_H\otimes\kappa_H\pi_H^*(\frac{dx}{x})\Big)\equiv
 f^*\omega_{H_0}\otimes f^*\Big(\sum_{H\in\Hcal_I-\Hcal_L} 
 \kappa_H\pi_H^*(\frac{dx}{x})\Big).
 \]
If $I$ is irreducible, then $\sum_{H\in\Hcal_I-\Hcal_L} \kappa_H\pi_H=
\kappa_I\pi_I-\kappa_L\pi_L$. Since $\pi_I^*$ and $\pi_L^*$ leave $dx$
invariant, it follows that 
$\sum_{H\in\Hcal_I-\Hcal_L} \kappa_H\pi_H^*(x^{-1}dx)=
\sum_{H\in\Hcal_I-\Hcal_L} (\kappa_I-\kappa_L)(x^{-1}dx)$, and hence
the image of this sum under $f^*$ lies in $\Omega_{\bl_L V,p}(\log D)$.
If $I$ is reducible, then it has two irreducible components $L$ and  $I(L)$.
In that case $\sum_{H\in\Hcal_I-\Hcal_L} \kappa_H\pi_H^*(x^{-1}dx)
=\kappa_{I(L)}\pi_{I(L)}^*(x^{-1}dx)$  and since $L^\perp$ is in the  kernel of 
$\pi_{I(L)}$   it follows that the latter is identically zero.
The proof of the lemma is now complete.
\end{proof}

Given $L\in \Lcal (\Hcal)$, then we say that a Dunkl connection on $V^\circ$
has \emph{semisimple holonomy around $L$} if the holonomy around the 
exceptional divisor of the blowup $\bl_LV$ 
has that property. It is a property we know is satisfied when 
$\kappa_L\notin\ZZ$ or $\kappa_L=0$. 

\begin{corollary}\label{cor:affineproj}
Suppose we have semisimple holonomy around $L\in \Lcal_{\irr}(\Hcal)$.
Then the conditions (and hence the conclusions)
of Proposition \ref{prop:localmodel} are satisfied in the generic point of the blow-up
$f:\bl_L V\to V$ of $L$ in $V$ with $\lambda=\kappa_L -1$. In particular,
we have a normal linearization in the generic point of
the exceptional divisor of  $\bl_L V$. 
\end{corollary} 

Here is a simple application.

\begin{corollary}\label{cor:invarform}
If no $\kappa_H$ is an integer
and $\kappa_0-1$ is not a negative integer, then
every flat $1$-form on $V^\circ$ is zero. (Equivalently, every cotangent vector
of $V^\circ$ which is invariant under the monodromy representation is zero.)
Similarly, if no $\kappa_H$ is a negative integer and $\kappa_0-1$ 
is not a positive 
integer, then every flat vector field on $V^\circ$ is zero.
\end{corollary}
\begin{proof}
We only prove the first assertion; the proof of the second is similar.
Let $\alpha$ be a flat $1$-form on $V^\circ$. Since the Dunkl connection is
torsion free, $\alpha$  is closed. Let us verify that under the 
assumptions of the statement, $\alpha$ is regular
in the generic point of $H\in\Hcal$.
Near the generic point of $H$ is a linear combination of the pull-back 
of a differential on the generic point of $H$ under the canonical retraction and  a
differential which is like $\phi^{-\kappa_H}d\phi$, where $\phi$ is a local
defining equation for $H$. So if the latter appears in $\alpha$ with nonzero 
coefficient, then $\kappa_H$ must be an integer and  this we excluded. 
So $\alpha$ is regular in the generic point of $H$.

Hence $\alpha$ is regular on all of $V$. On the other hand,
$\alpha$ will be homogeneous of degree $1-\kappa_0$. So if $\alpha$ is 
nonzero, then $1-\kappa_0$ is a positive  integer.  But this we excluded 
also.
\end{proof}

Let  now $L_0>\cdots >L_k>L_{k+1}=V$ be a flag in $\Lcal_\irr (\Hcal)$ and 
let $f: W\to V$ be the iterated blowup of these subspaces
in the correct order: starting with $L_0$ and ending with $L_k$. Denote the 
exceptional divisor over $L_i$ by $E_i$, so that the $E_i$'s make up a  normal 
crossing divisor.  The common intersection $S$ of the $E_i$'s
has a product decomposition
\[
S\cong L_0\times \PP (L_1/L_0)\times\cdots\times \PP (V/L_k).
\]

\begin{proposition}\label{prop:normalformctd} 
Let $z=(z_0,\dots ,z_{k+1})$ be a general point of $S$.
If we have semisimple  holonomy around every $L_i$, then
there exist a local equation  $t_i$ for $E_i$ and a morphism 
$(F_1,\dots ,F_{k+1}): W_z\to T_1\times\cdots \times T_{k+1}$ 
to a product of linear spaces such that 
\begin{enumerate}
\item[(i)] $F_i|S_z$ factors through a local isomorphism $\PP (L_i/L_{i-1})_{z_i}\to T_i$
(and hence the system $(pr_{L_0},t_0,F_1,\dots ,t_k,F_{k+1})$ is chart for $W_z$),
\item[(ii)] the developing map at $z$ is affine equivalent to the multivalued map
$W_z\to L_0\times (\CC\times T_1)\times \cdots\times (\CC\times T_{k+1})$ 
given by
\[
\Big( pr_{L_0}, t_0^{1-\kappa_0}(1,F_1),t_0^{1-\kappa_0}t_1^{1-\kappa_1} (1,F_2)
,\dots ,t_0^{1-\kappa_0}t_1^{1-\kappa_1}\cdots t_k^{1-\kappa_k}(1,F_{k+1})\Big) ,
\]
\end{enumerate}
where $\kappa_i$ stands for $\kappa_{L_i}$.

If $\kappa_k=0$, but the holonomy around $L_i$ is semisimple for $i<k$, then 
then there exist a local equation  $t_i$ for $E_i$ and a morphism 
$(F_1,\dots ,F_{k}): W_z\to T_1\times\cdots \times T_{k}$ to a 
product of linear spaces such that 
\begin{enumerate}
\item[(i)] $F_i|S_z$ factors through a local isomorphism $\PP (L_i/L_{i-1})_{z_i}\to T_i$
if $i<k$, whereas $F_k|S_z$ factors through a local isomorphism 
$\PP (L_k/L_{k-1})\times \PP(V/L_k)_{(z_k,z_{k+1})}\to T_k$,
\item[(ii)] the developing map at $z$ is affine equivalent to the multivalued map
$W_z\to L_0\times (\CC\times T_1)\times \cdots\times (\CC\times T_k\times\CC)$
given by 
\[
\Big( pr_{L_0}, t_0^{1-\kappa_0}(1,F_1),t_0^{1-\kappa_0}t_1^{1-\kappa_1} (1,F_2)
,\dots ,t_0^{1-\kappa_0}t_1^{1-\kappa_1}\cdots t_{k-1}^{1-\kappa_{k-1}}
(1,F_k,\log t_{k})\Big).
\]
\end{enumerate}
\end{proposition}
\begin{proof}
This is a straightforward application of Proposition \ref{prop:localmodel2}. 
To see that this
applies indeed, we notice that the formation of the
affine quotient of $E_0$ is its projection to $L_0$, hence defined everywhere
on $E_0$. Likewise, the formation of the affine quotient of $E_i$ is defined
away from the union $\cup_{j<i}E_j$ of exceptional divisors of previous blowups
and given by the projection $E_i-\cup_{j<i}E_j\to L_i-L_{i-1}$. 
\end{proof}

\subsection{A classification of Dunkl forms for reflection
arrangements}\label{subsect:reflarr}
Let be given be a complex vector space $V$  in which acts
a finite complex irreducible reflection group $G\subset \GL (V)$.
We suppose that the action is essential so that $V^G=\{ 0\}$.
Let $\Hcal$ be the collection reflection hyperplanes  of $G$ in $V$.
We want to describe the space of Dunkl connections on $V^\circ$,
where we regard the inner product as unknown. So we wish
to classify the pairs $(\la\, ,\, \ra,\kappa)$, where $\la\, ,\, \ra$ 
is an inner product on $V$ and $\kappa\in \CC^\Hcal$ is such that 
$\sum_{H\in\Hcal} \omega_H\otimes\kappa_H\pi_H$ is a Dunkl form
(with $\pi_H$ being the projection with kernel $H$ that is orthogonal 
relative to $\la\, ,\, \ra$). We shall see that in case
$G$ is a Coxeter group of rank $\ge 3$, any such Dunkl system is $G$-invariant 
and hence of the type investigated in Subsection \ref{subsect:hecke}, 
unless $G$ is of type $A$ or $B$. We begin with a lemma.

\begin{lemma}\label{lemma:a2}
Let $V$ be a complex inner product space of dimension two
and let $\Hcal$ be a collection of lines in $V$. 
\begin{enumerate}
\item[(i)] If $\Hcal$ consists of two distinct elements, then a compatible
Dunkl system exists if and only if the lines are perpendicular. 
\item[(ii)] If $\Hcal$ consists of three distinct elements, then a
compatible Dunkl form exists if and only if the corresponding three points in 
$\PP(V)$ lie on a geodesic (with respect to the Fubini-Study metric).
Such a form is unique up to scalar. 
\item[(iii)] Let $(\phi_1,\phi_2)$ be a basis of $V^*$ such that $\Hcal$ consists
of the lines $H_1,H_2,H', H''$ defined by the linear forms 
$\phi_1,\phi_2, \phi':=\phi_1+\phi_2,\phi'':=\phi_1-\phi_2$. Suppose that $\la\, ,\, \ra$ is an
inner product on $V$ for which $H_1$ and $H_2$ are perpendicular. Let
$\mu_i$ be the square norm of $\phi_i$ relative to the inverse inner product on $V^*$.
Then  for every system $(\kappa_1,\kappa_2,\kappa',\kappa'')$
of exponents of a compatible Dunkl system there exist 
$a,b\in\CC$ such that $\kappa'=\kappa''=b(\mu_1+\mu_2)$
and $\kappa_i=a+2b\mu_i$ for $i=1,2$.
\end{enumerate}
\end{lemma}
\begin{proof} The proofs are simple calculations.
The first statement is easy and left to the reader. To prove the second:
let $H_1,H_2,H_3$ be the three members of $\Hcal$. Choose a defining linear form 
$\phi_i\in V^*$
for $H_i$ in such a way that $\phi_1+\phi_2+\phi_3=0$. The triple 
$(\phi_1,\phi_2,\phi_3)$ is then defined up to a common scalar factor.
Let $V(\RR)$ be the set of $v$ on which each $\phi_i$ is $\RR$-valued.
This is a real form of $V$ and the image $P$ of $V(\RR)-\{ 0\}$ in $\PP(V)$
is the unique real projective line which contains the three points defined by 
$H_i$'s. The funcions $\phi_1^2, \phi_2^2,\phi_3^2$ form a basis of
the space of quadratic forms on $V$ and so 
if $\la \, ,\, \ra$ is an inner product on $V$, then its real part restricted to
$V(\RR)$ is the restriction of  $\sum_i a_i\phi_i^2$ for unique $a_i\in\RR$.
Then $P$ is a geodesic for the associated Fubini-Study metric on $\PP(V)$
if and only if complex conjugation with respect to $V(\RR)$ interchanges the arguments
of the inner product. The latter just means that 
$\la \, ,\, \ra=\sum_i a_i\phi_i\otimes \overline{\phi_i}$.
According to Example \ref{ex:dunkl2} this is equivalent to: 
$\la\, ,\, \ra$ is part of a Dunkl system with
$\kappa_i=a_i|\phi _i(v)|^2/\la v,v\ra$, where $v$ is a generator of $H_i^\perp$
(and any other triple 
$(\kappa_1,\kappa_2,\kappa_3)$ is necessarily proportional to this one).

To prove the last statement, let $(e_1,e_2)$ be the basis of $V$ 
dual to $(\phi_1,\phi_2)$. Since $e_1\pm e_2$ has square length 
$\mu_1^{-1}+\mu_2^{-1}$, a quadruple 
$(\kappa_1,\kappa_2,\kappa',\kappa'')$ is a system of exponents 
if and only if there exist a $\lambda\in\CC$ such for all $v\in V$:
\begin{multline*}
\lambda v=\mu_1\kappa_1 \la v,e_1\ra e_1+\mu_2\kappa_2 \la v,e_2\ra e_2
+\kappa'\frac{\mu_1\mu_2}{\mu_1+\mu_2}\la v, e_1+e_2\ra (e_1+e_2)\\
+\kappa''\frac{\mu_1\mu_2}{\mu_1+\mu_2}\la v, e_1-e_2\ra (e_1-e_2).
\end{multline*}
Subsituting $e_1$ and $e_2$ for $v$ shows that this amounts to:
\begin{align*}
\kappa'=\kappa'',\quad 
\lambda=\kappa_1+\frac {\mu_2(\kappa'+\kappa'')}{\mu_1+\mu_2}=
\kappa_2+\frac {\mu_1(\kappa'+\kappa'')}{\mu_1+\mu_2}.
\end{align*}
Now put $b:=\kappa'(\mu_1+\mu_2)^{-1}=
\kappa''(\mu_1+\mu_2)^{-1}$ so that 
$\kappa_1+2b\mu_2=\kappa_2+2b\mu_1$. The assertion follows
with $a:=\kappa_1-2b\mu_1=\kappa_2-2b\mu_2$.
\end{proof}

Recall that on $A_n$, we have the Lauricella systems: 
for positive real $\mu_0,\dots ,\mu_n$ we define an
inner product $\la\, ,\, \ra$ on $\CC^{n+1}$ by $\la e_i,e_j\ra=\mu_i\delta_{i,j}$  
and  the hyperplanes $H_{i,j}=(z_i=z_j)$, $0\le i<j\le n$, restricted to the 
orthogonal complement $V=(\sum_i\mu_i z_i=0)$ of the main diagonal, 
then make up a Dunkl  system with $\kappa_{i,j}=\mu_i+\mu_j$. 
Is is convenient to switch to $\phi_i:=\mu_i z_i$ so that
$\sum_i\phi_i$ vanishes on $V$ and each $n$-element subset of is
a coordinate system. The group $G$ permutes the $\phi_i$'s (it is the full
permutation group on them) and the inner product is now $\sum_i
\mu_i^{-1}\phi_i\otimes\overline{\phi_i}$. There are
choices for the $\mu_i$'s that are not all positive for  which
$\sum_i \mu_i\phi_i\otimes\overline{\phi_i}$ is nevertheless positive definite 
on $V$. We then still have a Dunkl system and in what follows we shall 
include such cases when we refer to the term \emph{Lauricella system}.

\begin{proposition}\label{prop:dunkla}
If $G$ is of type $A_n$,  $n\ge 2$, then any Dunkl form is proportional to 
a Lauricella form.
\end{proposition}
\begin{proof}
For the case $n=2$, it easily follows from Lemma \ref{lemma:a2} that the
Lauricella systems exhaust all examples. So assume $n\ge 3$ and 
consider the space $\HH (V)$ of hermitian forms on $V$ and
regard it as a real representation of   $G=\Scal_{n+1}$. Its decomposition
into its irreducible subrepresentations has three summands: one trivial representation, 
one isomorphic to the natural real
form of $V$, and another indexed by the numerical partition $(n-1,2)$ of $n+1$.
The hermitian forms with the property that for any $A_1\times A_1$ 
subsystem the two summands are perpendicular make up a subrepresentation
of $\HH (V)$; it is in fact the sum of the trivial representation and the one isomorphic 
to $V$: these are the forms
$\sum_{i=0}^n c_i |\phi_i|^2$ with $c_i\in\RR$ restricted to the hyperplane
$\sum_{i=0}^n \phi_i=0$. The inner products in this subset are those of
Lauricella type (with $\mu_i=c_i^{-1}$). According to Lemma \ref{lemma:a2} 
such an inner product 
determines $\kappa$ on every $A_2$-subsystem up to scalar. Hence 
it  determines $\kappa$ globally up to scalar. This implies
that  the Dunkl form is proportional to one of Lauricella type.
\end{proof}

Let now $G$ be of type $B_n$ with $n\ge 3$. We use the standard set
of positive roots: in terms of the basis $e_1,\dots ,e_n$ of $\CC^n$ 
these are the basis elements themselves $e_1,\dots ,e_n$ and the 
$e_i\pm e_j$, $1\le i<j\le n$. 

\begin{proposition}\label{prop:dunklb}
Let $\mu_1,\dots ,\mu_n$ be positive real
numbers and let $a\in\CC$.
Then relative to this hyperplane system of type $B_n$ and the inner product
defined by $\la e_i,e_j\ra=\mu_i^{-1}\delta_{i,j}$, 
the exponents $\kappa_{i,\pm j}:=\mu_i+\mu_j$, $\kappa_i:=a+2\mu_i$
define a Dunkl form. In this case, $\kappa_0=a+2\sum_i\kappa_i$.
Any Dunkl form  is proportional to one of this kind for certain $\mu_1,\dots ,\mu_n;a$.
In particular, it is always invariant under reflection in the mirrors of the short roots. 
\end{proposition}
\begin{proof}
The  Dunkl property is verified for the given data by means of
Proposition \ref{prop:dunkl}-iv and the computation of $\kappa_0$ is
straightforward.

Suppose now that we are given a 
Dunkl form defined by the inner product $\la\,,\, \ra$ and the system 
$(\kappa_i,\kappa_{i,\pm j})$.
For $1\le i<j<n$ and $\varepsilon\in\{1,-1\}$ the hyperplanes 
$z_i+\varepsilon z_j=0$ and $z_n=0$ make up a
$A_1\times A_1$ system that is saturated (i.e., not contained in a 
larger system of rank two). So these hyperplanes are orthogonal. 
By letting $i$ and $j$ vary,
we find that $\la e_i,e_n\ra =0$ for all $i<n$. This generalizes to: 
$\la e_i,e_j\ra =0$ when $i\not=j$. Hence the inner product has the 
stated form. For every pair of indices $1\le i<j\le n$ we have 
a subsystem of type $B_2$ with positive roots $e_i,e_j, e_i\pm e_j$.
We can apply \ref{lemma:a2}-iii to that subsystem and find that 
there exist $a_{ij},b_{ij}\in\CC$ such that 
$\kappa_{i,j}=\kappa_{i,-j}=b_{ij}(\mu_i+\mu_j)$
and $\kappa_i=a_{ij}+2b_{ij}\mu_i$ and $\kappa_j=a_{ij}+2b_{ij}\mu_j$.
It remains to show that both $a_{ij}$ and $b_{ij}$ do not depend on 
their indices. For the $b_{ij}$'s this follows by considering a 
subsystem of type $A_2$ defined by $z_1=z_2=z_3$: our treatment of 
that case implies that we must have $b_{12}=b_{13}=b_{23}$ and this 
generalizes to arbitrary index pairs. If we denote the common value 
of the $b_{ij}$ by $b$, then we find that 
$a_{ij}=\kappa_i-2b\mu_i=\kappa_j-2b\mu_j$. This implies that 
$a_{ij}$ is also independent of its indices.
\end{proof}

\begin{corollary}\label{cor:ab}
A Dunkl system of type $B_n$ in $\CC^n$, $n\ge 3$, has $A_1^n$-symmetry 
and the quotient by this group is a Dunkl system of type $A_n$. 
If the parameters of $B_n$-system (as in
Proposition \ref{prop:dunklb}) are given by $(\mu_0,\dots ,\mu_n;a)$, then those 
of the quotient $A_n$-system are $(\mu_0,\mu_1\dots ,\mu_n)$ with $\mu_0=\half (a+1)$. 
\end{corollary}
\begin{proof}
The quotient of the Dunkl connection by the symmetry group in question
will be a flat connection on $\CC^n_u$ with logarithmic poles and 
is $\CC^\times$-invariant.
So by Corollary \ref{cor:dunklchar}, its the connection form has the shape 
$\sum_{H\in\Hcal}\omega_H\otimes\rho_H$, with $\rho_H$ a linear map.
A little computation shows that the nonzero eigenspace of 
$\rho_{(z_i-z_j=0)}$ is spanned by $e_i-e_j$ with eigen value $\mu_i+u_j$.
\end{proof}

\begin{remark}\label{rem:ab}
A $B_n$-arrangement appears in a $A_{2n}$-arrangement as the 
restriction to a linear subspace not contained in a 
$A_{2n}$-hyperplane as follows. Index the standard basis of $\CC^{2n+1}$
by the integers from $-n$ through $n$: $e_{-n},\dots , e_{n}$ and let $V$ be the
hyperplane in $\CC^{2n}$ defined by $\sum_{i=-n}^n z_i= 0$.
An arrangement $\Hcal$ of type $A_{2n}$ in $V$ is given by the 
hyperplanes in $V$ defined by $z_i=z_j$, $-n\le i<j\le n$.
The involution $\iota$ of $\CC^{2n+1}$ which interchanges $e_{-i}$ and $-e_i$
(and so sends $e_0$ to $-e_0$) leaves $V$ and the arrangement invariant; its fixed point 
subspace in $V$ is  parametrized by $\CC^n$ by: 
$(w_1,\dots ,w_n)\mapsto (-w_n,\dots ,-w_1,0,w_1,\dots ,w_n)$.
The members of $\Hcal$ meet $V^\iota$ as follows:
for $1\le i<j\le n$, $w_i=w_j$ is the trace of the 
$A_1\times A_1$-subsystem
$\{ z_i=z_j,z_{-i}=z_{-j}\}$ on $V^\iota$, likewise 
$w_i=-w_j$ is the trace
for $\{ z_i=z_{-j},z_{-i}=z_j\}$, and $w_i=0$ is the trace 
of the $A_2$-system $z_{-i}=z_i=z_0$. 
This shows that  $\Hcal |V^\iota$ is of type $B_n$. 
Suppose that we are given a Dunkl form on $V$ which is invariant 
under $\iota$.  This implies that $V^\circ $ contains $V^\circ\cap 
V^\iota$ as a flat subspace, so that the Dunkl connection on $V$ 
induces one on $V^\iota$. The values of $\kappa$ on the hyperplanes
of $V^\iota$ are easily determined: since the inner product on $V$ 
comes from an inner product on $\CC^{2n}$ in diagonal form:
$\la e_i,e_j\ra =\mu_i^{-1}\delta_{i,j}$ for certain 
positive numbers $\mu_{\pm i}$, $i=1,\dots ,n$, we must have 
$\mu_{-i}=\mu_i$. Up to scalar factor we have 
$\kappa_{(z_i=z_j)}=\mu_i+\mu_j$
for $-n\le i<j \le n$. So with that proviso,
$\kappa_{(w_i\pm w_j=0)}=\mu_i+\mu_j$, $1\le i<j\le n$
and $\kappa_{(w_i=0)}=2\mu_i+\mu_0$, which shows that we get
the Dunkl form described in Proposition \ref{prop:dunklb} 
with $a=\mu_0$. 
\end{remark}

We complete our discussion of the Coxeter case with

\begin{proposition}\label{prop:uniquedunkl}
Suppose that $G$ is a finite Coxeter group of rank $\ge 3$ 
which is not of type $A$ or $B$. Then every Dunkl system 
with the reflection hyperplanes of $G$ as its polar arrangement is 
$G$-invariant.
\end{proposition}

We shall see in Subsection \ref{subsect:hecke} that the 
local system associated to such a Dunkl system can be explicitly 
described in terms of the Hecke algebra of $G$.

We first prove:

\begin{lemma}\label{lemma:d4} 
If the complex reflection group $G$ contains a reflection subgroup 
of type $D_4$, but not one of type $B_4$, then any Dunkl form relative 
to $\Hcal$ is necessarily $G$-invariant.
\end{lemma}
\begin{proof}
We prove this with induction on the dimension of $V$. To start this off,
let us first assume that $G$ is of type $D_4$. We use the standard
root basis $(e_1-e_2, e_2-e_3, e_3-e_4, e_3+e_4)$ in \cite{bourbaki}
The four roots $\{e_1\pm e_2, e_3\pm e_4\}$ define a subsystem 
of type $(A_1)^4$. So by the first clause of Lemma \ref{lemma:a2}, these
roots are mutually perpendicular: the inner product on $V$ has the shape
\[
\la v,v\ra =a |v_1-v_2|^2+b |v_1+v_2|^2 +c |v_3-v_4|^2+d |v_3+v_4|^2
\]
for certain positive $a,b,c,d$. Any $g\in G$ sends  a $(A_1)^4$-subsystem 
to another such, and so must transform $\la \, ,\,\ra$ into a form of 
the same type (with possibly different constants $a,\dots, d$). 
From this we easily see that $a=b=c=d$, so that $\la v,v\ra =a\sum_i |v_i|^2$.
This form is $G$-invariant.  If we apply \ref{lemma:a2} to any subsystem 
of type $A_2$ and find that $\kappa$ is constant on such subsystem.
Since the $\Hcal$ is connected by its  $A_2$-subsystems, it follows
that $\kappa$ is constant.

In the general case, let $L\in\Lcal_\irr (\Hcal)$ be such that its normal system 
contains a system of type $D_4$. By our induction hypothesis,  the  Dunkl system 
transversal to $L$ is invariant under the subgroup of $g\in G$ which 
stabilizes $L$ pointwise. An inner product is already determined by its
restriction to three distinct hyperplanes; since we at least three such $L$,
it follows that the inner product is $G$-invariant. The $A_2$-connectivity
of $\Hcal$ implies that $\kappa$ is constant.
\end{proof}

\begin{proof}[Proof of Proposition \ref{prop:uniquedunkl}]
By Lemma \ref{lemma:d4} this is so when $G$ contains a subsystem of 
type $D_4$. The remaining cases are those of type $F_4$, $H_3$ and $H_4$.
In each case the essential part of the proof is to show 
that the inner product $\la\, ,\, \ra$ is $G$-invariant. Let us first do 
the case $F_4$. If we have two perpendicular roots of different length, then
they generate a saturated $A_1\times A_1$ subsystem.
So the corresponding coroots must be perpendicular for the inverse 
inner product. It is easily checked that any such an inner product must be
$G$-invariant. Lemma \ref{lemma:a2} then shows see that the exponents 
are constant on any subsystem of type $A_2$. Since a $G$-orbit 
of reflection hyperplanes is connected by its  $A_2$ subsystems, it follows
that the Dunkl form is $G$-invariant.

The cases $H_3$ and $H_4$ are dealt with in a similar fashion: any 
inner product with the property that the summands of a 
$A_1\times A_1$ subsystem (all are automatically saturated)
are orthogonal must be $G$-invariant. The $A_2$-connectivity
of the set of reflection hyperplanes implies that every such hyperplane
has the same exponent. 
\end{proof}

\section{From Dunkl to Levi-Civita}\label{subs:metric}

\subsection{The admissible range}\label{subs:admrange}
According to Lemma \ref{lemma:uniqueform}, the inner product $\la \, ,\, \ra$ is unique up to a 
scalar factor.  An inner product on $V$ determines a 
(Fubini-Study) metric on  $\PP (V)$ and two inner products determine the same
metric if and only if they are proportional. So we are then basically 
prescribing a Fubini-Study metric on  $\PP (V)$.

The inner product $\la\, ,\,\ra$ defines a translation invariant (K\"ahler) 
metric on the tangent bundle of $V$; its restriction to $V^\circ$ (which we
shall denote by $h^0$) has $\nabla^0$ as Levi-Civita connection. We shall see
that we can often deform $h^0$  with the connection. 

The main results of this subsection are

\begin{theorem}\label{thm:semipos}
Let $\dim V\ge 2$,  $\kappa\in (0,1]^{\Hcal,\Flat}$ and 
let $h$ be a hermitian form on $V^\circ$ flat for $\nabla^\kappa$ 
with at least one positive eigenvalue.
Then $h$ is positive definite if and only if $\kappa_0<1$ and for $\kappa_0=1$,
$h$ is positive semidefinite with kernel  spanned  by the Euler vector field.
 \end{theorem}
  
\begin{thmdef}\label{thm:indefinite} 
Let $\dim V\ge 2$ and  $\kappa\in (0,1]^{\Hcal,\Flat}$  be such that
$\kappa_0=1$. Assume  we are given for every $s\ge 0$ a nonzero hermitian form 
$h_s$ which is flat for $\nabla^{s\kappa}$  and such that $h_s$ 
depends real-analytically on $s$. Then there is a $m>1$  such that  for all
$s\in (1,m)$, $h_s$ is of hyperbolic signature and $h_s(E_V,E_V)$ is negative 
everywhere. The supremum $m_\hyp$ of such $m$ has the property that 
when it is finite, $h_{m_\hyp}$ is degenerate. 
We call this supremum the \emph{hyperbolic exponent} of the family.
\end{thmdef}

\begin{remark}\label{rem:erfvorm}
If $V^\circ$ has a nonzero hermitian form $h$ which is flat relative to 
$\nabla^\kappa$, and $L\in\Lcal (\Hcal)$, then such a form
is often inherited by the  transversal and longitudinal system associated to $L$.
For instance, if $L$ is irreducible and such that
$\kappa_L$ is not an integer, then the monodromy around $L$ has the two distinct 
eigenvalues $1$ and $e^{2\pi\sqrt{-1}\kappa_L}$. These decompose the tangent 
space of a point near $L^\circ$ into two eigenspaces. This decomposition is orthogonal 
relative to $h$, since the latter is preserved  by the monodromy.
Both decompositions are flat and hence are integrable to foliations.
It follows that the transversal system on $V/L$ and the longitudinal system
on $L$ inherit from $h$ a flat form. (But we cannot exclude the possibility that 
one of  these is identically zero. )
\end{remark}

The proofs of the two theorems above require some preparation. We begin with a lemma.

\begin{lemma}\label{lemma:flatsubbundle}
Let $\kappa\in (0,\infty)^{\Hcal,\Flat}$ and let $\Fcal$ be a vector subbundle
of rank $r$ of the holomorphic tangent bundle of $V^\circ$ 
which is flat for $\nabla^\kappa$. Let $\Hcal (\Fcal)$ denote the set of 
$H\in\Hcal$ for which the connection on $\Fcal$ becomes singular (relative to
its natural extension across the generic point of $H$ as a line 
subbundle of the tangent bundle). 
Then there exists an $r$-vector field $X$ on $V$ with the following properties:
\begin{enumerate}
\item[(i)] $X\vert V^\circ$ defines $\Fcal$ and the zero set of $X$ is
contained in the union of the codimension two intersections from $\Hcal$,
\item[(ii)] $X$ is homogeneous of degree
$r(\kappa_0-1)-\sum_{H\in \Hcal (\Fcal)}\kappa_H$ and multiplication 
of $X$ by  $\prod_{H\in\Hcal(\Fcal)}\phi_H^{\kappa_H}$ yields
a flat multivalued form.
\end{enumerate}
In particular, $\sum_{H\in\Hcal (\Fcal)}\kappa_H\le r\kappa_0$, so that
$\Hcal (\Fcal)\not=\Hcal$.
Moreover, in the case of a line bundle  ($r=1$), the degree of $X$
is nonnegative and is zero only when $\Fcal$ is spanned by the Euler field 
of $V$.

Likewise there exists a regular $(\dim V-r)$-form $\eta$ on $V$ satisfying similar
properties relative to the annihilator of $\Fcal$:
\begin{enumerate}
\item[(iii)] $\eta\vert V^\circ$ defines the annihilator of $\Fcal$ and the zero set of 
$\eta$ is contained in the union of the codimension two intersections from $\Hcal$,
\item[(iv)] $\eta$ is homogeneous of 
degree $(\dim V-r)(1-\kappa_0)+\sum_{H\in \Hcal-\Hcal (\Fcal)}\kappa_H$ and 
multiplication of $\eta$ by  $\prod_{H\in\Hcal-\Hcal(\Fcal)}\phi_H^{-\kappa_H}$ yields
a flat multivalued form. 
\end{enumerate}
\end{lemma} 

\begin{remark}
We will use this lemma in the first instance only in the  case of a line bundle. 
When $r=\dim V$, then clearly  $\Hcal (\Fcal)=\Hcal$ and so the lemma then tells us
that for any  translation invariant  $\dim V$-vector $X$ (i.e., one which is defined
by a generator of $\wedge^{\dim V} V$), 
$\prod_{H\in\Hcal}\phi_H^{\kappa_H} . X$ is flat for $\nabla^\kappa$.
\end{remark}

\begin{proof}[Proof of Lemma \ref{lemma:flatsubbundle}]
Let us first observe that $\Fcal$ will be invariant under scalar multiplication. 
It extends as an analytic vector subbundle of the tangent bundle over the complement
of the union of the codimension two intersections from $\Hcal$ and 
it is there given by a section $X$ of the $r$th exterior power  of the tangent bundle of $V$.
Since $\Fcal$ is invariant under scalar multiplication, we can $X$ to be homogeneous.
The local form \ref{prop:localmodel} of 
$\nabla^\kappa$ along the generic point of $H\in\Hcal $ implies that $\Fcal$  
is in this point either tangent or perpendicular to $H$. In the first case 
the connection $\nabla^\kappa$ restricted to $\Fcal$  is regular there, whereas in the 
second case it has there a logarithmic singularity with residue $-\kappa_H$. 
So if $D\pi_H$ denotes the action of $\pi_H$ on polyvectors as a derivation
(i.e., it sends an $r$-polyvector $X_1\wedge\cdots\wedge X_r$ to 
$\sum_i  X_1\wedge\cdots\wedge\pi_{H*} X_i\wedge\cdots \wedge X_r$), then
$\phi_H$ divides
$D\pi_{H }(X)$ or $D\pi_{H }(X)-X$ according to whether $H\in \Hcal-\Hcal (\Fcal)$ or 
$H\in \Hcal (\Fcal)$. Consider the multivalued function $\Phi
:= \prod_{H\in\Hcal (\Fcal)} 
\phi_H^{\kappa_H}$ on $V^\circ$. Locally we can find a holomorphic
function $f$ on $V^\circ$ 
such that $f\Phi X$ is flat for $\nabla^\kappa$; we then  have
\begin{multline*}
-\frac{df}{f}\otimes X=\nabla^0(X)-\sum_{H\in\Hcal-\Hcal (\Fcal)}
\kappa_H d\phi_H\otimes \phi_H^{-1}D\pi_{H }(X)\\
-\sum_{H\in \Hcal (\Fcal)}\kappa_H d\phi_H\otimes \phi_H^{-1}(D\pi_{H }(X)-X).
\end{multline*}
We have arranged things in such a manner that the right-hand side of this 
identity is regular. Hence so is the left-hand side. Since $X$ is nonzero 
in codimension one, it follows that $df/f$ is the restriction of a regular, globally  defined 
(closed) differential on $V$. This can only happen if $f$ is a nonzero constant.
Hence $e^{-a}\Phi X$ is a flat 
multivalued $r$-vector field on $V^\circ$. Such a field 
must be homogeneous of degree $r(\kappa_0-1)$. Since $\Phi X$ is homogeneous,
so is $e^a$. It follows that $a$ is a scalar and that the degree of $X$ is
$r(\kappa_0-1)-\sum_{H\in\Hcal (\Fcal)}\kappa_H$.
The fact that $X$ must have a degree of homogeneity 
at least $\dim V-r$ implies that 
$\sum_{H\in\Hcal (\Fcal)}\kappa_H\le r\kappa_0$. 

The assertions regarding the annihilator of $\Fcal$ are proved in a 
similar fashion. 

Now assume $r=1$ so that $X$ is a vector field.
Its degree cannot be $-1$, for then $X$ would be a constant vector field, 
that is, given by some nonzero $v\in V$. But then  
$v\in \cap_{H\in\Hcal-\Hcal (\Fcal)}H$, whereas $\Hcal (\Fcal)$ is empty or 
consists of $v^\perp$,  and this contradicts the irreducibility of $\Hcal$.

If $X$ is homogeneous of degree zero, then clearly
$\Hcal (\Fcal)=\emptyset$
(in other words, $X$ is tangent to each member of $\Hcal$)  and $\kappa_0=1$.
If we think of $X$ as a linear 
endomorphism $\Xi$ of $V$, then the tangency property amounts to 
$\Xi^*\in\End (V^*)$ 
leaving each line in $V^*$  invariant which is the annihilator of some 
$H\in\Hcal$. Since $\Hcal$ is irreducible, there are $1+\dim V$ such lines in
general position and so $\Xi^*$ is must be a scalar. This means that
$X$ is proportional to the Euler vector field of $V$. 
\end{proof}

\begin{proof}[Proof of Theorem \ref{thm:semipos}]
We first consider the case when $\dim V=2$.
Assume that $\kappa_0\le 1$. If $h$ is degenerate, then the 
kernel of $h$ is a flat line subbundle and according to Lemma
\ref{lemma:flatsubbundle} we then must  have $\kappa_0 = 1$ and this 
kernel is spanned by the  Euler vector field. For $\kappa_0 = 1$,
the Euler field is in the kernel of $h$ indeed:
if that kernel were trivial, then the orthogonal complement
of the Euler field (relative to $h^\kappa$) is also a flat subbundle $\Lcal^\perp$ of
the tangent bundle. But we have just seen that such a bundle must be generated
by the Euler field and so we have a contradiction.

Suppose now that $h\ge 0$  with kernel trivial or spanned by the Euler field.
Then $h$ induces on the punctured Riemann sphere 
$\PP (V^\circ)$ a constant curvature metric. This metric is spherical or flat
depending on whether $h>0$. 
The punctures are indexed by $\Hcal$ and at 
a puncture $p_H$, $H\in\Hcal$, the metric  has a simple type of singularity:
it is locally obtained by identifying the sides of a geodesic sector of total angle 
$2\pi (1-\kappa_H)$. The Gau\ss -Bonnet theorem (applied for instance to a geodesic 
triangulation of $\PP (V)$ whose vertices include the punctures) says that the
curvature integral $4\pi -2\pi \sum_i \kappa_i= 4\pi (1-\kappa_0)$. This implies in 
particular that $\kappa_0<1$
when $h$ is positive definite. This settles \ref{thm:semipos} in this case.

We now verify the theorem by induction on $\dim V$. So suppose $\dim V>2$.
According to Lemma \ref{lemma:irrflag}
there exists an irreducible member of $\Lcal (\Hcal)$ of dimension one.
If $\kappa_0\le 1$, then we have $\kappa_L<1$ by the monotonicity property of $\kappa$.
By Corollary \ref{cor:affineproj} we have an affine retraction of the 
germ of $L^\circ$ in $V^\circ$ and by our induction hypotheses, $h$ will be 
definite on the fibers of this retraction. It follows, that  if $h$ is degenerate,
then its kernel is of dimension one; this defines flat line subbundle
and we conclude as before that this can only happen when $\kappa_0=1$ and
the kernel is spanned by the Euler vector field.

It remains to show that if   $h$ is positive definite, then $\kappa_0<1$.
Our induction assumption implies that then $\kappa_L<1$ for all
$L\in\Lcal_\irr (\Hcal)$ different from $\{ 0\}$. Now let $H\in\Hcal$. There exists 
by Corollary \ref{cor:affineproj} an affine retraction of the 
germ of $H^\circ$ in $V^\circ$ and the restriction of $h$ to the
tangent vectors invariant under monodromy defines a form
on $H^\circ$ which is flat for the longitudinal connection.
So the Dunkl system on $H$ leaves invariant a positive definite form.
But the exponent of $\{ 0\}$ viewed as a member of $\Lcal_\irr (\Hcal^H)$
is $\kappa_0$ and so we must have $\kappa_0<1$.
\end{proof}

For the proof of Theorem \ref{thm:indefinite} we need:

\begin{lemma}\label{lemma:formchange}
Let $T$ be a finite dimensional complex vector space, $L\subset T$ a line
and $s\in (-\eps ,\eps)\mapsto H_s$ a
real-analytic family of hermitian forms on $T$ such that $H_s>0$  if
and only if $s<0$ and  $H_0\ge 0$ with kernel $L$. Then  for $s>0$, $H_s$ is 
of hyperbolic type and negative on $L$.
\end{lemma}
\begin{proof}
Let $T'\subset T$ be a supplement of $L$ in $T$. Then $H_0$ is positive definite
on $T'$. By making $\eps$ smaller, we can assume that every $H_s$ restricted to $T'$
is positive. A Gramm-Schmid process then produces
an orthonormal basis $(e_1(s),\dots, e_m(s))$ for $H_s$ restricted  
to $T'$ which depends real-analytically on $s$. Let $e\in T$ generate $L$,
so that $(e,e_1(s),\dots, e_m(s))$ is a basis for $T$.
The determinant of $H_s$ with respect this basis is easily calculated to be
$H_s(e,e)-\sum_{i=1}^m |H_s(e,e_i(s))|^2$. We know that this determinant changes 
sign at $s=0$. This can only happen if $H_s(e,e)$ is the dominating term
and (hence) changes sign at $s=0$.
\end{proof}

\begin{proof}[Proof of \ref{thm:indefinite}]
If $p\in V^\circ$, then Theorem \ref{thm:semipos} and Lemma 
\ref{lemma:formchange} applied to the restriction $h_s(p)$ of
$h_s$ to $T_pV$,  imply that there exists an $m_p>1$ such that
for  all $s\in (1,m_p)$, $h_s(p)$ is of hyperbolic signature and 
$h_s(E_V(p),E_V(p))<0$. In particular, 
$h_s$ is of hyperbolic type for $s$ in a nonempty interval $s\in (1,m')$.
We take $m$ be the supremum of the values $m'$ for which this is true
(so $h_m$ will be degenerate if $m$ is finite). This proves 
part of \ref{thm:indefinite}. The remainder amounts to the assertion
that we can take $m_p=m$ for every $p\in V^\circ$. 
For this we note that since $h_s(E_V(p),E_V(p))$ is homogeneous (of degree
$1-2s$), it suffices to verify this on the intersection of $V^\circ$
with the unit sphere $V_1$ (with respect to $\la\, ,\,\ra$).

Let us first investigate the situation near $H^\circ$, $H\in\Hcal$, for
$s$ slightly larger than $1$ (certainly such that $s\kappa_H <1$).
According to Proposition \ref{prop:localmodel}  we have a natural affine local retraction 
$r_s: V_{H^\circ}\to H^\circ$.
The naturality implies that it sends the Euler field of $V$ to the Euler field of $H$.
The naturality also accounts for the fact that $r_s$ depends real-analytically on $s$. 
The retraction $r_s$ is compatible with $h_s$ in the
sense that $h_s$ determines a hermitian form $h'_s$ on $H^\circ$
which is (i) flat for the longitudinal connection associated to $\nabla^{s\kappa}$ and 
(ii) is such that $r_s^*h_s$ and $h$ coincide on the $h_s$-orthogonal complement 
of the relative tangent space of $r_s$. In particular, $h_s$ is nonzero on the kernel of 
$dr_s$. Since $h_s$ is nondegenerate, so is $h'_s$. We know that for $s$ slightly larger
than $1$,  $h'_s$ will be of hyperbolic type. So $h_s$ must be positive on the kernel
of $dr_s$.  The Euler field $E_V$ is tangent to $H$ and we see that on $V_{H^\circ}$,
$h_s(E_V,E_V)\le (r_s^*h'_s) (E_V,E_V)=r_s^*(h'_s(E_H,E_H))$. This proves that
for every $p\in H^\circ$, there exist an $m_p>1$ and a neighborhood $U_p$ of
$p$ such that for $s\in (1, m_p)$,  $h_s(E_V,E_V)$ is negative  on $U_p\cap V^\circ$.

Now let $p\in V_1^\circ$ be arbitrary. Choose a linear subspace of dimension two 
$P\subset V$ through $p$  which is in general position with respect to $\Hcal$ in
the sense that it is not contained in a member of $\Hcal$ 
and no point of $P-\{ 0\}$ is contained in two distinct members of $\Hcal$.
Let $P_1:= P\cap V_1$ and consider the function 
\[
P_1^\circ\times (1,m)\to \RR,\quad  (p,s)\mapsto h_s(E_V,E_V)(p).
\] 
Since every point of $P_1$ is either in $V^\circ$ or in some $H^\circ$, it follows
from the  preceding discussion (and the compactness of $P_1$) 
that there exists a $m'_P\in (1,m]$
such that the above function is negative on $P^\circ\times (1,m'_P)$. Let $m_P$
be the supremum of the $m'_P$ for which this is true. 
It remains to prove that $m_P=m$. Suppose that this not the case
and assume that $m_P<m$. Then for $s=m_P$, $h_s$ is of hyperbolic type and 
$h_s(E_V,E_V)\vert P_1^\circ$ has  $0$ as maximal value.
This means that the developing map for $\nabla^{s\kappa}$ is affine-equivalent
to a morphism from a cover of  $P_1^\circ$ to the subset of $\CC^n$ defined 
by $|z_1|^2+\cdots |z_{n-1}|^2-|z_n|^2\le 0$, and  such that the inequality 
is an equality at some point. This, however, contradicts a convexity property of this subset
as is shown by the following lemma.
\end{proof}

\begin{lemma}
Let  $f=(f_1,\dots ,f_n): U\to \CC^n$ be a holomorphic map from
a connected complex manifold $U$ such that 
$|f_1|^2+\cdots +|f_{n-1}|^2\le |f_n|^2$. Then the latter inequality is strict
unless $f$ maps to a line.
\end{lemma}
\begin{proof} We may assume that  $f_n$ is not  constant equal to zero so that
each $g_i:= f_i/f_n$ is a meromorphic function. Since $g:=(g_2,\dots ,g_n)$
takes values in the closed unit ball, it is holomorphic. It is well-known that such a
map takes values in the open unit ball unless it is constant. This yields the lemma.
\end{proof}

\subsection{The Lauricella integrand as a rank two example}\label{subsect:thurston}
We do not know whether a Dunkl system with real exponents always
admits a nontrivial flat hermitian form, not even in the case $\dim V=2$.
However, if $\dim V=2$ and $\kappa_0=1$, then there is natural choice.
In order to avoid conflicting notation, let us write $P$ instead of $V$,
let $H_0,\dots ,H_{n+1}$ be the distinct elements of $\Hcal$ (so that
$|\Hcal |=n+2$) and write $\mu_i$ for $\kappa_{H_i}$ (so that $\sum_i \mu_i=2$).
Recall from Lemma \ref{lemma:flatsubbundle} that if $\alpha$ is a translation invariant
$2$-form, then $(\prod_{H\in\Hcal} \phi_H^{-\kappa_H})\alpha$ is a flat 
multivalued $2$-form. Since $\kappa_0=1$, the Euler field $E_P$ is flat,
and so if $\omega$ denotes the $1$-form obtained by
taking the inner product of $E_P$ with $\alpha$, then 
$(\prod_{i=0}^{n+1} \phi_i^{-\mu_i})\omega$ is a flat multivalued $1$-form.
Hence its absolute value,
\[
h:=|\phi_0|^{-2\mu_0}\cdots |\phi_{n+1}|^{-2\mu_{n+1}} |\omega |^2,
\]
is then a nontrivial flat hermitian form. It is positive semidefinite with
kernel spanned by the Euler field. 

This is intimately connected with an observation due to Thurston
\cite{thurston}, about which we will have more to say later on. 
Since $\kappa_0=1$,
the punctured Riemann sphere $\PP (P^\circ)$ acquires an affine structure.
The form $h$ is a pull-back from $\PP (P^\circ)$  so that $\PP (P^\circ)$
has in fact a Euclidean (parabolic) structure. If we
assume that $\mu_i\in (0,1)$ for all $i$, then $\PP(P)$
is a euclidean \emph{cone manifold} in Thurston's sense:
at the point $p_i\in \PP(P)$ defined by $H_i$, the metric is conical with 
total angle $2\pi(1-\mu_i)$. In such a point is concentrated  
a certain amount of curvature, its \emph{apex curvature} $2\pi\mu_i$,
which is its contribution to the Gau\ss-Bonnet formula (the sum of these
is indeed $4\pi$, the area of the unit sphere). On the other hand,
the multivalued form $(\prod_{H\in\Hcal} \phi_H^{-\kappa_H})\omega$
is directly related to the Lauricella integrand. To see this, choose an 
affine coordinate $z$ on $\PP(V)$ such that if $z_i:=z(p_{H_i})$, then $z_{n+1}=\infty$.
Then $(\prod_{i=0}^{n+1}  \phi_i^{-\mu_i})\omega$ is up to a constant factor 
the pull-back of a constant times  $\prod_{i=0}^n (z_i-\zeta )^{-\mu_i}dz$, which we recognize
as the Lauricella integrand. 

Of course, the $(n+1)$-tuple $(z_0,\dots ,z_n)\in \CC^{n+1}$ is defined
only up to an affine-linear transformation of $\CC$. This means that if
$V$ is the quotient of $\CC^{n+1}$ by its main diagonal (as in 
Subsection \ref{subsect:lauricella}), then only the image of $(z_0,\dots ,z_n)$ in
$\PP(V^\circ)$ matters.  Thus $\PP(V^\circ)$ can be understood as the moduli 
space of Euclidean metrics on the sphere with $n+2$ conical singularities 
which are indexed by $0,\dots ,n+1$ with prescribed apex curvature $2\pi\mu_i$
at the $i$th point.

\subsection{Flat hermitian forms for reflection arrangements}
The following theorem produces plenty of interesting situations
to which the results of Subsection \ref{subs:admrange} apply. It may
very well hold in a much greater generality.

\begin{theorem}\label{thm:herm}
Suppose that $\Hcal$ is the reflection arrangement 
of a finite complex reflection group $G$.
Then there exists a map from $(\RR^{\Hcal})^G$ to the
space of nonzero hermitian forms on the tangent bundle of $V^\circ$
(denoted $\kappa\mapsto h^\kappa$) 
with the following properties: for every $\kappa\in (\RR^{\Hcal})^G$,
\begin{enumerate}
\item[(i)] $h^\kappa$ is flat for $\nabla^{\kappa}$ and invariant under $G$.
\item[(ii)] $t\in\RR\mapsto h^{t\kappa}$ is smooth 
(notice that $h^0$ was already defined)
and the associated curve of projectivized forms,
$t\mapsto [h^{t\kappa}]$ is real-analytic.
\end{enumerate} 
Moreover this map is unique up to multiplication by a (not necessarily 
continuous) function $(\RR^{\Hcal})^G\to (0,\infty)$.

Likewise there is a map from $(\RR^{\Hcal})^G$ to the
space of nonzero hermitian forms on the cotangent bundle of $V^\circ$
(denoted $\kappa\mapsto \check{h}^\kappa$) with analogous properties.
\end{theorem}

\begin{example}\label{ex:dimone} 
For $V=\CC$ and $\Omega =\kappa z^{-1}dz$, we can take  
$h^\kappa(z):=|z|^{-2\kappa }|dz|^2$. Notice that we can expand
this in powers of $\kappa$ as 
\[
h^\kappa(z)=\sum_{k=0}^\infty \kappa^k\frac{(-\log |z|^2)^k}{k!}|dz|^2.
\]
\end{example}

We shall first prove that in the situation of Theorem  \ref{thm:herm}
we can find such an $h^\kappa$ formally at $\kappa=0$.
For this we need the following notion, suggested by Example \ref{ex:dimone}.
Let be given a complex manifold $M$ and a smooth hypersurface $D\subset M$. 
We have the real-oriented blowup of $D$ in $M$; this is a real-analytic 
manifold with boundary. If $(\phi,z_1,\dots ,z_n)$ is a 
coordinate system at $p\in D$ such
that $D$ is given by $\phi =0$, then $r:=|\phi|,\theta:= arg (\phi), 
x_i:=\re (z_i), y_i:= \im (z_i)$ are coordinates for this blowup,
where of course $\theta$ is given modulo $2\pi$ and the boundary is given by
$r=0$. We say that a function on a neighborhood of $p$ in $M-D$
is \emph{mildly singular along} $D$ if it can be written as a 
polynomial in $\log r$ with certain continuous coefficients:
we want these coefficients to be real-analytic on the
real-oriented blowup of $D$ at $p$ (and so constant on its boundary). 
Since $\phi $ is unique up to a unit factor, $\log r$ is unique up to 
an analytic function in the coordinates, and so
this notion is independent of the coordinate system.

Likewise, we say that a differential on a neighborhood of $p$
in $M-D$ is \emph{mildly singular along} $D$ if it is a linear combination
by mildly singular functions at $p$ of real-analytic forms
on the real-oriented blowup whose restriction to the boundary \emph{as a 
form} is zero. So this is a module over the ring of mildly singular functions
at $p$ and as such generated by $dr, rd\theta$ and $dx_i, dy_i$, 
$i=1,\dots ,n$.

\begin{lemma}\label{lemma:mild}
In this situation we have:
\begin{enumerate}
\item[(i)] $\log r$ is algebraically independent over the ring 
of real-analytic functions on the  real-oriented blowup of $D$ over $p$.
\item[(ii)] Any mildly singular differential at $p$ that is closed is 
the differential of a mildly singular function at $p$. 
\end{enumerate}
\end{lemma}
\begin{proof}
For the proof of (i), suppose that we have a nontrivial
relation: 
\[
\sum_{k=0}^N f_k(r,\theta, x,y)
(\log r)^k=0,
\]
with each $f_k$ analytic (and periodic in $\theta$).  
Divide then by the highest power of $r$ which divides each $f_k$, so that
now not all $f_k(0,\theta ,x,y)$ vanish identically. 
If we substitute $r:=e^{-1/\rho}$, with $\rho$ small, then 
$\sum_{k=1}^N f_k(0,\theta ,x,y)(-\rho)^{-k}$ will 
be a flat function at $\rho =0$. This can only be the case if each 
$f_k(0,\theta ,x,y)$ is identically zero, which contradicts our assumption.

For the proof of (ii) we note that if $f$ is mildly singular at
$p$, and $\eta$ is one of the module generators $dr, rd\theta, dx_i, dy_i$,
then the integral of $f\eta$ over the circle $r=\varepsilon$, 
$x=y=0$ tends to zero with $\varepsilon$.
So if $\omega$ is a closed differential that is mildly singular at $p$, then
it can be integrated to a function $f$ on the complement of $D$ in 
a neighborhood of $p$. This function will there be real-analytic.
It is a straightforward to verify that $f$ is mildly singular at $p$. 
\end{proof}

\begin{lemma}\label{lemma:formaldef}
In the situation of Theorem  \ref{thm:herm}, let $\kappa\in (\RR^\Hcal)^G$.
Then there  exists a formal expansion $h^{s\kappa}=\sum_{k=0}^\infty s^kh_k$ 
in $G$-invariant hermitian forms that
are mildly singular along the smooth part of the
arrangement with initial coefficient $h_0=h^0$,
and with the property that $h^{s\kappa}$  is flat for $\nabla^{s\kappa}$.  
\end{lemma}
\begin{proof}
The flatness of  $h^{s\kappa}$ means that for every pair $v,v'\in V$ 
(thought of as translation invariant vector fields on $V$) we have 
\[
d(h^{s\kappa}(v,v'))=-s h^{s\kappa}(\Omega^\kappa (v),v')
-s h^{s\kappa}(v, \Omega^{\kappa} (v')), 
\]
where $\Omega^\kappa (v)=\sum_H \kappa_H\pi_H(v)\otimes\omega_H$, 
which boils down to 
\[
(*)\quad d(h_{k+1} (v,v'))=  -h_{k}(\Omega^\kappa  (v),v')
-h_{k}(v, \Omega^{\kappa}  (v')),\quad k=0,1,2,\dots 
\]
In other words, we must show that we can solve
(*) inductively by $G$-invariant forms. In case we can solve (*), then it is clear
that a solution will be unique up to a constant. 

The first step is easy: if we choose our defining equation $\phi_H\in V^*$ for $H$
to be such that $\la \phi_H,\phi_H\ra=1$, then
\[
h_1(v,v'):= 
-\kappa_H\sum_H \la \pi_H(v),\pi_H (v')\ra \log |\phi_H |^2.
\]
will do.
Suppose that for some $k\ge 1$ the forms $h_0,\dots ,h_k$ have 
been constructed. In order that $(*)$
has a solution for $h_{k+1}$ we want the right-hand side (which we shall denote
by $\eta_k (v,v')$) to be exact. 
It is certainly closed: if we agree that 
$h(\omega\otimes v,\omega'\otimes v')$ stands for 
$h(v,v')\omega\wedge\overline{\omega'}$, then
\begin{multline*}
d\eta_k(v,v')
=h_{k-1}(\Omega^\kappa\wedge\Omega^\kappa (v),v')-h_{k-1}(\Omega^\kappa (v),
\Omega^{\kappa}(v'))+\\
+h_{k-1}(\Omega^\kappa (v),
\Omega^{\kappa} (v'))+h_{k-1}(v,\Omega^{\kappa} 
\wedge\Omega^{\kappa} (v'))=\\
=h_{k-1}(\Omega^\kappa\wedge\Omega^{\kappa} (v),v')
+h_{k-1}(v,\Omega^{\kappa}\wedge\Omega^{\kappa} (v'))=0
\end{multline*}
(since $\Omega^\kappa\wedge \Omega^\kappa =0$).
So in order to complete the induction step, it suffices by Lemma 
\ref{lemma:mild} that to prove that 
$\eta_k$ is mildy singular along the arrangement: since the complement in $V$  
of the singular part of the arrangement  is simply connected, we then write
$\eta_k$ as the differential of a hermitian form $h_{k+1}$ on $V$ that is mildly
singular along the arrangement and averaging such $h_{k+1}$ over its $G$-transforms
makes it $G$-invariant as well.

Our induction assumption says that near $H^\circ$ we can expand $h_k$ 
in $\log|\phi |$ as:
\[
 h_k=\sum_{i=0}^N (\log |\phi_H|)^ih_{k,i}
\]
with $h_{k,i}$ a continous hermitian form on $TV$ near $H^\circ$ which
becomes real-analytic on the on the real-oriented blowup of $H^\circ$.
We claim that the projection $\pi_H$ restricted to $TV|H^\circ$ is selfadjoint 
relative to each term $h_{k,i}$.
For $h_k$ is $G$-invariant and hence invariant  under a nontrivial complex 
reflection $g\in G$ with mirror $H$. Since $|\phi_H|$ is also invariant under $g$
and since the above expansion is unique by Lemma \ref{lemma:mild}-i, it follows 
that this property is inherited by each term $h_{k,i}$. In particular, the 
restriction of $h_{k,i}$ to $TV|H^\circ$ is invariant under $g$. Since $\pi_H$ \
is the projection on an eigenspace of $g$, the claim follows. 
Now $\eta_k$ is near $H^\circ$ modulo a mildly singular form equal to 
\[
-\kappa_H\sum_{i=0}^N 
\Big( \omega_H h_{k,i}(\pi_H(v),v')+\overline{\omega_H}
h_{k,i}(v, \pi_H(v')\Big)(\log |\phi_H|)^k.
\]
The selfadjointness property of $\pi_H$ implies that this, in turn, is modulo 
a mildly singular form equal to
\[
-2\kappa_H\sum_{i=0}^N h_{k,i}(\pi_H(v),v')(\log |\phi_H|)^k d(\log |\phi_H|),
\]
showing that $\eta_k$ is mildly singular along $H^\circ$ as desired.
\end{proof}

In order to prove Theorem \ref{thm:herm}, we begin with a few generalities 
regarding conjugate 
complex structures. Denote by $V^\dagger$ the complex vector space $V$ with 
its conjugate complex structure: scalar multiplication by $\lambda\in\CC$ 
acts on $V^\dagger$ as scalar multiplication by 
$\overline{\lambda}\in\CC$ in $V$. Then $V\oplus V^\dagger$ has a natural
real structure for which complex conjugation is simply interchanging
arguments. The ensuing conjugation on  
$\GL (V\oplus V^\dagger)$ is, when restricted to 
$\GL (V)\times \GL(V^\dagger)$, also interchanging arguments, whereas 
on the space of bilinear forms on $V\times V^\dagger$, it is given
by $h^\dagger(v,v'):=\overline{h(v',v)}$. So a real point of 
$(V\otimes V^\dagger)^*$
is just a hermitian form on $V$.

Fix a base point $*\in V^\circ$ and identify $T_*V^\circ$ with $V$. For 
$\kappa\in (\CC^{\Hcal})^G$, we denote the monodromy representation of 
$\nabla^\kappa$
by $\rho^\kappa\in \Hom (\pi_1(V^\circ,*), \GL (V))$. Notice that
$\rho^\kappa$ depends complex-analytically on $\kappa$. Then the  same 
property must hold for 
\[
\kappa\in(\CC^{\Hcal})^G\mapsto (\rho^{\overline 
{\kappa}})^\dagger\in \Hom (\pi_1(V^\circ,*), \GL (V^\dagger)).
\]

Recall from \ref{cor:monotony} that $(\CC^{\Hcal})^G$ is invariant 
under complex
conjugation. 

\begin{lemma}\label{lemma:hermalg}
Let $\HH$ be the set of pairs
$(\kappa ,[h])\in (\CC^{\Hcal})^G\times \PP((V\otimes V^\dagger)^*)$,
where  $h\in V\times V^\dagger\to \CC$ is invariant under
$\rho^\kappa \otimes  (\rho^{\overline {\kappa}})^\dagger$ and let 
$p_1:\HH\to (\CC^{\Hcal})^G$ be the projection. Then $\HH$ resp.\ 
$p_1(\HH)$ is a complex-analytic set defined over $\RR$  (in 
$(\CC^{\Hcal})^G\times \PP((V\otimes V^\dagger)^*)$ resp.\ 
$(\CC^{\Hcal})^G$) and  we have $p_1(\HH(\RR))=(\RR^{\Hcal})^G$.
\end{lemma}
\begin{proof} That  $\HH$ is  complex-analytic and defined over 
$\RR$ is clear. Since $p_1$ is proper and defined over $\RR$, 
$p_1(\HH)$ is also complex-analytic and defined over $\RR$.
If $\kappa\in(\RR^{\Hcal})^G$ is in the image of $\HH$, then 
there exists a nonzero bilinear map $h:V\times V^\dagger\to\CC$ invariant 
under 
$\rho^\kappa \otimes  (\rho^{\overline {\kappa}})^\dagger$. But then both 
the `real part' 
$\frac{1}{2} (h+h^\dagger)$
and the `imaginary part' $\frac{1}{2\sqrt{-1}} (h-h^\dagger)$ of
$h$ are hermitian forms invariant under $\rho^\kappa$
and clearly one of them will be nonzero. The lemma follows.
\end{proof}

\begin{proof}[Proof of Theorem \ref{thm:herm}]
Now let $L\subset (\CC^{\Hcal})^G$ be a line defined over $\RR$.
By the preceding discussion, there is a unique  irreducible component 
$\tilde L$ of 
the preimage of $L$ in $\HH$ which contains $(0,[h^0])$. The
map $\tilde L\to L$ is proper and the preimage of $0$ is a singleton.
Hence $\tilde L\to L$ is an analytic isomorphism. Since $L$ is defined
over $\RR$, so are $\tilde L$ and  the isomorphism $\tilde L\to L$. 
The forms parametrized  by $\tilde L (\RR)$ define a real line bundle
over $L(\RR)$.  Such a line bundle is trivial in the smooth category and
hence admits a smooth generating section with prescribed value in $0$.
We thus find a map $\kappa\mapsto h^\kappa$ with the stated properties.
The proof for the map $\kappa\mapsto \check{h}^\kappa$ is similar.
\end{proof}

If $h$ is a nondegenerate hermitian form on the tangent bundle of 
$V^{\circ}$ which is flat for the Dunkl connection, then $\nabla$ must be 
its Levi-Civita connection of $h$ (for $\nabla$ is torsion free); in particular, 
$h$ determines $\nabla$. 
Notice that to give a flat hermitian form $h$ amounts to giving a monodromy 
invariant hermitian form on the translation space of $A$. So $h$ will be 
homogeneous in the sense that the pull-back of $h$ under scalar  multiplication 
on $V^\circ$ by $\lambda\in\CC^\times$ is 
$|\lambda |^{2-2\re (\kappa_0)} h$.

\subsection{The hyperbolic exponent of a complex reflection group}
\label{subsect:hypexp}
In case $\Hcal$ is a complex reflection arrangement of a finite reflection group 
$G$, we can estimate the hyperbolic exponent. According 
to Chevalley, the graded algebra of $G$-invariants $\CC[V]^G$ is a polynomial algebra. 
Choose  a set of homogeneous generators, $f_1,\dots ,f_n$,  ordered
by their degrees $\deg (f_1)\le \cdots \le \deg (f_n)$. Although the generators are
not unique, their degrees are. We put $d_i:=\deg (f_i)$.
The number $m_i:=d_i-1$, which is the degree of the coefficients of $df_i$ on a basis of
constant differentials on $V$, is called the $i$th \emph{exponent} of $G$.
It is known that the  subalgebra of $G$-invariants in the exterior 
algebra $\CC[V]\otimes\wedge^\pt V^*$ of regular forms on $V$ is 
generated as such by $df_1,\dots ,df_n$ \cite{solomon}. In particular
any invariant $n$-form is proportional to $df_1\wedge\cdots \wedge df_n$. 

The geometric content of Chevalley's theorem is the assertion that the orbit space 
$G\bs V$ is an affine space, a fact which never stops to surprise us. 
The union of the members
of $\Hcal$ is also the union of the irregular orbits and hence is the
singular locus of the orbit map $\pi: V\to G\bs V$. The image of this
orbit map is a hypersurface in $G\bs V$, the \emph{discriminant} of $G$.
It is defined by a suitable power of the jacobian of $(f_1,\dots ,f_n)$.
 
A vector field on  $G\bs V$ lifts to $V$ precisely when it is tangent to the discriminant
and in this manner we get all the $G$-invariant vector fields on $V$.
The $G$-invariant regular vector fields  make up a graded
$\CC[V]^G$ -module and it is known \cite{orlikterao} that this module is free. 
As with the Chevalley generators, we choose a 
system of homogeneous  generators $X_1,\dots ,X_n$ 
ordered by their degree: $\deg (X_1)\le \cdots \le \deg (X_n)$.  We put 
$d^*_i:=\deg (X_i)$ and $m ^*_i:=1+\deg (X_i)$ (so that
$m^*_i$ is the degree of the coefficients of $X_i$ on a basis of
constant vector fields on $V$).   The generator of smallest degree is proportional
to the Euler field. Hence $d^*_1=0$ and $m^*_1=1$.
The number $m_i^*$ is called the $i$th \emph{co-exponent} of $G$.  It usually 
differs from $m_i$, but when $G$ is a Coxeter group they are equal,
because the defining representation of $G$ is self-dual.

A polyvector field on $G\bs V$ lifts to $V$ if and only if it does so in codimension
one (that is, in the generic points of the discriminant) and we thus obtain all
the  $G$-invariant polyvector fields  on $V$.  For reasons similar 
to the case of forms, the subalgebra of $G$-invariants in the exterior 
algebra $\CC[V]\otimes\wedge^\pt V$ of regular polyvector fields on $V$ is 
generated as such by $X_1,\dots ,X_n$.

\begin{theorem}\label{thm:hyperbolicexponent}
Suppose that $\Hcal$ is the reflection arrangement 
of a finite complex reflection group $G$ which is transitive on $\Hcal$.
Then the hyperbolic 
exponent for the ray $((0,\infty)^\Hcal)^G$ (which is defined
in view of Theorem \ref{thm:herm}) is $\ge m^*_2$.
\end{theorem}
\begin{proof}
Let $\kappa\in ((0,1)^\Hcal)^G$ be such that $\kappa_0=1$ and let 
$h_s$ be the family of hermitian forms on the tangent bundle of $V^\circ$
whose existence is asserted by Theorem \ref{thm:herm}. Let $m\in (1,\infty ]$
be its hyperbolic exponent.  If $m=\infty$ there is nothing to show, so 
let us assume that $m<\infty$. This means that $h_m$ is degenerate.
So its kernel defines a nontrivial subbundle $\Fcal$ of the tangent bundle of 
$V^\circ$ (of rank $r$, say) which is flat for $\nabla^{m\kappa}$. 
This bundle is $G$-invariant.
So the developing map maps to a vector space $A$ endowed with
a monodromy invariant hermitian form $H_m$ with a kernel of dimension
$r$. Since $H_m$ is nontrivial, so is $H_m(E_A,E_A)$
and hence so is $h_m(E_V,E_V)$. In other words, $\Fcal$ does not contain the 
Euler field.

Let $X$ be the associated $r$-vector field 
on $V$ as in Lemma \ref{lemma:flatsubbundle}. That lemma asserts that
$\Hcal(\Fcal)\not=\Hcal$. Since $\Hcal (\Fcal)$ is $G$-invariant, this 
implies that $\Hcal (\Fcal)=\emptyset$ so that $X$ has degree $r(m-1)$. 
We prove that $X$ is $G$-invariant.
Since $X$ is unique up to a constant factor 
it will transform under $G$ by means of a character. 
For this it is enough to show that $X$ is left invariant 
under any complex reflection. Let $H\in\Hcal$. The splitting $V=H\oplus H^\perp$ 
defines one of $\wedge^rV$: $\wedge^rV=\wedge^rH \oplus 
(H^\perp\otimes\wedge^{r-1}H)$. This splitting is the eigenspace decomposition 
for the action of the cyclic group $G_H$ of $g\in G$ which leave $H$ pointwise fixed.
It is clear from the way $X$ and $\Hcal (\Fcal)$ are defined that the value of $X$ on $H$
will be a section of the first summand so that $X$ is invariant under $G_H$ indeed.
Now write $X$ out in terms of our generators:
\[
X=\sum_{1\le i_1<\dots <i_r\le n} a_{i_1,\dots ,i_r}
X_{i_1}\wedge\cdots \wedge X_{i_r}, \quad a_{i_1,\dots ,i_r}\in \CC [V]^G.
\]
Since $\Fcal$ does not contain the Euler field,  $X$ is not divisible by $X_1$ and so a term
with $i_1\ge 2$ appears with nonzero coefficient.  This means that the
degree of $X$ will be at least $d_2^*+\cdots +d_{r+1}^*\ge 
r(d_2^*)=r(m_2^* -1)$.  It follows that $m\ge m_2^*$, as asserted.
\end{proof}

\begin{remark}
There are only two primitive complex reflection groups of rank greater than two 
which the hypothesis of Theorem \ref{thm:hyperbolicexponent} excludes: 
type $F_4$ and the extended Hesse group (no. 26 in the Shepherd-Todd list).
The former is a Coxeter group and the latter is an arrangement known to have the
same discriminant as the Coxeter group of type $B_3$ (in the sense of Corollary 
\ref{cor:discriminant}). Since we deal with Coxeter groups in a more concrete manner in
the next Subsection \ref{subsect:hecke}, we shall have covered these cases as well.
\end{remark}

\subsection{A Hecke algebra approach to the case with a 
Coxeter symmetry}\label{subsect:hecke}
The monodromy representation of $\nabla^\kappa$ and its invariant  
form $h^\kappa$ can be determined up to equivalence 
in case the  Dunkl connection is associated to a finite Coxeter group.

Let $W$ be an irreducible finite reflection group in a real vector space 
$V(\RR)$ without a nonzero fixed vector.
We take for $\Hcal$ the collection of reflection hyperplanes of $W$ in $V$
and for $H\in \Hcal$, we let $\pi_H=\half (\one-s_H)$, where $s_H$ is 
the reflection in $H$.
Choose $\kappa\in\RR^\Hcal$ to be $W$-invariant. We know that then 
$\nabla^\kappa$ is a flat $W$-invariant connection. We account 
for the $W$-invariance by regarding  $\nabla^\kappa$ as a connection 
on the tangent bundle
of $V^\circ_W$ (the group $W$ acts freely on $V^\circ$). So if we fix 
a base point $*\in V/W^\circ$, then we 
have a monodromy  representation 
$\rho^\mon\in\Hom (\pi_1(V_W^\circ,*), GL(V))$. 
It is convenient to let the base point be the image of a real point 
$x\in V(\RR)^\circ$.
So $x$ lies in a chamber $C$ of $W$. Let $I$ be a set that
labels the (distinct) supporting hyperplanes of $C$: $\{ H_i\}_{i\in I}$
and let us write $s_i$ for $s_{H_i}$. 
Then $I$ has $\dim V$ elements. Let $m_{i,j}$ denote the order of 
$(s_is_j)$, so that $M:=(m_{i,j})_{i,j}$ is the \emph{Coxeter matrix} 
of $W$. Then the \emph{Artin group} $\ar (M)$ associated to $M$ has a 
generating set 
$(\sigma_i)_{i\in I}$ with defining relations (the \emph{Artin relations})
\begin{equation*}
\underbrace{\sigma_i\sigma_j\sigma_i\cdots }_{m_{i,j}}=
\underbrace{\sigma_j\sigma_i\sigma_j\cdots }_{m_{i,j}},
\end{equation*}
where both members are words comprising $m_{i,j}$ letters.
The Coxeter group $W$ arises as a quotient of $\ar (M)$ by introducing the 
additional 
relations $\sigma_i^2=1$; $\sigma_i$ then maps to $s_i$. 
According to Brieskorn \cite{briesk} this lifts to an isomorphism of groups
$\ar (M)\to \pi_1(V_W^\circ,*)$ which 
sends $\sigma_i$ to the loop  is represented by the path in $V^\circ$ 
from $x$ to $s_i(x)$ which stays in the contractible set 
$V^\circ \cap (V(\RR)+\sqrt{-1}\bar C)$.  

As long as $|\kappa_i|<1$, $\rho^\mon(\sigma_i)$ is semisimple
and acts as a complex reflection over an angle $\pi (1+\kappa_i)$. So if we put
$t_i:=\exp (\half\pi\kappa_i\sqrt{-1})$, then $\sigma_i$ satisfies the identity
$(\sigma-1)(\sigma+t_i^2)=0$. Although the monodromy need not be semisimple for
$\kappa_i=1$, this equation still holds (for $t^2_i=-1$) .
In other words,  when $-1 <\kappa_i\le 1$, $\rho^\mon$
factors through the quotient of 
of the group algebra $\CC [\ar (M)]$ by the two-sided ideal
generated by the elements $(\sigma_i-1)(\sigma_i+t_i^2)$, $i\in I$.
These relations are called the \emph{Hecke relations} and the 
algebra thus defined is known as the \emph{Hecke algebra} attached to 
the matrix $M$ with parameters $t=(t_i)_i$. (It is more 
traditional to use the elements $-\sigma_i$ as generators;
for these the Artin relations remain valid, but  the 
Hecke relations take the form  $(\sigma_i+1)(\sigma_i-t_i^2)=0$.)
If the $t_i$'s are considered as variables (with $t_i=t_j$ if and only 
if $s_i$ and $s_j$ are conjugate in $W$), then this is an algebra
over the polynomial ring $\CC [t_i  |\,  i\in I]$. 

There are at most two conjugacy classes of reflections in $W$.
This results in a partition of $I$ into at most two subsets; we
denote by $J\subset I$ a nonempty part.
We have two conjugacy classes (i.e., $J\not=I$) only for 
a Coxeter group of type $I_2(\even )$,
$F_4$ and $B_{l\ge 3}$. We denote the associated variables $t$ 
and $t'$ (when the latter is defined).

If we put all $t_i=1$, then the Hecke algebra reduces to the group algebra 
$\CC[W]$, which is why the Hecke algebra for arbitrary parameters
can be regarded as a deformation of this group algebra.

For us is relevant the \emph{reflection representation} of the Hecke algebra
introduced in \cite{cik}. Since we want the reflections to be unitary
relative some nontrivial hermitian form we need to adapt this discussion 
for our purposes. We will work over the domain $R$ obtained from 
$\CC [t_i  |\,  i\in I]$ by adjoining the square root of $(t_it_j)^{-1}$ 
for each pair $i,j\in I$. So either $R=\CC [t,{t}^{-1}]$ or
$R:=\CC [t, t', (tt')^{-1/2}]$, depending on whether $W$ has one or two
conjugacy classes of reflections. So $R$ contains $t_i^kt_j^l$ 
if $k$ and $l$ are half integers which differ by an integer.
So $T:=\spec (R)$ is a torus of dimension one or two. 
Complex conjugation in $\CC$ extends to an 
anti-involution $r\in  R\mapsto \overline{r}\in R$ which sends  $t_i$ to 
$t_i^{-1}$ and $(t_it_j)^{1/2}$ to $(t_it_j)^{-1/2}$. 
This gives $T$ a real structure for which $T$
is anisotropic (i.e., $T(\RR)$ is compact). We denote by
$\Re :R\to R$  `taking the real part': $\Re (r):=\half (r+\overline{r})$.

Let $\Hcal (M)$ stand for the Hecke algebra as defined above  with 
coefficients taken in $R$ (so this is a quotient of $R[\ar (M)]$).
For $i,j\in I$ distinct, we define a real element of $R$:
\[
\lambda_{i,j}:=\Re \big(\exp (\pi\sqrt{-1} /m_{i,j})
t_i^{1/2}t_j^{-1/2}\big).
\]
Notice that $\lambda_{i,j}=\cos (\pi/m_{i,j})$ if $t_i=t_j$. 
If $W$ has two orbits in $\Hcal$, then there is a unique pair 
$(j_0,j_1)\in J\times (I-J)$ with $m_{j_0,j_1}\not=2$. Then $m_{j_0,j_1}$ must
be even and at least $4$ and we write $m$ for $m_{j_0,j_1}$,
and $\lambda$ resp.\ $\lambda'$ for 
$\lambda_{j,j'}$ resp.\ $\lambda_{j',j}$. So $\lambda= 
\Re (\exp (\pi\sqrt{-1}/m)t^{1/2}t'{}^{-1/2})$ and
$\lambda'=\Re(\exp (\pi\sqrt{-1}/m)t^{-1/2}t'{}^{1/2})$.

Define for every $i\in I$ a linear form $l_i: R^I\to R$ by
\[
l_i(e_j)=\begin{cases}
1+t_i^2 &\text{ if  $i=j$,}\\
-2\lambda_{i,j}t_i &\text{ if  $i\not= j$.}
\end{cases}
\]
Let $\rho^\refl(\sigma_i)$ be the pseudoreflection  in $R^I$ defined by 
\[
\rho^\refl(\sigma_i)(z)= z-l_i(z)e_i.
\]
We claim that this defines a representation of $\Hcal (M)$. 
First observe that the minimal polynomial of $\rho^\refl(\sigma_i)$ is $(X-1)(X+t_i^2)$.
For $i\not=j$, we readily verify that
\[
l_i(e_j)l_j(e_i)=t_i^2+t_j^2+2t_it_j\cos(2\pi /m_{i,j}),
\]
This implies that the trace of
$\rho^\refl (\sigma_i)\rho^\refl(\sigma_j)$ on the plane spanned by
$e_i$ and $e_j$ is equal to $2t_it_j\cos(2\pi/m_{i,j})$. Since its 
determinant is $t_i^2t_j^2$, it follows that the eigenvalues of 
$\rho^\refl (\sigma_i)\rho^\refl (\sigma_j)$ in this plane are $t_it_j\exp(2\pi\sqrt{-1}/m_{i,j})$ and
$t_it_j\exp(-2\pi\sqrt{-1}/m_{i,j})$. In particular
$\rho^\refl(\sigma_i)$ and $\rho^\refl (\sigma_j)$ 
satisfy the Artin relation.  So $\rho^\refl$ defines a 
representation of $\Hcal (M)$.

\begin{lemma}
Fix a $p\in T$ and consider the reflection representation of
the corresponding specialization $\Hcal(M)(p)$ on $\CC^I$.
Then $(\CC^I)^{\Hcal(M)(p)}$ is the kernel of the associated 
linear map $(l_i)_i :\CC^I\to \CC^I$. Moreover, if $K$ is a proper 
invariant subspace of $\CC^I$ which is not contained in 
$(\CC^I)^{\Hcal(M)(p)}$, then
$J\not= I$ and $\lambda\lambda'=0$ and $K$ equals $\CC^J$ resp.\ $\CC^{I-J}$
modulo $(\CC^I)^{\Hcal(M)(p)}$ when $\lambda'=0$ resp.\ $\lambda=0$.
\end{lemma}
\begin{proof}
The first statement is clear.

Since $K\not\subset (\CC^I)^{\Hcal(M)(p)}$, some $l_i$ with 
will be nonzero on $K$; suppose this happens for $i\in J$.
Let $z\in K$ be such $l_i(z)\not= 0$. From
$z-\rho^\refl(\sigma_i)(z)=l_i(z)e_i$ it follows that $e_i\in K$. 
Since $t\not=0$, our formulas imply that then $K\supset\CC^{J}$. 
Since $K$ is a proper subspace of $\CC^I$, $J\not=I$ and 
$l_j$ vanishes on $K$ for all $j\in I-J$ (otherwise the same argument 
shows that $K\supset \CC^{I-J}$). This implies in particular that 
$\lambda'=0$.
\end{proof}

By sending $\kappa_H$ to $e^{\half\pi\sqrt{-1}\kappa_H}$ we obtain a 
universal covering
\[
\tau : (\CC^\Hcal)^W\to T.
\]
Let $\Delta\subset (\CC^\Hcal)^W$ denote the locally finite union
of affine hyperplanes defined by: $\kappa_H\in\ZZ$
and $\kappa_0\in \{0,-1,-2,\dots\}$.

\begin{proposition}\label{prop:heckedunkl}
The map $\tau$ lifts to a holomorphic intertwining morphism
$\tilde \tau$ from the monodromy representation $\rho^\mon$ of $\ar (M)$ 
to the reflection representation $\rho^\refl$ of $\Hcal (M)$ in such a manner 
that it is an isomorphism away from $\Delta$ and nonzero away from
a codimension two subvariety $(\CC^\Hcal)^W$ contained in $\Delta$.
\end{proposition}
\begin{proof}
Suppose first $\kappa\notin\Delta$.

Since each $\kappa_H$ is nonintegral, $\rho^\mon(\sigma_i)$ is 
semisimple and acts in $V$ as a complex reflection (over an angle 
$\pi (1+\kappa_i)$). Hence $\one - \rho^\mon(\sigma_i)$ is of the form
$v_i\otimes f_i$ for some $v_i\in V$ and $f_i\in V^*$.
The individual $f_i$ and $v_i$ are not unique, only their tensor product is.
But we have  $f_i(v_i)=1+t_i^2=l_i(e_i)$ and  
the fact that $\sigma_i$ and $\sigma_j$ satisfy the Artin relation implies
that $f_i(v_j)f_j(v_i)=t_i^2+t_j^2+2t_it_j\cos(2\pi /m_{i,j})=
l_i(e_j)l_j(e_i)$. 

We claim that the $v_i$'s are then independent and hence form a basis of $V$.
For if that were not the case, then there would exist a nonzero 
$\phi\in V^*$ which vanishes on all the $v_i$'s. This $\phi$ will be clearly 
invariant under the monodromy representation. But this is prohibited
by  Corollary \ref{cor:invarform} which says that then $\kappa_0-1$ 
must be a negative integer. 

Since the Coxeter graph is a tree, we can put a total order on $I$ such that
that if $i\in I$ is not the smallest element, there is precisely one
$j<i$ with $m_{i,j}\not= 2$.
Our assumption implies that whenever $m_{i,j}\not=2$, at least one of 
$\lambda_{i,j}$ and $\lambda_{j,i}$ is nonzero.
This means that in such a case one of $l_i(e_j)$ and $l_j(e_i)$
is nonzero. On the other hand, it is clear that $l_i(e_j)= 0$ when 
$m_{i,j}=2$. We can now choose $f_i$ and $e_i$ in such a manner that 
$f_i(v_j)=l_i(e_j)$ for all $i,j$: proceed by induction on $i$:
The fact that for exactly one $j<i$ we have
that one of $l_i(e_j)$ and $l_j(e_i)$ is nonzero can be used to fix
$v_i$ or $f_i$ and since $v_i\otimes f_i$ is given, one 
determines the other. This prescription is unambiguous in case both 
$l_i(e_j)$ and $l_j(e_i)$ are nonzero, for as we have seen, 
$f_i(v_j)f_j(v_i)=l_i(e_j)l_j(e_i)$. 

We thus obtain an intertwining isomorphism
$\tilde\tau (\kappa): V\to\CC^I$,  $e_i\mapsto v_i$, which 
depends holomorphically on $\kappa$ and is meromorphic along $\Delta$.
Since we are free to multiply $\tilde\tau$ by a meromorphic function on 
$(\CC^\Hcal)^W$, we can arrange that $\tilde\tau$ extends holomorphically and
nontrivially over the generic point of each irreducible component of 
$\Delta$. 
\end{proof}

\begin{remark} 
With a little more work, one can actually show that the preceding
proposition remains valid if we alter the definition of $\Delta$ by
letting $\kappa_H$ only be an \emph{odd} integer. 
\end{remark}

We define a hermitian form $H$ on $R^I$ (relative to our 
anti-involution) preserved by $\rho^\refl$. This last condition
means that we want that for all $i\in I$,
\[
l_i(z)H(e_i ,e_i)=(1+t_i^2)H(z ,e_i ).
\]
In case all the reflections of $W$ belong to a single conjugacy 
class so that all $t_i$ take the same value $t$, then the form defined by
\[
H(e_i ,e_j) :=
\begin{cases}
\Re (t)&\text{ if  $i=j$,}\\
-\cos(\pi/m_{i,j})& \text{ if  $i\not=j$}  
\end{cases}
\]
is as desired. In case we have two conjugacy classes of reflections, then
\[
H(e_i ,e_j)= 
\begin{cases}
\lambda'\Re(t)&\text{if  $i=j\in J$,}\\
\lambda\Re(t')&\text{if  $i=j\in I-J$,}\\
-\lambda'\cos(\pi/m_{i,j})& \text{if $i,j\in J$  are distinct,}\\
-\lambda\cos(\pi/m_{i,j})& \text{if $i,j\in I-J$  are distinct,}\\
-\lambda\lambda'\cos (\pi/m_{i,j})&\text{otherwise.}
\end{cases}
\]
will do. If we specialize in some $p\in T$, then
the kernel of $H$ is of course $\Hcal (M)(p)$-invariant.  
If $\lambda'(p)=0$ resp.\ $\lambda (p)=0$, then the formulas show 
that this kernel contains $\CC^J$ resp.\ $\CC^{I-J}$.
The zero loci of $\lambda'$ and $\lambda$ are disjoint 
and so no specialization of $H$ is trivial, unless $I$ is a singleton
and $t^2=-1$.

We conclude from Proposition \ref{prop:heckedunkl}:

\begin{corollary}
Suppose that $\kappa$ takes values in $(0,1)$. Then the
monodromy representation is isomorphic to the reflection
representation and thus comes via such an isomorphism with a
nonzero $W$-invariant hermitian form.
\end{corollary}

At points where all the $t_i$'s take the same value
(so this is all of $T$ in case $J=I$ and the locus
defined by $t=t'$ otherwise), 
there is a neat formula for the determinant of $H$, which goes back 
to Coxeter and appears as Exercise 4 of Ch.\ V, \S\ 6 in Bourbaki 
\cite{bourbaki}:  
\[
\det (H(e_i,e_j)_{i,j})=\prod_{j=1}^{\vert I\vert}(\Re (t)-\cos (\pi m_j/h)),
\]
where $h$ is the Coxeter number  of $W$ and the $m_j$'s are the exponents of $W$.
Since $\re (t)=\cos (\half\pi\kappa)$. So if $t=\exp (\half\sqrt{-1}\pi\kappa)$,
we see that $H$ is degenerate precisely 
when $\kappa/4\equiv m_j/2h \pmod{\ZZ }$ for some $m_j$. Since the $m_j$'s
are distinct and in the interval $\{1,\dots ,h-1\}$, the nullity of $H$ is $1$ in
that case. The cardinality of $\Hcal$ is $h\vert I\vert /2$ (\cite{bourbaki}, Ch.\  V,\S\ 6,
no.\ 2, Th.\ 1), so that $\kappa_0=h\kappa/2$. Hence
$H$ is degenerate precisely when $\kappa_0\equiv m_j\pmod{2h\ZZ}$.
If we combine this with the results of Subsection \ref{subs:admrange} and 
\ref{prop:heckedunkl}, we find:

\begin{corollary}\label{cor:criticalexponents}
In case $\kappa:\Hcal\to (0,1)$  is constant, then the flat hermitian form of 
the associated Dunkl connection is degenerate precisely when  $\kappa_0$ equals 
some exponent $m_j$. In particular, $m_2$ is the hyperbolic exponent.
\end{corollary}

This raises the following

\begin{question} Assuming that $I$ is not a singleton, 
can we find a system of generators $X_1,\dots , X_{\vert I\vert}$ 
of the $\CC [W]$-module of $W$-invariant vector fields on $V$ of the correct degrees
$(m_1-1,\dots ,m_{\vert I\vert}-1)$ such that the ones in degree $m_j$
generate the kernel of the flat hermitian metric we found for the constant
map $\kappa:\Hcal\to (0,1)$ characterized by $\kappa_0=m_j$?

It makes sense to ask this question more generally for a  
complex reflection group (where we should then take the co-exponents as the appropriate
generalization). (We checked by an entirely different technique that the hermitian form
attached to a constant  map $\kappa :\Hcal\to (0,1)$ is degenerate precisely when 
$\kappa_0$ is a co-exponent, at least when  the group is primitive of rank at least 
three.)
\end{question}

\subsection{A flat hermitian form for the 
Lauricella system}\label{subsect:lauricellaflat}
Let $H$ be a monodromy invariant hermitian form on the translation space of
$A$ and denote by $h$ the corresponding flat hermitian form on $V^\circ$. 
Suppose that $\kappa_0\not=1$, so that we can think of $H$ as a 
hermitian form on the vector space $(A,O)$. Then the  associated `norm squared'
function, $H(a,a)$, evidently  determines $H$. So if we
view $H$ as a translation invariant form on $A$, then we can express as:
$\frac{1}{2}\sqrt{-1}\partial\bar\partial (H(E_A,E_A))=\im (H)$,
where  $E_A$ is the Euler vector field on $(A,O)$. Since the developing map
sends $E_V$ to $(1-\kappa_0)E_A$,
this property is transfered to $V^\circ$ as: if $N: V^\circ\to\RR$ 
is defined by $N:=h(E_V,E_V)$, then
\[
\frac{\sqrt{-1}}{2}\partial\bar\partial N= |1-\kappa_0 |^2\im (h).
\]
So if $h$ is nondegenerate, then the Dunkl connection is also determined by 
$N$. It would be interesting to find $N$ explicitly, or at least to characterize
the functions $N$ on $V^\circ$ that are thus obtained. We can do this
for the Lauricella example:

We consider the Lauricella system \ref{subsect:lauricella}.
For the moment we choose all the parameters 
$\mu_i\in (0,1)$ as usual, but we now also require that $\mu_0+\cdots +\mu_n>1$ 
(recall that here $\mu_0+\cdots +\mu_n=\kappa_0$).
We abbreviate the multivalued Lauricella differential by
\[
\eta:=(z_0-\zeta )^{-\mu_0}\cdots (z_n-\zeta )^{-\mu_n}d\zeta.
\]
Notice that $\eta\wedge \overline{\eta}$ is univalued $2$-form and that the conditions
imposed on the $\mu_i$'s guarantee that it is integrable, provided that
$(z_0,\dots ,z_n)\in V^\circ$. Since $\frac{\sqrt{-1}}{2}d\zeta \wedge d\bar \zeta$ is 
the area element of $\CC$, 
\[
N(z_0,\dots ,z_n):=-\frac{\sqrt{-1}}{2}\int_\CC \eta\wedge \overline{\eta}
\]
is \emph{negative}.
We will show that $N$ is a hermitian form in Lauricella functions. 
This implies that the Levi form of $N$ is flat and hence defines a flat hermitian 
form on $V^\circ$. 

For this purpose, let  $\gamma$ be an smoothly embedded oriented interval on 
the Riemann sphere which connects $z_0$ with $\infty$ and passes through 
$z_1, \dots ,z_{n}$ (in this order). On the complement of $\gamma$,
$\eta$ is representable by a holomorphic univalued differential which we extend to 
$\CC-\{ z_0,\dots ,z_n\}$ by taking on $\gamma$ the limit `from the left'. 
We continue to denote this differential by $\eta$, but this now makes $\eta$ discontinuous
along $\gamma$: 
its limit from the right on the stretch $\gamma_k$ 
from $z_{k-1}$ to $z_k$ (read $\infty$ for $z_{n+1}$) is easily 
seen to be $e^{-2\pi\sqrt{-1}(\mu_0+\cdots +\mu_{k-1})}$
times $\eta$. We find it convenient to put $w_{0}=1$ and $w_k:=
e^{\pi\sqrt{-1}(\mu_0+\cdots +\mu_{k-1})}$ for $k=1,\dots ,n$
so that the limit in question can be
written $\bar w_k^{2}\eta$. We put
\[
F(\zeta ):=\int_{\infty}^\zeta  \eta,
\]
where the path of integration is not allowed to cross $\gamma$.
So $F$ is holomorphic on $\CC-\gamma$ and continuous along $\gamma$ from 
the left. In case $z_0, \dots ,z_{n}$ are all real and ordered by size, then
a natural choice for $\gamma$ is the straight line on the real 
axis which goes from $z_0$ in the positive direction to $\infty$. Then
on $\gamma_k$ (the positively oriented interval $[z_{k-1},z_k]$) a natural
choice of determination of the integrand is the one which is real and positive:
$\eta_k:=(\zeta -z_0)^{-\mu_0}\cdots (\zeta -z_{k-1})^{-\mu_{k-1}}
(z_k-\zeta )^{-\mu_k}\cdots (z_n-\zeta )^{-\mu_n}d\zeta$. As 
$\eta_k=\bar w_k\eta$, this suggests to introduce 
\[
F_k := \bar w_k\int_{\gamma_k}\eta, \quad k=1,\dots ,n+1,
\]
in general. This is a Lauricella function  (up to scalar factor), and so is $F_0:=F(z_0)$.
For $\zeta\in\gamma_k$, $k=1,\dots ,n+1$, we have
\[
F(\zeta )=\sum_{j=0}^{k-1}w_{j}F_j(z)+ \int_{z_{k-1}}^\zeta  \eta .
\]

\begin{lemma}\label{lemma:lauricellavolume} 
Under the above assumptions (so $\mu_k\in (0,1)$ for all $k$ and 
$\sum_{k=0}^n \mu_k >1$) the Lauricella functions $F_k$ satisfy the linear relation
$\sum_{k=1}^{n+1} \im (w_k)F_k=0$ and we have 
$N(z)=\sum_{1\le j<k\le n+1} \im(w_j \bar w_k)\bar F_jF_k$.
\end{lemma}
\begin{proof}
If $\zeta \in\gamma_k$, then the limiting value of $F$ in $\zeta $ 
from the right 
is equal to
\[
\sum_{j=0}^{k-1}\bar w_jF_j+ \bar w_k^2\int_{z_{k-1}}^\zeta  \eta.
\]
The fact that the value of $F$ at $\infty$ is thus calculated 
in two ways yields the identity
$\sum_{k=1}^{n+1} \bar w_kF_k=\sum_{k=1}^{n+1} w_kF_k$ or
what amounts to the same $\sum_{k=1}^{n+1} \im (w_k)F_k=0$.

Now $N(z_0,\dots ,z_n)$ is the integral of the exterior derivative
of the $1$-form $\frac{\sqrt{-1}}{2}\bar{F}\eta$. 
If $\eta$ is the $1$-form on $\gamma$ which is the difference between 
$\bar{F}\eta$ and its limiting value from the right, then 
the theorem of Stokes implies that $N(z_0,\dots ,z_n)=
\frac{\sqrt{-1}}{2}\int_\gamma\eta$. 
The above computations show that on $\gamma_k$, $\eta$ is equal to
\begin{multline*}
\Big (\sum_{j=1}^{k-1} \bar w_j\bar F_j+\int_{z_{k-1}}^\zeta  
\bar\eta\Big)\eta -
\Big (\sum_{j=1}^{k-1}w_j\bar F_j+ \int_{z_{k-1}}^\zeta  
w_k^2\bar\eta\Big)\bar w_k^2\eta =\\
=\sum_{j=1}^{k-1}(\bar w_j-w_j\bar w_k^2)\bar F_j\eta=
-2\sqrt{-1}\sum_{j=1}^{k-1}\im (w_j\bar w_k)\bar F_j\bar w_k\eta
\end{multline*} 
and hence 
\[
N(z_0,\dots ,z_n)= \frac{\sqrt{-1}}{2}\int_\gamma\eta
=\sum_{1\le j<k\le n+1} \im(w_j \bar w_k)\bar F_jF_k
\]
\end{proof}

Let us think of $F_1,\dots ,F_{n+1}$ as linear functions on the receiving space $A$
of the developing map that satisfy the linear relation 
$\sum_{i=1}^{n+1} \im(w_k)F_k=0$. The preceding lemma
tells us that $N$ defines a hermitian form on $A$ that is invariant under the 
holonomy group. This suggests to consider for any $(n+1)$-tuple
$w=(w_1,\dots ,w_{n+1})$ of complex numbers of norm one 
that are not all real, the hyperplane $A_w$ of $\RR^{n+1}$ 
with equation $\sum_{k=1}^{n+1} \im(w_k)a_k=0$ and the quadratic form on 
$\RR^{n+1}$ defined by
\[
Q_w(a):=\sum_{1\le j<k\le n+1}\im(w_j\bar w_k)(a_ja_k).
\]
We determine the signature of $Q_w$. 

\begin{lemma}\label{lemma:lauricellaform}
Let us represent $w_1,\dots ,w_{n+1}$ by real numbers
$\mu_0,\dots ,\mu_n$ as before, so  $w_k=e^{\pi\sqrt{-1}(\mu_0+\cdots +\mu_{k-1})}$.
Then the \emph{nullity} (that is, the number of zero eigenvalues)  
of $Q_w$ on $A_w$
is equal to the  number of integers in the sequence $\mu_0,\dots ,\mu_n,
\sum_{i=0}^n\mu_i$ and
its \emph{index} (that is, the number of negative eigenvalues) is equal to
$[\sum_{i=0}^n\mu_i]-\sum_{i=0}^n [\mu_i]$.
\end{lemma}
\begin{proof}
It is clear that $(A_w,Q_w)$ only depends on the reduction of  
$\mu_0,\dots ,\mu_n$
modulo $2$, but the isomorphism 
type of $(A_w,Q_w)$ only depends on their reduction modulo $1$:  if we 
replace $\mu_k$ by $\mu_k+1$, then the new values $w'_j$ of $w_j$ are:
$w'_j= w_j$ for $j\le k$ and  $w'_j= -w_j$ for $j>k$ and we note that
$(a_1,\dots ,a_{n+1})\mapsto (a_1,\dots ,a_k,-a_{k+1},\dots ,-a_{n+1})$ 
turns $(A_w,Q_w)$ into $(A_{w'},Q_{w'})$. So without loss of generality 
we may assume that $0\le\mu_k<1$ for all $k$. 

We proceed by induction on $n\ge 0$. As the case $n=0$ is trivial, 
we suppose $n\ge 1$ and the lemma proved for smaller values of $n$.
This allows us to restrict ourselves to the case when
$0<\mu_k<1$ for all $k$:  if $\mu_k=0$, then $w_k=w_{k+1}$ and so if
$w':=(w_1,\dots ,w_k, w_{k+2},\dots ,w_n)$, then $(A_w,Q_w)$ is the pull-back of
$(A_{w'},Q_{w'})$ under $(a_1,\dots ,a_n)\mapsto 
(a_1,\dots ,a_{k-1},a_k+a_{k+1}, a_{k+2},\dots ,a_n)$.

We now let $w':= (w_1,\dots ,w_n)$.
First assume that $w_n\notin\RR$ so that $\sum_{k=0}^{n-1}\mu_k\notin\ZZ$.
According to our induction hypothesis this means
that $(A_{w'},Q_{w'})$ is nondegenerate of index $[\sum_{i=0}^{n-1}\mu_i]$.
There exist unique $s,t\in\RR$ such that $w_{n+1}=s w_n+t$. The fact
that $0<\mu_n<1$ implies that $t\not= 0$.
We set $a':=(a_1,\dots ,a_{n-1},a_n+sa_{n+1})$. Then we have
\[
\sum_{k=1}^{n+1} \im (w_k)a_k -\sum_{k=1}^n \im (w'_k)a'_k=
\im (w_{n+1}a_{n+1}-w_na_{n+1}s)=\im (ta_{n+1})=0
\]
so that $a\in A_w$ if and only if $a'\in A_{w'}$. A similar calculation shows that
\[
Q_w(a)=Q_{w'}(a')- t\im (w_{n+1})a_{n+1}^2,\quad a\in A_w.
\]
If $w_{n+1}\notin\RR$, then 
from the equality $t=-sw_n+w_{n+1}$ and the fact that $-w_n$ makes a positive angle
(less than $\pi$) with $w_{n+1}$,  we see that $t \im (w_{n+1})> 0$ if and only if 
$\im (w_n)$ and $\im(w_{n+1})$ have different sign.  The latter amounts
to $[\mu_0+\cdots +\mu_n]=[\mu_0+\cdots +\mu_{n-1}]+1$, and so here 
the induction hypothesis yields the lemma for $(A_{w},Q_{w})$. This is also
the case when $w_{n+1}\in\RR$, for then $\sum_{i=0}^n \mu_i\in\ZZ$.

Suppose $w_n\in\RR$, in other words, that $\sum_{i=0}^{n-1}\mu_i\in\ZZ$.
If we let $w''=(w_1,\dots ,w_{n-1})$, then 
$Q_{w'}(a_1,\dots ,a_n)=Q_{w''}(a_1,\dots ,a_{n-1})$. 
We  may assume that $n\ge 2$, so that $A_{w''}$ is defined. 
By induction, $(A_{w''}, Q_{w''})$ is
nondegenerate of index $[\sum_{i=0}^{n-2}\mu_i]$. It is now easy to check that
$(A_w, Q_w)$ is isomorphic to the direct sum of $(A_{w''}(\RR), Q_{w''})$
and a hyperbolic plane. Hence $(A_w, Q_w)$ is nondegenerate of index
$[\sum_{i=0}^{n-2}\mu_i]+1$. This last integer is equal to 
$\sum_{i=0}^{n-1}\mu_i$ and hence also equal to $[\sum_{i=0}^n\mu_i]$.
\end{proof}

\begin{corollary}
The function $N$ defines an invariant hermitian form on the Lauricella system
whose isomorphism  type is given by Lemma \ref{lemma:lauricellaform}. If
$0< \mu_k<1$ for all $k$, then the form is admissible of elliptic, parabolic, 
hyperbolic type for $\kappa_0<1$, $\kappa_0=1$, $1<\kappa_0<2$ respectively.
\end{corollary}
\begin{proof}
All the assertions follow from Lemma's \ref{lemma:lauricellaform}, except 
the admissibility statement. For the hyperbolic range $1<\kappa_0<2$, admissibility
follows from the fact that $N$ is negative in that case. For $\kappa=1$, 
Lemma \ref{lemma:lauricellaform} says that the hermitian form is 
positive semidefinite with nullity one. According to Theorem \ref{thm:semipos} this kernel is then 
spanned by the Euler vector field and so we have admissibility in this case,
too.
\end{proof}

\begin{remark}
In the hyperbolic case: $\mu_i\in (0,1)$ for all $i$ and $\sum_i \mu_i \in (1,2)$, 
we observed with Thurston in Subsection \ref{subsect:thurston}
that $\PP(V^\circ)$ can be understood as the moduli space of Euclidean 
metrics on the sphere with $n+2$ conical singularities with a prescribed total angle. 
The hyperbolic form induces a  
natural complex hyperbolic metric on $\PP(V^\circ)$.
The modular interpretation  persists on the metric completion of $\PP(V^\circ)$:
in this case  we allow some of the singular points to collide, that is, we may include 
some the diagonal strata. This  metric completion is quite special and is of the 
same nature as the objects it parametrizes: it is what Thurston calls a cone manifold. 
\end{remark}

\begin{remark}
If each $\mu_i$ is positive and rational, then the associated 
Lauricella system with its hermitian form can 
also be obtained as follows. Let $q$ be a common 
denominator, so that $\mu_i:=p_i/q$ for some positive integer $p_i$, and put 
$p:=\sum_i p_i$. Consider the Dunkl system on the Coxeter arrangement of type
$A_{p-1}$ defined by the diagonal hyperplanes in the hyperplane $V_p$ in
$\CC^p$ defined by $\sum_{i=1}^q z_i=0$
and with $\kappa$ constant equal to $1/q$. Let $V_P\subset V_p$ be the intersection
of hyperplanes defined by the  partition $P:=(p_0,p_1,\dots ,p_n)$ of $p$.
Then the Lauricella system can be identified with longitudinal system on $V_P$.
The hermitian form that we have on the  ambient system via the Hecke algebra approach 
\ref{subsect:hecke} is inherited by $V_P$ (as a flat hermitian form). This 
approach is taken (and consistently followed) by B.~Doran in his thesis
\cite{doran:thesis}.  
\end{remark}

\subsection{The degenerate hyperbolic case}\label{subsect:degenhyp}
By a \emph{degenerate hyperbolic form} on a vector space we simply mean a degenerate 
hermitian form which is a hyperbolic form on the quotient of this vector space by 
kernel of the form. If $H$ is such a form on the vector space $A$ with kernel $K$,
then the subset $\BB\subset \PP(A)$ defined by $H(a,a)<0$ is best understood as 
follows: since $H$ induces a nondegenerate form $H'$ on $A':=A/K$, there is a
ball $\BB'$ defined in $\PP (A')$ by $H'(a',a')<0$. The projection $A\to A'$ 
induces a morphism $\pi: \BB\to \BB'$ whose fibers are affine spaces
of the same dimension as $K$. The vector group $\Hom (A',K)$
acts as a group of bundle automorphisms of $\pi$ which act as 
the identity on $\BB'$ but this action is not proper. So 
if the holonomy preserves a form of this type it might not act 
properly on $\BB$.

Let us see what happens in the Lauricella case. 
We return to the situation of Subsection \ref{subsect:lauricellaflat} and 
choose 
$\mu_i\in (0,1)$ for $i=0,\dots ,n$  and such that $\sum_i\mu_i=2$.
We also let $w=(w_k:=e^{\pi\sqrt{-1}(\mu_0+\cdots +\mu_{k+1})})_{k=1}^{n+1}$,
$A_w\subset\RR^{n+1}$, the hyperplane defined by $\sum_i \im (w_i)a_i=0$,
and  $Q_w: A_w\to\RR$, 
$Q_w(a):=\sum_{1\le i<j\le n+1}\im (w_j\bar w_k)a_ja_k$ be as before.
Notice that $w_{n+1}=1$. 
According to Lemma \ref{lemma:lauricellaform}, $Q_w$ 
has a one dimensional kernel. In fact, if $w':=(w_1,\dots ,w_n)$, then 
omission of the last coordinate, 
$a=(a_1,\dots ,a_{n+1})\mapsto a':=(a_1,\dots ,a_n)$, defines a
projection $A_w\to A_{w'}$,
we have $Q_w(a)=Q_{w'}(a')$ and $Q_{w'}$ is nondegenerate of hyperbolic 
signature (see the proof of Lemma \ref{lemma:lauricellaform}).
This describes the situation at the receiving end of the developing map.
Now let us interpret this in the domain. 
The projection $A_w\to A_{w'}$ amounts to ignoring the Lauricella function
$F_{n+1}$; this is the only one among the $F_1,\dots ,F_{n+1}$ which involves
an integral with $\infty$ as end point. Observe that the condition 
$\sum_i\mu_i=2$ implies that $\infty$ is not a singular point of the 
Lauricella form 
$\eta =(z_0-\zeta )^{-\mu_0}\cdots (z_n-\zeta )^{-\mu_n}d\zeta$.
This suggests an invariance property with respect to M\"obius 
transformations.
This is indeed the case: a little exercise shows that 
\begin{math}
\left( 
\begin{smallmatrix}
a&b\\
c&d
\end{smallmatrix}
\right)\in \SL (2,\CC)
\end{math}
transforms $\eta$ into $(cz_0+d)^{\mu_0}\cdots (cz_n+d)^{\mu_n}\eta$.
Hence the first $n$ coordinates of the developing map
$(F_1,\dots ,F_{n+1})$ (with values in $A_w\otimes\CC$) 
all get multiplied by the same factor:
for $k=1,\dots ,n$ we have
\[
F_k\Big( \frac{az_0+b}{cz_0+d},\dots ,\frac{az_n+b}{cz_n+d} \Big) =
(cz_0+d)^{\mu_0}\cdots (cz_n+d)^{\mu_n}F_k(z_0,\dots ,z_n).
\]
In geometric terms this comes down to the following. Embed 
$\CC^{n+1}$ in $(\PP^1)^{n+1}$ in the 
obvious manner and let the M\"obius group $\PSL (2,\CC)$ act on 
$(\PP^1)^{n+1}$ diagonally.
This defines a birational action of  $\PSL (2,\CC)$ on $(\CC^{n+1})^\circ$.
Recall that $V^\circ$ stands for the quotient of $(\CC^{n+1})^\circ$ by the main
diagonal. The obvious map $(\CC^{n+1})^\circ\to\PP(V^\circ)$ is 
the formation of the orbit space with respect to the 
group of affine-linear transformations of $\CC$. Hence a 
$\PSL (2,\CC)$-orbit in $(\PP^1)^{n+1}$ which meets $(\CC^{n+1})^\circ$ 
maps to a rational curve in $\PP(V^\circ)$. Thus the fibration of
$\PP(V^\circ)$ can (and should) be thought of as the forgetful morphism 
$\Mcal_{0,n+2}\to \Mcal_{0,n+1}$ which ignores the last point: 
it is represented by 
$(\PP^1; z_0,\dots ,z_n,\infty)\mapsto (\PP^1; z_0,\dots ,z_n)$.
In particular, the fiber is an $(n+1)$-pointed rational curve;
it can be understood as the curve on which is naturally defined
the Lauricella form $\eta$ (up to a scalar multiple). Thus we have before us
the universal family for the Lauricella integral. We conclude:  

\begin{proposition}\label{prop:laurinfty}
The fibration $\Mcal_{0,n+2}\to \Mcal_{0,n+1}$ integrates the distribution 
defined by the kernel of the flat hermitian form so that we have a commutative 
diagram 
\[
\begin{CD}
\widetilde{\Mcal}_{0,n+2}@>>>\BB_w\\
@VVV @VVV\\
\widetilde{\Mcal}_{0,n+1}@>>> \BB_{w'}
\end{CD}
\]
where on the left we have the holonomy cover of 
$\Mcal_{0,n+2}\to \Mcal_{0,n+1}$ and on the right $\BB_w$ and $\BB_{w'}$
are the open subsets of $\PP (A_w\otimes \CC)$ resp.\ 
$\PP (A_{w'}\otimes\CC)$ defined by the hermitian forms.
\end{proposition}

The holonomy along a fiber of $\Mcal_{0,n+2}\to \Mcal_{0,n+1}$ is understood
as follows. Let $C:=\PP^1-\{z_0,\dots ,z_n\}$  
represent a point of $\Mcal_{0,n+1}$. The map $H_1(C;\ZZ)\to \RR$
which assigns to a small circle centered at $z_i$ the value $\mu_i$ defines
an abelian covering of $C$; it is a  covering on which the
Lauricella integrand becomes single valued.
Yet another abelian cover may be needed to make this single valued form exact. 
The resulting nilpotent cover $\widetilde{C}\to C$ appears as a fiber of 
$\widetilde{\Mcal}_{0,n+2}\to \widetilde{\Mcal}_{0,n+1}$ and the developing map
restricted to this fiber is essentially the function $\widetilde{C}\to\CC$  which 
integrates the Lauricella integrand.

\section{The Schwarz conditions}\label{sect:wellbehave}

\subsection{The Schwarz symmetry groups}\label{subsect:schwarz}
We begin with the  simple, but basic

\begin{example}\label{example:simple} 
Take $V$ of dimension $1$ so that $\Hcal$ consists of the origin.
If $z$ is a coordinate for $V$, then $\Omega = \kappa z^{-1}dz$ for 
some $\kappa\in\CC$. The new affine structure on $V^\circ=V-\{ 0\}$
is given by $z^{1-\kappa}$ ($\kappa\not=1$) or $\log z$ ($\kappa=1$). 
So in case $\kappa$ is irrational or equal to $1$, then the 
developing map defines an isomorphism of the universal cover of
$V-\{ 0\}$ onto an affine line. 

Suppose now $\kappa\in \QQ$, but distinct from $1$, and write 
$1-\kappa=p/q$ with $p,q$ relatively prime integers 
and $q>0$. The holonomy cover
extends with ramification over the origin as the $q$-fold cover
$\tilde V\to V$ defined by $w^q=z$. The developing map is the 
essentially given by $w\mapsto w^p$. So it extends across the origin 
only if 
$p>0$, that is, if $\kappa <1$, and it is injective only if $k=\pm 1$. 
This is why it would have been better if $V$ had been equipped
with the group of $p$th roots of unity $\mu_p$ as a symmetry group.
For then we can pass  to the orbit space of $V$ by this 
group: the $\mu_p$-orbit space of $V^{\circ}$ is covered by the 
$\mu_p$-orbit space of $\widetilde{V^{\circ}}$ and the developing 
map  factors through the latter as an open embedding.
This motivates the definition below. 
\end{example}

\begin{definition}\label{def:schwarz}
Given a Dunkl system for which $\kappa$ takes values in $\QQ$,
then we say that $L\in\Lcal_\irr (\Hcal)$
satisfies the \emph{Schwarz condition} if $1-\kappa_L$ is zero or
a nonzero rational number with the following property: if we write $1-\kappa_L=p_L/q_L$ 
with $p_L,q_L$  relatively prime and $q_L>0$, then the Dunkl system is invariant 
under the group $G_L$ of unitary transformations of $V$ which fix $L$ 
pointwise and act as scalar multiplication in $L^\perp$ by a $|p_L|$th root of unity.
We call $G_L$  the \emph{Schwarz rotation group} of $L$.
The \emph{Schwarz symmetry group} is the subgroup of 
the unitary group of $V$ generated by the Schwarz rotation groups $G_L$ of the
$L\in\Lcal_\irr (\Hcal)$ which satisy the Schwarz condition; we will usually denote 
it by $G$. We say that the Dunkl system \emph{satisfies the Schwarz condition
in codimension one} if every member of  $\Hcal$  satisfies the Schwarz condition.
We say that  the \emph{Dunkl system satisfies the Schwarz condition} if 
every $L\in \Lcal_\irr(\Hcal)$ satisfies the Schwarz condition.
\end{definition}

Notice that the Schwarz symmetry group is finite: this follows from the fact that
the group of projective-linear transformations of $\PP(V)$  which leave $\Hcal$ 
invariant is finite (since $\Hcal$ is irreducible) and 
the fact that the determinants of the generators of $G$ are roots of unity. 
This group may be trivial or be reducible  nontrivial  (despite
the irreducibility of $\Hcal$).  If the Schwarz symmetry group
is generated in codimension one, then according to Chevalley's theorem, the orbit space 
$G\bs V$ is isomorphic to affine space. 

It it clear that  $\{ 0\}$ always satisfies the Schwarz condition.

\begin{example}\label{example:int}
For the Lauricella system discussed in Subsection \ref{subsect:lauricella},
the Schwarz condition in codimension one amounts to: for $0\le i<j\le n$,
$1-\mu_i-\mu_j$ is a positive rational number with numerator $1$ or $2$ 
with $2$ only allowed if $\mu_i=\mu_j$. This last
possibility is precisely Mostow's $\Sigma$INT-condition \cite{mostow:int}.
\end{example}

Let $L\in\Lcal_\irr (\Hcal)$. If a Dunkl system satisfies the  Schwarz 
condition, then this property is clearly inherited by both the $L$-transversal 
Dunkl system. This is also true for the $L$-longitudinal Dunkl system:

\begin{lemma}
Suppose that the Dunkl system satisfies the Schwarz condition.
Then for every $L\in\Lcal_\irr (\Hcal)$, the longitudinal Dunkl system on
$L^\circ$ also satisfies  the Schwarz condition.
\end{lemma}
\begin{proof}
Let $M\in \Lcal_\irr (\Hcal^L)$. Either $M$ is irreducible in $\Hcal$ or
$M$ is reducible with two components $L$ and $M'$. The 
exponent of $M$ relative to $\Hcal^L$ is then  $\kappa_M$ and $\kappa_{M'}$
respectively. It is clear that the Schwarz symmetry group of $M$ resp.\ $M'$
preserves $L$.   
\end{proof}

\subsection{An extension of the developing map}\label{subsect:extension}
Every point  of $V$ determines a conjugacy class of subgroups in the 
fundamental group of
$V^{\circ}$ (namely the image of  the map on fundamental groups of
the inclusion in $V^{\circ}$ of the trace on $V^{\circ}$ of a small 
convex neighborhood of that point), hence also determines a conjugacy 
class in $\G$. If the latter is a conjugacy class of finite 
subgroups we say that we have \emph{finite holonomy at this point}.
The set $V^f\subset V$ of the points at which we have finite holonomy
is a union of $\Hcal$-strata which contains $V^{\circ}$ and is open in 
$V$ (the subscript $f$ stands for finite).
We denote the corresponding subset of $\Lcal(\Hcal)$ by
$\Lcal^f(\Hcal)$. Notice that the holonomy covering extends uniquely to a 
ramified $\G$-covering $\widetilde{V^f}\to V^f$.

If each $\kappa_{H}$ is rational $\not=1$, 
then $\Lcal^f(\Hcal)$ contains $\Hcal$ and so $V-V^f$ is everywhere of 
codimension $\ge 2$.

\begin{theorem}\label{thm:devmap}
Assume that $\kappa$ takes values in the rational numbers.
Then the Schwarz symmetry group $G$ acts freely on $V^{\circ}$ and 
lifts naturally to one on $\widetilde{V^f}$. The latter action commutes 
with the $\G$-action and the developing map is constant on $G$-orbits:
it factors through  a morphism $\ev_G :G\bs \widetilde{V^\circ}\to A$.

If $\kappa_0\not=1$ and $1-\kappa_0$ is written as a fraction
$p_0/q_0$ with $p_0,q_0$ relatively prime and $q_0>0$ as usual, 
then $G\cap \CC^{\times}$ consists of the $p_0$-th roots of unity and
both  $\widetilde{V^f}$ and $G\bs \widetilde{V^f}$ come with natural effective 
$\CC^{\times}$-actions such that $\widetilde{V^f}\to V^f$ is homogeneous of  
degree $q_0$, $\widetilde{V^\circ}\to G\bs \widetilde{V^\circ}$ is 
homogeneous of degree $p_0$ and  
$\ev_G: G\bs \widetilde{V^\circ}\to A$ is homogeneous of degree one.

In case $\kappa_0=1$, then the lift of the Euler vector field 
generates a free action of $\CC^{+}$ on $G\bs \widetilde{V^\circ}$  
such that $\ev_G$ is equivariant with respect to a 
one-dimensional translation subgroup of $A$.
\end{theorem}
\begin{proof}
Since $G$ preserves the Dunkl connection, it preserves the local 
system $\Aff_{V^\circ}$.
So $G$ determines an automorphism group $\G_G$ of 
$\widetilde{V^{\circ}}$ (with its affine structure) which contains 
the holonomy group $\G$ and has $G$ as quotient acting in the given 
manner on $V^{\circ}$. This group  acts on 
$A$ as a group of affine-linear transformations. 
Denote by $K$ the kernel of this representation. 
Since $\G$ acts faithfully on $A$, $K\cap\G=\{ 1\}$ and so
the map $K\to G$ is injective. On the other hand, 
if $L\in\Lcal_\irr\Hcal$ satisfies the Schwarz condition, then
the local model near the blowup of $L$ in $V$ shows that the
developing map is near $L$ constant on the $G_L$-orbits.
So  $G_L\subset K$ and hence $G\subset K$. This proves that
$\G_{G}$ is in fact the direct product of $\G$ and $G$.
It is now also clear that the developing map factors as asserted. 
Since the developing map is a local isomorphism on 
$\widetilde{V^{\circ}}$, the action of $G$ on $\widetilde{V^{\circ}}$ must be free. 

Suppose now $\kappa_0\not=1$.
The holonomy of $\Aff_{V^\circ}$ along a $\CC^{\times}$-orbit in $V^\circ$ is of 
order $q_0$ and so $\widetilde{V^{\circ}}$ comes with an 
effective $\CC^{\times}$-action  for which its projection to 
$V^{\circ}$ is homogeneous of degree $q_0$.
The developing map $\ev :\widetilde{V^{\circ}}\to A$ is constant on 
the orbits of the order $p_0$ subgroup of $\CC^{\times}$, but 
not for any larger subgroup. 
The infinitesimal generators of the 
$\CC^{\times}$-actions on $\widetilde{V^{\circ}}$ and $A$ are 
compatible and so $\ev$ is homogeneous of degree $p_0$
and there is a (unique) effective $\CC^{\times}$-action on 
$G\bs \widetilde{V^{\circ}}$ which makes $\widetilde{V^{\circ}}\to 
G\bs \widetilde{V^{\circ}}$ homogeneous of degree $p_0$.
Then $\ev_G: G\bs \widetilde{V^{\circ}}\to A$ will be homogeneous of degree 
one. These actions extend
to $\widetilde{V^f}$ and $G\bs \widetilde{V^f}$ respectively.

The last assertion follows from the fact that the holonomy along
a $\CC^\times$-orbit in $V$ is a nontrivial translation.
\end{proof}

\begin{theorem}\label{thm:devmap2}
Suppose that every $\kappa_H$ is a rational number smaller than $1$ and that 
the Dunkl system satisfies the  Schwarz condition in codimension one. 
Then the developing map $\widetilde{V^{\circ}}\to A$ extends to 
$\widetilde{V^f}$ and this extension drops to a local isomorphism 
$\ev_G :G\bs \widetilde{V^f}\to A$. In particular, $G\bs\widetilde{V^f}$ is smooth
and the $G$-stabilizer of a point of $\widetilde{V^f}$ acts near that
point as a complex reflection group.
Moreover,  every $L\in \Lcal_\irr(\Hcal)\cap \Lcal^f(\Hcal)$
satisfies the Schwarz condition and has $\kappa_{L}<1$.
\end{theorem}
\begin{proof}
The local model of the connection near the generic point of $H\in\Hcal$ 
shows that $H^\circ\subset V^f$ and  that the developing map extends over $H^\circ$ 
and becomes a local isomorphism if we pass to the $G_H$-orbit space.
So the developing map extends to $\widetilde{V^f}$ in codimension one. Hence
it extends to all of $\widetilde{V^f}$ and the resulting  extension of $\ev_G$
to $G\bs \widetilde{V^f}$ will even be a local isomorphism.

Now let $L\in \Lcal_\irr(\Hcal)\cap \Lcal^f(\Hcal)$. 
Then  the composite of $\ev$ with a generic 
morphism $(\CC ,0)\to (V,L^{\circ})$ is of the form 
$z\mapsto z^{1-\kappa_{L}}$ plus higher order terms
(for $\kappa_{L}\not=1$) or $z\mapsto \log z$ plus higher order terms
(for $\kappa_{L}=1$). As the developing map
extends over $L^\circ$,  we must have $\kappa_{L}<1$.
Since the developing map is in fact a local isomorphism at $L^\circ$, 
$L$ must satisfy the Schwarz condition.
\end{proof}

\begin{remark}
The orbit spaces  $G\bs V$ and $G\bs \widetilde{V^f}$ are both smooth. 
Notice that $G\bs V^f$ underlies two affine orbifold  structures.
One regards $G\bs V^f$ as a finite quotient of $V^f$ and has orbifold 
fundamental group $G$. Another inherits this structure from the
Dunkl connection, has $\ev_G: G\bs \widetilde{V^f}\to A$ as developing 
map and $\G$ as orbifold fundamental group.
\end{remark}

\section{Geometric structures of elliptic and parabolic type}
\label{sect:elliptic}

\subsection{Dunkl connections with finite holonomy}
In case $\G$ is finite, then the vector space $(A,O)$ admits a 
$\G$-invariant hermitian positive definite inner product. 
In particular, the tangent bundle of $V^{\circ}$ admits a 
positive definite inner product invariant under the holonomy group
of the Dunkl connection. Since the Dunkl connection is torsion free,
the latter is then the Levi-Civita connection of this metric.
Conversely:

\begin{theorem}\label{thm:finite2}
Suppose that $\kappa\in (0,1)^\Hcal$, that the Dunkl system satisfies the  
Schwarz condition in codimension one 
and that there is flat positive  definite hermitian form.
Then the holonomy of the affine structure defined by the Dunkl 
connection is finite and so we are in the situation where $\ev_G$ is a 
$\G$-equivariant isomorphism of $G\bs \widetilde{V}$ onto $A$
and $\kappa_0<1$. In particular, this map descends to an isomorphism
of orbit spaces of reflection groups $G\bs V\to \G\bs A$ via which
$\PP(G\bs V)$ acquires another structure as a complete elliptic orbifold.
\end{theorem}

The proof of Theorem \ref{thm:finite2} uses the following topological
lemma. We state it in a form that makes it applicable to other cases
of interest.

\begin{lemma}\label{lemma:cov}
Let $f: X\to Y$ be an continuous  map with discrete fibers 
between locally compact Hausdorff spaces and let $Y'\subset Y$ be an open subset
of which the topology is given by a metric.
Suppose that there is a symmetry group $\G$ of this situation 
(i.e., $\G$ acts on $X$ and $Y$, $f$ is $\G$- equivariant and $\G$ preserves $Y'$ 
and acts there as a group of isometries) for which the following properties hold:
\begin{enumerate}
\item[(i)] The action of $\G$ on $X$ is cocompact.
\item[(ii)] For every $y\in Y$  and  neighborhood $V$  of $y$ in $Y$ there exists
an $\varepsilon>0$ and a neighborhood  $V'$ of $y$ such that the
$\varepsilon$-neighborhood of $V'\cap Y'$ is contained in $V$. 
\end{enumerate}
Then there exists an $\varepsilon>0$ such that every 
$x\in f^{-1}Y'$ has a neighborhood which is proper over the
$\varepsilon$-ball in $Y'$ centered at $f(x)$.
In particular, if $f$ is a local homeomorphism over $Y'$ and $Y'$ is connected 
and locally connected, then $f$ is a covering projection over $Y'$.
\end{lemma}
\begin{proof}
Let $x\in X$. Since the fiber through $x$ is discrete,
we can find  a compact neighborhood $K$ of $x$
such that $f(x)\notin f(\partial K)$. 
Put $U_x :=K\setminus f^{-1}f(\partial K)$ and $V_x:=Y-f(\partial K)$
so that $U_x$ is a neighborhood of $x$, $V_x$  a neighborhood of $f(x)$
and $f$ maps $U_x$ properly to  $V_x$. By (ii) there exist
a neighborhood $V'_x$ of $f(x)$ and a $\varepsilon_x>0$ such that
such that for every $y\in V'_x\cap Y'$ the $\varepsilon_x$-neighborhood of 
$y$ is contained in $V_x$. We let $U'_x$ be the preimage of $V'_x$ in $U_x$.
It has the property that any $\varepsilon_x$-ball 
centered at a point of $f(U'_x)\cap Y'$ has a preimage in $U_x$ that is proper
over that ball.

Let $C\subset X$ be compact and such that $\G .C=X$.
Then $C$ is covered by $U'_{x_1},\dots ,U'_{x_N}$, say. We claim that
$\varepsilon :=\min_{i=1}^N\{\varepsilon_{x_i}\}$ has the required property.
Given any $x\in f^{-1} Y'$, then  $\g x\in U'_{x_i}$ for
some $i$ and $\g\in\G$. By construction, the
$\varepsilon$-ball centered at $f(\g x)$ is contained in $V_{x_i}$ and its
preimage in $U_{x_i}$ is proper over that ball. Now take the translate over $\g^{-1}$
and we get the desired property at $x$.
\end{proof}

\begin{proof}[Proof of Theorem \ref{thm:finite2}]
We have already verified this when $\dim (V)=1$. So we
take $\dim (V)\ge 2 $ and assume inductively the theorem proved
for lower values of $\dim (V)$. The induction hypothesis 
implies that $V^f$ contains $V'=V-\{ 0\}$. By 
Theorem \ref{thm:devmap2}  $\ev_G$ is then a local isomorphism on 
preimage $G\bs \widetilde{V'}$. On  $G\bs \widetilde{V'}$ we have an effective 
$\CC^\times$-action for which $\ev_G$ is homogeneous of nonzero
degree. Since $\ev_G$ is a local isomorphism, it maps $G\bs \widetilde{V'}$ 
to $A-\{ O\}$ and is the $\CC^{\times}$-action on $G\bs \widetilde{V'}$  
without fixed points.
So $\ev_G$ induces a local isomorphism  of $\CC^{\times}$-orbit spaces
$G\bs \PP (\widetilde{V'})\to \PP(A)$.
The action of $\G$ on $G\bs \PP(\widetilde{V'})$ is 
discrete and the orbit space of this action is a finite quotient of 
$\PP (V)$ and hence compact. So $G\bs \PP (\widetilde{V'})\to \PP(A)$ 
satisfies the hypotheses of Lemma \ref{lemma:cov} (with $Y'=Y=\PP (A)$), hence is a covering map.
Then $\ev_G :G\bs \widetilde{V'}\to A-\{ O\}$ is also a 
covering  map. But $A-\{ O\}$ is simply connected and so this must be an 
isomorphism. Such a map extends across the origin and so the degree
of homogeneity is positive: $1-\kappa_0>0$. It also follows
that the subgroup $\G$ of $\GL (A)$ acts properly discretely on 
$A-\{ O\}$ so that $\G$ is finite.
\end{proof}

\subsection{A remarkable duality}\label{subs:duality}
Suppose that the holonomy of the Dunkl connection is finite.
Then according to Theorem  \ref{thm:devmap2}, we have 
$\kappa_{L}<1$ for all $L\in \Lcal_{\irr}(\Hcal)$ and the developing 
map defines a isomorphism of  $G\bs V$ onto $A_{\G}$. 
So $G\bs V$ has two orbifold structures, one with
orbifold fundamental group $G$ , another with  $\G$. 

There is a simple relation between the invariant theory of
the groups $G$ and $\G$, which was observed earlier by
Orlik and Solomon \cite{orliksol} in a somewhat different
and more special setting.

The  $\CC^\times$-action on 
$(A,O)$ descends to a $\CC^\times$-action on $A_{\G}$ with 
kernel $\G\cap \CC^\times$. 
Let $1\le d_1(\G)\le d_2(\G)\le\cdots \le d_{\dim A}(\G)$ be the set of 
weights of this action, ordered  by size. The degrees $>1$ 
are the degrees of the basic invariants of $\G$. Their product 
$\prod_i d_i(\G)$ is the degree of $A\to A_{\G}$, that is, the order 
of $\G$. The situation for the $G$-action is likewise. The isomorphism 
between  the two orbit 
spaces is $\CC^{\times}$-equivariant once we pass to the corresponding
effective actions. This implies that the weights of these groups 
are proportional:
\[
d_i(\G)=(1-\kappa_0)^{-1}d_i(G ),\quad i=1,\dots ,\dim V.
\]
So the degrees of $\G$
are readily computed from the pair $(\kappa ,G)$.
In particular, we find that
\[
|\G |= (1-\kappa_0)^{-\dim V} |G|.
\]
The isomorphism $G\bs V\to \G\bs A$ maps the $G$-orbit space
of the union of the hyperplanes from $\Hcal$
onto a hypersurface in $A$ whose preimage in $A$ is a 
$\G$-invariant union of hyperplanes containing the reflection
hyperplanes of $\G$. If we denote that linear arrangement in
$A$ by $\Hcal'$, then we have bijection between the
$G$-orbits in $\Hcal$ and the $\G$-orbits in the $\Hcal'$. 

\smallskip
We can also go in the opposite direction, that is, start with the 
finite reflection group $\G$ on $A$ and define a compatible 
$\G$-invariant Dunkl connection 
on $A$ whose holonomy group is $G$ has a developing map 
equal to the inverse of the developing map of for the Dunkl connection 
on $V$. The following theorem exhibits the symmetry of the 
situation. At the same time it shows that all pairs of reflection groups
with isomorphic discriminants arise from Dunkl connections.

\begin{theorem}\label{thm:duality}
Let for $i=1,2$, $G_i\subset \GL (V_i)$ be a finite complex 
reflection group and $D_i\subset V_i$ its union of reflection 
hyperplanes. Then any isomorphism of orbit spaces 
$f: G_1\bs V_1\to G_2\bs V_2$ which maps $G_1\bs D_1$ onto $G_2\bs 
D_2$ and is $\CC^\times$-equivariant relative the natural \emph{effective}  
$\CC^\times$-actions on range and domain is obtained from 
the developing map of a $G_1$-invariant Dunkl connection on $V_1-D_1$ (and 
then likewise for $f^{-1}$, of course). 
\end{theorem}
\begin{proof}
The ordinary (translation invariant) flat connection on $V_2$ descends
to a flat connection on  $G_2\bs (V_2-D_2)$. Pull this back via $f$ 
to a flat connection on $G_1\bs (V_1-D_1)$ and lift the latter to a 
$G$-invariant flat connection $\nabla$ on  
$V_1-D_1$. It is clear that $\nabla$ is $\CC^\times$-invariant.
A straightforward local computation at the generic point of 
a member of the arrangement shows that $\nabla$ extends to the tangent bundle
of $V_1$ with a logarithmic poles and semisimple residues. So by
Corollary \ref{cor:dunklchar} it is a Dunkl connection.
It is clear that $f$ realizes its developing map.
\end{proof}

\begin{corollary}\label{cor:discriminant}
Let for $i=1,2$, $G_i\subset \GL (V_i)$ be a finite complex 
reflection group and $D_i\subset V_i$ its union of reflection hyperplanes.
If the germs of  $G_1\bs D_1$ and $G_2\bs D_2$
at their respective origins are isomorphic, then the two are related by
the above construction: one is obtained from the other by means of
the developing map of a Dunkl connection.
\end{corollary}
\begin{proof}
Any isomorphism
of germs $f: G_1\bs(V_1,D_1,0)\to G_2\bs (V_2,D_2,0)$ takes the 
effective $\CC^\times$-action on $G_1\bs V_1$ to an effective  
$\CC^\times$-action on the germ $G_2\bs (V_2,D_2,0)$. 
A finite cover of this action lifts to an effective action on the germ
$(V_2,D_2,0)$ which commutes with the action of $G_2$.  
Restrict this action to the tangent space of $V_2$ at the origin.
The fact that it preserves $D_2$ implies that it is just scalar 
multiplication
in $T_0V_2$. So if we identify this tangent space with $V$, then we 
get another isomorphism $f_0: (G_1\bs V_1,D_1,0)\to (G_2\bs V_2,D_2,0)$
which is $\CC^\times$-equivariant (and hence extends globally as 
such). Now apply Theorem \ref{thm:duality}     
\end{proof}  

\begin{remark}
The group $G_{i}$ acts on $\Lcal(\Hcal_{i})$ as a group of poset 
automorphisms and we
have a quotient poset $G_{i}\bs \Lcal(\Hcal_{i})$. The ramification 
function induces $\kappa_{i}: G_{i}\bs \Lcal_{\irr}(\Hcal_{i})\to \QQ$. 
If $z_{i}$ is the function on $G_{i}\bs \Lcal_{\irr}(\Hcal_{i})$ 
which assigns to $L\in \Lcal_{\irr}(\Hcal_{i})$
the order of the group of scalars in the image of $Z_{G_i}(L)$
in $V_{i}/L$, then  the
isomorphism $f$ of this theorem induces an isomorphism
of posets $G_1\bs \Lcal(\Hcal_1)\cong G_2\bs \Lcal(\Hcal_2)$
which takes  $z_2$ to $(1-\kappa_1)z_1$ and
$z_1$ to $(1-\kappa_2)z_2$. 
\end{remark}  

\subsection{Dunkl connections with finite holonomy (continued)}\label{subs:finiteholctd}
In this subsection we concentrate on a situation where we want to
establish finite holonomy without the hypothesis that $\kappa_H<1$ for all $H\in\Hcal$. 
We denote the collection of $L\in\Lcal_\irr (\Hcal)$ for which $\kappa_L-1$ is 
negative, zero, positive by  $\Lcal^- (\Hcal)$, $\Lcal^0 (\Hcal)$, $\Lcal^+ (\Hcal)$
respectively. Since $\kappa$ is monotonic, the union $V^{-}$  of the members of 
$\Lcal^{-}(\Hcal)$ is an open subset of $V$. 

The result that we are aiming  at is the following. 
It will be used when we treat the hyperbolic case. 

\begin{theorem}\label{thm:finitehol2}
Let be given a Dunkl system which has a flat positive definite hermitian form. 
Suppose that $\Lcal^0(\Hcal)$ is empty and that the following two 
conditions are satisfied:
\begin{enumerate}
\item[(i)] every $H\in\Hcal$ with $\kappa_L<1$ and every 
line $L$ in $\Lcal_\irr (\Lcal)$ with $\kappa_L>1$ satisfies the  
Schwarz condition and 
\item[(ii)] the intersection of any two distinct members of 
$\Lcal_\irr (\Lcal)$ with $\kappa_L>1$ is irreducible.
\end{enumerate}
Then the system has a finite holonomy group, satisfies the Schwarz condition,
and the developing map induces an isomorphism $G\bs V^-\cong 
\G\bs  A^\circ$, where  $A^\circ$ is a linear arrangement complement in $A$. 
This gives $\PP (G\bs V^-)$ the structure of
an elliptic orbifold whose completion can be identified with $\G\bs \PP (A)$.
\end{theorem}

\begin{remark}
Observe that we are not making the assertion here that the 
developing map extends across a cover of  $V$. 
In fact, if we  projectivize, so that we get a Fubini metric on $\PP (V^\circ)$,
then we will see that the metric completion of $\PP (V^\circ)$ may involve some 
blowing up and blowing down on $\PP (V)$. The modification of $\PP (V)$ 
that is involved here is discussed below in a somewhat more general setting. 
After that we take up the proof of the theorem.
\end{remark}

\begin{discussion}\label{disc:devmap} 
Let be given a Dunkl system with semisimple holonomy around the members of
$\Lcal(\Hcal)$ and for which  $\Lcal^0(\Hcal)$ is empty, but
$\Lcal^+(\Hcal)$ is nonempty (so that $\kappa_0>1$).
We further assume that $\kappa$ takes values in $\QQ$
and that the Schwarz condition is satisfied by all  members
of $\Lcal^-(\Hcal)$ of codimension one (hyperplanes) and 
and all  members of $\Lcal^+(\Hcal)$  of dimension one (lines).
It follows from Theorem \ref{thm:finite2} that the holonomy cover extends to a normal 
cover  $\widetilde{V^{-}}\to V^{-}$ 
and that the  developing map extends to that cover and factors through a local 
isomorphism $G\bs\widetilde{V^{-}}\to A$. 
Let $f: V^+\to V$ be obtained by the blowing up the members of 
$L\in \Lcal^{+}(\Hcal)$ in the order defined by the partial order
(so starting with the origin first).  We shall identify $V^{-}$ with its 
preimage in $V^+$. Notice that the group $G$ naturally acts on  $V^+$.

Every $L\in \Lcal^{+}(\Hcal)$ defines an exceptional divisor $E(L)$ and 
these exceptional divisors intersect normally.
If we write $1-\kappa_L=p_L/q_L$ as usual (so $p_L$ and $q_L$ are relatively prime
integers with $q_L>0$ and hence $p_L<0$), then 
the holonomy around $E(L)$ is of finite order $q_L$.
So the holonomy covering extends to a ramified covering
$\widetilde{V^+}\to V^+$. The preimage of $\cup_L E(L)$ in $\widetilde{V^+}$ is 
also a normal crossing divisor. According to Lemma 
\ref{lemma:infsimple} the  affine structure on $V^\circ$ degenerates simply 
along $E(L)$ with logarithmic exponent $\kappa_L-1$ and the associated
affine foliation is given by its projection onto $L$.

The divisors $E(L)$ determine a simple type of stratification of $V^+$.
Let us describe the strata explicitly. For $L\in \Lcal^{+}(\Hcal)$ we put 
\[
L^{-}:= L- \cup \{M\; :\; M\in\Lcal^{+}(\Hcal) , L<M\}.
\]
So every $M\in\Lcal_\irr(\Hcal)$ which meets $L^{-}$ but does not contain $L$
belongs to $\Lcal^{-}(\Hcal)$. The preimage of  $L^{-}$ in $V^+$ is a union
of strata and trivial as a stratified space over $L^{-}$.
It has a unique open-dense stratum which can be identified with
the product $L^{-}\times \PP((V/L)^{-})$.

An arbitrary stratum is described inductively: 
the  collection of divisors defined by a subset of $\Lcal^{+}(\Hcal)$ has a nonempty
intersection if and only if that subset  makes up a flag:
$L_\pt: L_0>L_1>\cdots >L_k>V$. Their common intersection decomposes as a product:
\[
E({L_\pt}):=L_0^{+}\times \PP((L_1/L_0)^{+})\times\cdots\times \PP((V/L_k)^{+})
\]
and contains a stratum $S(L_\pt)$ as an open-dense subset,
which decomposes accordingly as:
\[
S(L_\pt)=L_0^{-}\times \PP((L_1/L_0)^{-})\times\cdots\times \PP((V/L_k)^{-}).
\]
The developing map will in general not extend to $\widetilde{V^+}$
(it will have a pole along the preimage of $\cup_L E(L)$), but 
things improve if we projectivize. That is why we shall focus on the 
central exceptional divisor $E_0$, which we will also denote 
by $\PP (V^+)$. Notice that $\PP (V^+)$ is a projective manifold
and that $\widetilde{V^+}\to V^+$ restricts to a $\G$-covering 
$\PP(\widetilde{V^+})\to \PP (V^+)$. Each $E(L)$ with $L\in \Lcal^{+} (\Hcal)-\{ 0\}$ meets
$\PP (V^+)$ in a smooth hypersurface $D(L)$ of $\PP (V^+)$ and these 
hypersurfaces intersect normally in  $\PP (V^+)$. The open dense stratum
of $\PP (V^+)$ is clearly $\PP (V^{-})$.
The group $\G$ acts on $\PP(\widetilde{V^+})$
properly discontinuously with compact orbit space $\PP (V^+)$.
We have a projectivized  developing map
\[
G\bs \PP (\widetilde{V^{-}})\to \PP(A)
\]
which is a local isomorphism. A stratum of  $\PP (V^+)$ is given by a flag $L_\pt$ as above with
$L_0=\{ 0\}$ and so will have the form:
\[
S(L_\pt)\cong
\PP(L_1^{-})\times\PP((L_2/L_1)^{-})\cdots \times\PP((V/L_k)^{-}).
\]
It is open-dense in $E(L_\pt)=\PP((L_1/L_0)^{+})\times\cdots\times \PP((V/L_k)^{+})$. 
Let us now write $E_i$ for $E_{L_i}$, $\kappa_i$ for $\kappa_{L_i}$ etc.
According to Proposition \ref{prop:normalformctd}, the developing map is then at 
$z=(z_1,\dots ,z_{k+1})\in S(L_\pt)$ linearly equivalent to a map of the form:
\[
V^+_z\to \prod_{i=1}^{k+1}(\CC\times T_i),\quad 
\Big( t_0^{1-\kappa_0}\cdots t_{i-1}^{1-\kappa_{i-1}}(1, F_i)\Big)_{i=1}^{k+1}.
\]
Here $F_i: V^+_z\to T_i$  is a morphism to a linear space $T_i$ whose restriction to 
$S(L_\pt)_z$ factors as the projection $S(L_\pt)_z\to\PP (L_i/L_{i_1})_{z_i}$ followed
by a local isomorphism $\PP (L_i/L_{i_1})_{z_i}\to T_i$, and 
$t_i$ is a local equation for $E_{L_i}$. So
$(t_{i-1},F_i)_{i=1}^{k+1}$ is a chart for $V^+$ at $z$.
If $\dim L_1=1$, then in terms of this chart, the group $G_{L_1}$ 
acts in the $t_1$-coordinate only (as  multiplication by $|p_1|$th roots of unity).

We restrict the projectivized developing map to $\PP (V^+)$ (which is defined by 
$t_0=0$). The preceding shows 
that this restriction is projectively equivalent to the map with coordinates 
\[
\Big(t_1^{\kappa_{1}-1}\cdots t_{k}^{\kappa_{k}-1}(1,F_1), 
\dots ,t_{k}^{\kappa_{k}-1} (1,F_k), (1, F_{k+1})\Big).
\]
(The component which is constant 1 reminds us that we are mapping to an affine space
which is to be viewed as an open subset of a projective space.) 
Let $\tilde z\in \PP (\widetilde{V^+})$ lie over $z$, put $D_i:=E_i\cap 
\PP(\widetilde{V^+})$
and denote by $\tilde D_i$ the irreducible component of the preimage
of $D_i$ which contains $\tilde z$ and by $\widetilde{S}(L_\pt)$ 
the stratum.  If $i>0$, then near $\tilde z$, $V^+$ is simply given by 
extracting the $q_i$th root of $t_i$: $\tau_i^{q_i}:=t_i$
Since we have semisimple holonomy around the members of $\Lcal(\Hcal)$,
the projectivized developing map is at $\tilde z$  given
in terms of this chart and an affine chart in $\PP(A)$ by
\[
\Big(\tau_1^{-p_1}\cdots \tau_{k}^{-p_k}(1,F_1), 
\dots ,\tau_{k}^{-p_k}(1,F_k), (1, F_{k+1})\Big).
\]
Recall that each $p_i$ is negative. So
this clearly shows that the projectivization defines a regular morphism
$\PP (\widetilde{V^+})\to \PP(A)$ and that its restriction to the preimage of 
$S(L_\pt)$ factors through a covering of the last factor $\PP((V/L_k)^{-})$.
The fiber through $\tilde z$ is here defined by putting  $\tau_k=0$
and $F_{k+1}$ constant. It follows that the connected component of this fiber lies
in $\tilde D_k$, more precisely, that it lies in a connected component of a fiber of the natural map 
$\tilde D_k\to D_k=\PP(L_k^+)\times \PP ((V/L_k)^+)\to \PP ((V/L_k)^+)$.
We also see that $\tilde z$ is isolated in its fiber if and only if the flag is reduced to 
$L_0=\{ 0 \}>L_1$ with $\dim L_1=1$; in that case, the map above is  simply given by 
$(\tau_1^{-p_1},F_1)$.  Since this is also
a chart for the orbit space $G_{L_1}\bs\PP (\widetilde{V^+})$, we see that 
the projectivized developing map modulo $G$ is then a local isomorphism at the image
of $\tilde z$. 
Since the holonomy near $S(L_\pt)$ decomposes as a product, a connected component
$\tilde S(L_\pt)$ of the preimage of $S(L_\pt)$ in $\PP (\widetilde{V^+})$ 
decomposes as a product as well:
$\tilde S(L_\pt)=
\PP(\widetilde{L_1^-})\times\PP(\widetilde{(L_2/L_1)^-})\cdots \times
\PP(\widetilde{(V/L_k)^-})$.
Its closure is an irreducible component of the preimage of $E(L_\pt)$; the normalisation
of that closure decomposes accordingly:
\[
\tilde E(L_\pt)=\PP(\widetilde{L_1^+})\times\PP(\widetilde{(L_2/L_1)^+})
\cdots \times\PP(\widetilde{(V/L_k)^+}).
\]
\end{discussion}

The proof  of \ref{thm:finitehol2} proceeds by induction on $\dim V$. 
The induction starts trivially.

Since the form is positive definite, we shall (by simple averaging) assume that it is
invariant under all the Schwarz symmetry group $G$.

\begin{lemma}\label{lemma:longhol}
For every $L\in \Lcal^{+}(\Hcal)-\{0\}$, the longitudinal holonomy in $L^\circ$ is finite.
\end{lemma}
\begin{proof} 
We verify that the affine structure on $L^\circ$  satisfies the hypotheses
of theorem that we want to prove, so that we can invoke the induction hypothesis.
The flat metric on $V^\circ$ determines one on $L^\circ$. It remains to show 
that every hyperplane  $I\in \Lcal^+(\Hcal^L)$ and every line $M\in \Lcal^-(\Hcal^L)$
satisfies the Schwarz condition. In the first case, $\kappa-1 $ must negative
on $\Hcal_I$ and so it follows from Theorem \ref{thm:devmap2} that $I$ 
satisfies the Schwarz condition.
We claim that in the second case, $M$ is  irreducible in $\Lcal (\Hcal)$ 
(so that the Schwarz condition holds).  For if that were not the case,
then by Lemma \ref{lemma:dich} $M$ has two 
irreducible components, $L$ and $M(L)$. The irreducible component $M(L)$ must be 
in $\Lcal^-(\Hcal)$ by assumption (ii) and since we have $\kappa^L_M=\kappa_{M(L)}$,
we would get a contradiction. 
\end{proof}

\begin{corollary}\label{cor:cfiber}
The connected components of the fibers of the projectivized developing map
$\PP (\widetilde{V^+})\to \PP(A)$ are compact.
\end{corollary}
\begin{proof}
Over $\PP (V^{-})$, the projectivized developing map is locally 
finite and so in these points the claim is clear.  Let us therefore examine
the situation over another stratum $S(L_\pt)$ (as in the Discussion \ref{disc:devmap}).
Since the stratum is not open, we have  $k\ge 1$.
We observed that the connected component of a fiber through $\tilde z$ lies in the fiber
over $z_k\in \PP ((V/L_k)^-)$ of the composite 
\[
\tilde D_k\to D_k=\PP (L_k^+)\times \PP ((V/L_k)^+)\to\PP ((V/L_k)^+)
\]
The holonomy in the last factor $\PP ((L_k)^+)$  is longitudinal and hence finite,
The implies that every irreducible component in $\tilde D_k$ over 
$\PP (L_k^+)\times\{ z_k\}$ is compact.
\end{proof}

A continuous map $f:X\to Y$ between topological spaces always 
has a \emph{topological Stein factorization}: this is the
factorization through the quotient $X\to
X_\stein$ of $X$ defined by  the partition of $X$ into connected 
components of fibers of $f$. So the latter map has 
then connected fibers
and the induced map $f_\stein :X_\stein\to Y$ has discrete fibers
in case the fibers of $f$ are locally connected. 
Here is a useful criterion for an analytic counterpart.

\begin{lemma}\label{lemma:stein}
Let $f: X\to Y$ be a morphism of connected normal analytic spaces. Suppose
that the connected components of the fibers of $f$ are compact. Then 
the Stein factorization of $f$, 
\[
\begin{CD}
f:X @>>> X_\stein @>{f_\stein}>> Y,
\end{CD}
\]
is in the analytic category. More precisely, $X\to X_\stein$ is a proper 
morphism with connected 
fibers to a normal analytic space $X_\stein$ and $f_\stein$ is a 
morphism with discrete fibers. If in addition, $Y$ is smooth, $f$ is a 
local isomorphism in every point that is isolated in its fiber and such 
points are dense
in $X$, then $f_\stein$ is a local isomorphism.
\end{lemma}
\begin{proof}
The first part is well-known and standard in case $f$ is proper. 
The second part perhaps less so, but we show that it is a 
consequence of the first part. Since $f:X\to Y$ is 
then a morphism from a normal analytic space  to a smooth space of the 
same dimension which contracts its singular 
locus, $f_\stein :X_\stein\to Y$ will be a local isomorphism outside
a subvariety of $X_\stein$ of codimension one. But then there is no 
ramification at all, since a ramified cover of a smooth variety has
as its ramification locus a hypersurface. 
 
 So it remains to show that we can reduce to the proper case.
We do this by showing that if $K\subset X$ is a connected component 
of the fiber  $f^{-1}(y)$,
then there exist open neighborhoods $U$ of $K$ in $X$ and $V$ of $y$ in 
$Y$ such that
$f(U)\subset V$ and $f: U\to V$ is proper. This indeed suffices, for if
$y'\in V$, then $f^{-1}(y')\cap U$ is open and closed in $f^{-1}(y')$, 
and hence a union of connected components of $f^{-1}(y')$.

Choose a compact neighborhood $C$ of $K$ which does not meet $f^{-1}(y)-K$.
Clearly, for every neighborhood $V$ of
$y$ in $Y$, $f: f^{-1}V\cap C\to V$ is proper. So it is enough to show that
$f^{-1}V\cap C$ is open in $X$ (equivalently, 
$f^{-1}V\cap \partial C=\emptyset$) for $V$ small enough. 
If that were not the case, then we could find a sequence of points 
$(x_i\in \partial C)_{i=1}^\infty$ whose image sequence converges to $y$.
Since $\partial C$ is compact, a subsequence will converge, to 
$x\in \partial C$, say.
But clearly $f(x)=y$ and so $x\in K$. This cannot be since   
$K\cap \partial C=\emptyset$.
\end{proof}

\begin{corollary}\label{cor:stein}
The Stein factorization  of $ G\bs \PP (\widetilde{V^+})\to \PP(V)$,
\[
\begin{CD}
G\bs \PP (\widetilde{V^+}) @>>> 
(G\bs \PP (\widetilde{V^+})_\stein @>>> \PP (V),
\end{CD}
\]
is analytic and the Stein factor $(G\bs\PP (\widetilde{V^+})_\stein\to\PP (V)$
is a  local isomorphism.
\end{corollary}
\begin{proof}
In Corollary \ref{cor:cfiber} and the Discussion \ref{disc:devmap} we established  that 
the conditions in both clauses of the Lemma \ref{lemma:stein} are satisfied.
\end{proof}

\begin{proof}[Proof of Theorem \ref{thm:finitehol2}]
We first prove that  $ \PP (G\bs\widetilde{V^+})_\stein\to\PP (V)$ is a 
$\G$-isomorphism. For this we verify that the hypotheses of Lemma \ref{lemma:cov} 
are verified for that  map with $Y'=Y=\PP (V)$. By Corollary \ref{cor:stein}
$\PP (G\bs\widetilde{V^+})_\stein\to \PP(V)$ is a local isomorphism.
We know that $\G$ acts properly discontinuously on
$\PP (\widetilde{V^+})$ with compact  fundamental domain. This is then also 
true for $\PP (G\bs\widetilde{V^+})_\stein$. 
Since $\G$ acts on $\PP(V)$ as a group of isometries, Condition (ii)
of \ref{lemma:cov} is fulfilled  as well.
So $\PP (G\bs \widetilde{V^+})_\stein\to\PP (V)$
is a covering projection. But $\PP (V)$ is simply connected, 
and so this must be an isomorphism. It follows that 
$\PP (\widetilde{V^+})$ is compact, so that $\G$ must be finite.

An irreducible component $\tilde D(L)$ over $D(L)$ gets contracted if $\dim L>1$,
with image in $\PP (V)$ a subspace of codimension equal to the dimension of $L$.
In particular, we get a divisor in case $\dim L=1$ and so 
the image of a covering of $\PP (V^-)$ is mapped to an arrangement complement,
$\PP (A^\circ)$, say.
So the developing map $\ev_G: G\bs\widetilde{V^-}\to A^\circ$ becomes an 
isomorphism if we pass to $\CC^\times$-orbit spaces. According to
Theorem \ref{thm:devmap}  $\ev_G$ is homogeneous of degree one. It follows 
that this map as well as the induced map $G\bs V^-\to \G\bs A^\circ$
are isomorphisms.

Finally we verify the Schwarz condition for any $L\in\Lcal_\irr(\Hcal)$.
We know already that this is the case when $L\in\Lcal^-(\Hcal)$.
For $L\in\Lcal^+(\Hcal)$ this is seen from
the simple form of the projectivized developing map at a general point
of $D(L)$: in terms of  a local chart $(\tau_1,F_1,F_2)$ of 
$\PP (\widetilde{V^+})$
at such a point it is given by $(\tau_1^{-p_1},\tau_1^{-p_1}F_1, F_2)$.
Since $(G\bs \PP (\widetilde{V^+})_\stein\to \PP(V)$ is an isomorphism,
$G$ must contain the group of $|p_L|$th roots of unity
acting on the transversal coordinate $\tau_1$. This just tells us that $L$ 
satisfies the  Schwarz condition.
\end{proof}

\subsection{Dunkl connections whose holonomy is almost a Heisenberg group}
\label{sect:crystal}

\begin{theorem}\label{thm:euclidean}
Let be given a Dunkl system with $\kappa\in (0,1)^\Hcal$ and $\kappa_0=1$, 
which satisfies the  Schwarz condition in codimension one 
and admits a nontrivial flat hermitian form. Then:
\begin{enumerate}
\item[(i)]  the flat hermitian form is semidefinite with kernel generated by the Euler field,
\item[(ii)]  $V^f=V-\{ 0\}$, the monodromy group $\G /\G_0$ of the
connection on $G\cdot\CC^\times\bs V^\circ$ is finite and 
$\G_0$ is an integral Heisenberg group,
\item[(iii)]  the developing map identifies the  $\G /\G_0$-cover of 
$G\bs V-\{ 0\}$ in a $\CC^\times$-equivariant fashion  with 
an anti-ample $\CC^\times$-bundle over an abelian variety,
\item[(iv)] $G\bs \widetilde{V-\{ 0\}}\to A$ is a $\G$-isomorphism and 
the Dunkl connection satisfies the  Schwarz condition.
\item[(v)] The hermitian form gives $\PP (G\bs V)$ the structure
of a complete parabolic orbifold: if $K$ is the kernel of the hermitian form on the 
translation space of $A$, then $\G$ acts in $K\bs A$ via 
a complex crystallographic space group and the developing map induces an isomorphism
between $\PP(G\bs V)$ and the latter's orbit space.
\end{enumerate}
\end{theorem}
\begin{proof}
The first assertion follows from Theorem \ref{thm:semipos}. 
Upon replacing  the flat form by its negative, we assume that 
it is positive semidefinite; we denote this form by $h$.
 The monodromy around every member of  $\Lcal_{\irr}(\Hcal)-\{ 0\}$
leaves invariant a positive definite form 
and hence is finite by Theorem \ref{thm:finite2}. 
 This implies that $V^f\supset V-\{ 0\}$; it also shows that the 
 monodromy of the connection is finite.
 Since $\kappa_0=1$, the Euler field $E_V$ is flat and determines
a nonzero translation $T$ $A$ such that $2\pi \sqrt{-1}T$ is the monodromy
around a $\CC^\times$-orbit in $V^\circ$. In particular, the monodromy
around such an orbit is not of finite order, so that $V^f= V-\{ 0\}$.

The Euler field resp.\ $T$ generate  a faithful $\CC^+$-actions on 
$\widetilde{V-\{ 0\}}$ resp.\ $A$ such that the developing map descends to a 
local isomorphism 
$(\CC^+\cdot G\bs \widetilde{V-\{ 0\}}\to \CC^{+}\bs A$.
Observe that the translation space of $\CC^+\bs A$ has a $\G$-invariant positive 
definite  hermitian form: if the kernel of
$h$ is spanned by $E_V$ this is clear and if $h$ is positive
definite we simply identify the translation space in question with the
orthogonal complement of $T$ in the translation space of $A$.
The group $\G/(2\pi \sqrt{-1}T)$ acts on $\CC^{+}.G\bs \widetilde{V-\{ 0\}}$ 
through  a group which acts properly  discretely. The orbit space of this action 
can be identified with $G\bs \PP (V)$, hence is compact. So the assumptions of Lemma 
\ref{lemma:cov} are fulfilled (with $Y'=Y=\CC^+\bs A$) and we conclude that 
\[
\CC^+.G\bs \widetilde{V-\{ 0\}}\to\CC^+\bs A
\]
is a covering. Since the range is an affine space (hence simply  connected), this must 
be an isomorphism. 
It follows that the action of $\G $ on $A$ is properly discrete and cocompact.
It also follows that the developing map defines  a $\G$-equivariant 
isomorphism of  $G\bs \widetilde{V-\{ 0\}}$ onto $A$.

Let $\G_{0}$ be the subgroup of $\g\in\G$ that act as a translation in 
$\CC^+\bs A$. This subgroup is of finite index in $\G$ and our assumption implies that 
$\G_{0}\bs A\to \G_{0}\cdot\CC^{+}\bs A$ has the structure 
of a flat $\CC^{\times}$-bundle over a complex torus. The developing map induces an 
isomorphism $\G_{0}\bs A\cong \G_{0}\cdot G\bs \widetilde{V-\{ 0\}}$;
the latter is finite over $G\bs V-\{ 0\}$ and extends therefore as a finite
cover over $G\bs V$. This means that the associated line bundle over the complex
torus has contractible zero section. Hence this line bundle is anti-ample and
$\G_{0}$ is a Heisenberg group. 

Property (iv) is almost immediate from Theorem \ref{thm:devmap2}.
\end{proof}

\section{Geometric structures of hyperbolic type}

In this section we consider Dunkl systems of admissible hyperbolic type. 
So the affine space $A$ in which the evaluation map takes its values 
is in fact a vector space (it comes with an origin) equipped with a 
nondegenerate hermitian form of hyperbolic signature. We denote by
$\LL^\times\subset A$ the set of vectors of negative self-product
and by $\BB:=\PP (\LL^\times)\subset\PP (A)$ its projectivization.
Notice that $\BB$ is a complex ball and that $\LL^\times$ can be thought
of as a $\CC^\times$-bundle over $\BB$. By adding $\BB$ at infinity
we obtain a line bundle $\LL$ over $\BB$ that has $\BB$ as the zero section.
The admissibility assumption means that the evaluation map takes its values
in $\LL^\times$ so that its projectivization takes its values in $\BB$.

\subsection{The compact hyperbolic case}
This is relatively simple case and for that reason we state and prove it 
separately. The result in question is the following.

\begin{theorem}\label{thm:comphyp}
Suppose that the Dunkl system is of admissible hyperbolic type,
satisfies the Schwarz condition in codimension one and is such that
$\kappa\in (0,1)^\Hcal$, $\kappa_L<1$ for all $L\in\Lcal_\irr(\Hcal)-\{ 0\}$. 
Then the Dunkl system satisfies the Schwarz condition,
$\G$ acts on $\BB$ discretely and with compact fundamental domain
and the developing map induces an isomorphism $G\bs V\cong \G\bs \LL^\times$. 
Thus $\PP (G\bs V)$ acquires the structure of a complete hyperbolic orbifold
isomorphic to $\G\bs \BB$.
\end{theorem}
\begin{proof}
Arguing as in the proof of Theorem \ref{thm:euclidean} we find that  
$V^{f}=V-\{ 0\}$. It follows from Theorem \ref{thm:devmap2} that 
the Dunkl system satisfies the Schwarz condition.
The developing map descends to a local isomorphism 
$G\bs \PP(\widetilde{V^f})\to \PP (A)$. It takes values in the 
complex ball $\BB$. The latter comes with a $\G$-invariant K\"ahler metric.
The orbit space of the $\G$-action on $G\bs \PP(\widetilde{V^f})$ can be 
identified with $G\bs \PP (V)$, hence is compact. 
So the assumptions of Lemma 
\ref{lemma:cov} are fulfilled and we conclude that   
$G\bs \PP (\widetilde{V^f})\to\BB$ is a covering. Since the range is simply
connected, this must be an isomorphism. In particular,  the action of
$\G $ on $\BB$ is properly discrete and cocompact.

It also follows that $G\bs V\cong \G\bs \LL^\times$ becomes an isomorphism
if we pass to $\CC^\times$-orbit spaces. It then follows that the map itself
is an isomorphism, because $G$ contains by definition all the scalars which leave 
the developing map invariant.
\end{proof} 

\subsection{Statement of the main theorem}\label{subsect:main}
The general hyperbolic case concerns the situation where the holonomy
group is of cofinite volume (rather than being cocompact) in the
automorphism group of a complex ball. This is substantially harder 
to deal with.
 
Given a Dunkl system for which the flat hermitian form
$h=h^\kappa$ is of hyperbolic type (i.e., nondegenerate of index one, so that $h$ defines
a complex ball $\BB$ in the projective space at infinity $\PP (A)$ of $A$).
If $L\in \Lcal_\irr (\Hcal)$ is such that $\kappa_L>1$, then
if we approach $L^\circ$ from $V^\circ$ along a curve, the image of a lift in 
$\widetilde{V^\circ}$
of this curve under the developing map tends to infinity with
limit a point of $\PP(A)$. These limit points lie in well-defined  $\G$-orbit of
linear subspaces of $\PP(A)$ of codimension $\dim (L)$.
We call such space  a \emph{special subspace} in $\PP (A)$ 
and its intersection with $\BB$ a \emph{special subball}.
We use the same terminilogy for the corresponding linear subspace of $A$.

The main goal of this section is to prove:

\begin{theorem}\label{thm:main}
Let be given a Dunkl system with $\kappa\in (0,1)^\Hcal$ which comes with a flat 
admissible form $h$ of hyperbolic type.  
Suppose that every hyperplane $H\in\Hcal$ with $\kappa_H<1$ and every line 
$L\in \Lcal_\irr(\Hcal)$ with $\kappa_L>1$ satisfies the Schwarz condition. 
Then:
\begin{enumerate}
\item[(i)] The system satisfies the Schwarz condition.
\item[(ii)] The collection of special hyperplanes is locally finite in $\LL^\times$ and
if $(\LL^\times)^-$ denotes the complement in $\LL^\times$ of the 
union of the special hyperplanes, then the projectivized developing map 
defines a $\G$-equivariant isomorphism $G\bs \widetilde{V^f}\to (\LL^\times)^-$.
\item[(iii)] The group $\G$, considered  as a subgroup of the unitary group $\U (h)$ 
of $h$, is discrete and has cofinite volume in $\U (h)$.
\item[(iv)] The developing map induces an isomorphism 
$G\bs V^f\to \G\bs (\LL^\times)^-$ of normal analytic spaces. 
\end{enumerate}
Thus if $\BB^-$ denotes
the complement in $\BB$ of the union of the special hyperplanes, then
$\PP (G\bs V^f)$ can be identified with $\G\bs \BB^-$ and
acquires the structure of a hyperbolic orbifold whose completion is $\G\bs \BB$.
\end{theorem}

\begin{remarks}
Our proof yields more precise information, for it tells us how $\PP(G\bs V)$
is obtained from the Baily-Borel compactification of $\G\bs \BB$ by a blowup 
followed by a blowdown. This is in fact an instance of the construction 
described in \cite{looij:ball}.

Couwenberg gives in his thesis \cite{couw} a (presumably complete) 
list of the cases for which  $\Hcal$ its Coxeter arrangement
and $G$ is the associated Coxeter group. The Schwarz condition
for the lines then amounts to: if 
$L$ is a line which is the fixed point subspace of an irreducible Coxeter subgroup
of $G$ and such that $\kappa_L>1$, then $(\kappa_L-1)^{-1}$
is an integer or, when $L^\perp\in\Hcal$, half an integer.
The fact that the list is substantial gives the theorem its merit.
In particular, it produces new examples of discrete hyperbolic groups of 
cofinite volume. 
\end{remarks}

\subsection{Connection with the work of Deligne-Mostow}
Theorem \ref{thm:main} implies one of the main results of Deligne-Mostow 
\cite{delmost1} and Mostow \cite{mostow:int}. 

\begin{theorem}[Deligne-Mostow]\label{thm:maindm}
Consider the Lauricella system with all of its  parameters 
$\mu_0,\dots ,\mu_n$  in $(0,1)$  and $\sum_{k=0}^n \mu_k\in (1,2)$ so that
$\mu_{n+1}:=2-\sum_{k=0}^n \mu_k\in (0,1)$ also. Suppose that for every
pair $0\le i<j\le n+1$ for which $1-\mu_i-\mu_j$ is positive,  $1-\mu_i-\mu_j$ is a 
rational number with numerator $1$ or $2$, allowing the latter only in case 
$j\le n$ and $\mu_i=\mu_j$.
 Then the system satisfies the Schwarz condition and the Schwarz symmetry group
is the group $G$ of permutations of $\{ 0,\dots ,n\}$ which preserves the weight 
function $\mu: \{ 0,1,\dots ,n\}\to \RR$, the collection of special hyperplanes is 
locally finite on $\BB$, $\G$ is a lattice in the unitary group of $A$ and the developing  
map identifies  $\PP (G\bs V^f)$ with $\G\bs \BB^-$. 
\end{theorem}
\begin{proof} We verify the hypotheses of Theorem \ref{thm:main}.
First of we all we want the Schwarz condition for every $H_{i,j}$ satisfied:  this means that
for every pair  $0\le i<j\le n$, we want $1-\mu_i-\mu_j$ to be positive rational number
with numerator $1$ or $2$, allowing the latter only in case $\mu_i=\mu_j$. 
We also want the Schwarz condition fulfilled at a line in $\Lcal^+(\Hcal)$.
Such a line is given by an $n$-element subset
of $\{ 0,\dots ,n\}$, say as the complement of the singleton $\{ i\}$, 
such that $\sum_{0\le j\le n, j\not= i} \mu_j>1$. The Schwarz condition is fulfilled at 
this line if $-1+\sum_{0\le j\le n, j\not= i} \mu_j$ is the reciprocal of an integer.
This amounts to: if $1-\mu_i-\mu_{n+1}$ is positive, then it is the reciprocal
of an integer. The rest follows from easily from Theorem \ref{thm:main}.
\end{proof}

\begin{remark}
The conditions imposed here imply Mostow's $\Sigma$INT-condition: 
this is the condition which says that for any pair 
$0\le i<j\le n+1$ such that $1-\mu_i-\mu_j$ is positive, we want this to be a 
rational number with numerator $1$ or $2$, allowing the latter only in case $\mu_i=\mu_j$. 
Clearly, this condition is more symmetric, because it does not attribute a special role
to $\mu_{n+1}$.  This symmetry is understood as follows.
We can regard of $\PP (V^\circ)$ as parametrizing
the collection of \emph{mutually distinct} $(n+1)$-tuples $(z_0,\dots ,z_n)$ in 
the affine line $\CC$ given up to an affine-linear
transformation.  But it is better to include $\infty$ and to think of 
$\PP (V^\circ)$ as the  moduli space of
\emph{mutually distinct} $(n+2)$-tuples $(z_0,\dots ,z_{n+1})$ on the  projective line 
$\PP^1$  given up to a projective-linear transformation, that is, to identify 
$\PP (V^\circ)$ with $\Mcal_{0,n+2}$.
This makes evident an action of the permutation group of $\{0,\dots ,n+1\}$ on
$\PP (V^\circ)$. 
It is conceivable that there are cases for which the $\Sigma$INT-condition
is satisfied and ours aren't, even after permutation.  
The table in \cite{thurstontable}, lists 94 systems
$(\mu_0\ge \mu_1\ge \cdots \ge\mu_{n+1}>0)$
satisfying the $\Sigma$INT-condition. Most likely, it is complete.
In this list, there is precisely one case
which escapes us and that is when $n+2=12$ and all $\mu_i$'s equal to 
$\frac{1}{6}$. With little extra effort, we can get around this 
(and at the same time avoid resorting to this list) if we let the group of 
permutations of $\{0,\dots ,n+1\}$ which leave $\mu: \{0,\dots ,n+1\}\to\QQ$ 
invariant act from the outset. This group contains $G$ and the elements not in $G$
act nonlinearly on $\PP (V^f)$. An alternative approach starts with analyzing
the developing map of a Dunkl system with a degenerate hyperbolic form 
(see Subsection \ref{subsect:degenhyp}), which  indeed is
a class worth studying its own right. 
\end{remark}

\begin{remark}\label{rem:delmost}
Deligne and Mostow show that there is a modular interpretation of the Baily-Borel 
compactification of $\G\bs \BB$.  Given positive rational numbers 
$\mu_0,\dots ,\mu_{n+1}$ with sum $2$, then let us say that an \emph{effective 
fractional anticanonical divisor} on $\PP^1$ of 
type $\mu$ is simply a given by a set of $n+2$ points endowed with the weights 
$\mu_0,\dots ,\mu_{n+1}$, given up to order.
We do not require the points to be distinct. 
So such a divisor  determines a \emph{support function} 
$\PP^1\to \QQ_+$  which is zero for all but finitely many points and whose 
sum  (over $\PP^1$) of its values is two.
It is said to be \emph{stable} (resp.\ \emph{semistable}) if this function is 
everywhere less than (resp.\ at most) one. 
The projective linear group acts on the variety of the semistable fractional divisors and 
this action is proper on the  (open) subvariety of the stable ones. So a stable orbit is 
always closed. Any other
minimal semistable orbit is represented by a fractional divisor whose support consists of 
two distinct points, each with weight $1$. The points of its Hilbert-Mumford quotient
are in bijective correspondence with the minimal semistable orbits. 
We thus get a projective  compactification $\Mcal_{0,n+2}\subset
\overline{\Mcal}^\mu_{0,n+2}$.
A period map enters the picture by imitating the familiar approach to the
elliptic integral, that is, by passing to a cyclic cover of $\PP^1$
on which the Lauricella integrand becomes a regular differential. Concretely,
write $\mu_i=m_i/m$ with $m_i,m$ positive integers such that the $m_i$'s
have no common divisor, and  write $\nu_i$ for the denominator of $\mu_i$.
Consider the cyclic cover $C\to \PP^1$ of order 
$m$ which has ramification over $z_i$ of order $\nu_i$. 
In affine coordinates, $C$ is given as the normalization of the curve defined by
\[
w^m=\prod_{i=0}^n (z_i-\zeta)^{m_i}.
\]
The Lauricella integrand pulls back to a regular differential $\tilde\eta$ on $C$, 
represented by $w^{-1}d\zeta$. Over $z_i\in\PP^1$ we have $m/\nu_i$ distinct
points in each of which $\tilde\eta$ has a zero of order $\nu_i (1-\mu_i)-1$.
This form transforms under the Galois group by a certain character $\chi$ and
up to a scalar factor, $\tilde\eta$ is the only regular form with that property:
$H^{1,0}(C)^\chi$ is a line spanned by $\tilde\eta$. It turns out that such Hodge data 
are uniformized by a complex ball. Although
the holonomy group need not map to an arithmetic group,  much of 
Shimura's theory applies here. Indeed, Shimura (see for instance \cite{shimura})
and Casselman \cite{casselman} (who was Shimura's student at the time) had investigated 
in detail the case for which $m$ is prime before Deligne and Mostow 
addressed the general situation.
A (if not \emph{the}) chief result of Deligne-Mostow \cite{delmost1} is 
a refined Torelli theorem: if their INT condition is satisfied, then
\begin{enumerate}
\item[(i)] the holonomy group maps to a subgroup of 
automorphisms of the Hodge period ball which is discrete and of cofinite volume, 
\item[(ii)] the corresponding orbit space admits a compactification of
 Baily-Borel type (this adds a finite number of points, the cusps),  
\item[(iii)] the map described above identifies $\overline{\Mcal}^\mu_{0,n+2}$
with this Baily-Borel compactification, making the  minimal 
semistable nonstable orbits correspond to the cusps. 
\end{enumerate}
This is essentially the content of their Theorem (10.18.2).
They also determine when the holonomy group is arithmetic 
(the systematic construction of such groups was in fact Mostow's original motivation).
\end{remark}

\subsection{The Borel-Serre extension}\label{subsect:borelserre}
Before we begin the prof the main theorem, we  first make a few observations
regarding the unitary group $\U (h)$ of $h$ (since $A$ has an origin, we
regard this as a group operating in $A$). Suppose we have a unipotent 
transformation $g\in \U (h)$ that is not the identity. Let $E\subset A$
be the fixed point space of $g$. Then $E^\perp$ is $g$-invariant and hence
contains eigenvectors. So $E\cap E^\perp$ is non trivial. In other words,
$E$ contains an isotropic line $I$. Now $g$ induces in $I^\perp/I$
a transformation that will preserve the form induced by $h$. Since
this form is positive definite and $g$ is unipotent, $g$ will
act trivially on $I^\perp/I$. 
The unitary transformations which respect the flag
$\{ 0\}\subset I\subset I^\perp\subset A$ and act trivially on the 
successive quotients form a Heisenberg group $N_I$ whose center
is parametrized as follows. Notice that the one-dimensional 
complex vector space $I\otimes \overline{I}$
has a natural real structure which is oriented: it is defined by
the `positive' ray of the elements $e\otimes e$, where $e$ runs 
over the generators of $I$. This line parametrizes a one parameter 
subgroup of $\GL (A)$: 
\[
\exp: I\otimes \overline{I}\to \GL (A), \quad 
\exp (\lambda e\otimes e): z\in A\mapsto z+ \lambda h(z,e)e, 
\quad e\in I, \lambda\in\CC.
\]
The transformation $\exp (\lambda e\otimes e)$ is unitary relative to $h$ 
if $\lambda$ is purely imaginary and so $\exp$ maps
$\sqrt{-1}I\otimes \overline{I}(\RR)$ to a one-parameter subgroup of 
$\U (h)$. This one-parameter subgroup is the center of the 
Heisenberg group $N_I$ above. The group $N_I$ is parametrized by pairs
$(a,e)\in I^\perp\times I$: any element of this group is written
\[
g_{a,e}: z\in A\mapsto z+ h(z,a)e -h(z,e)a-\tfrac{1}{2}h(a,a)h(z,e)e.
\] 
This is not quite unique since $g_{a+\lambda e,e}=g_{a,e}$ when 
$\lambda\in\RR$. But apart from that we have uniqueness: $N_I$ 
modulo its center can be identified with vector group 
$I^\perp/I\otimes \overline{I}$ by 
assigning to $(a,e)$ its image in $I^\perp/I\otimes \overline{I}$.

Let $T$ be a subspace of $A$ on which $h$ is degenerate with kernel $I$:
so $I\subset T\subset I^\perp$. We suppose that $T\not= I$. 
Clearly, $N_I$ preserves $T$.
Suppose that $g$ acts trivially on $A/T$ and induces in the
fibers of $A\to A/T$ a translation. So if we write $g$ in the above 
form: $g=g_{a,e}$, then we see that $a$ must be proportional to $e$:
$a=\lambda e$ with $\lambda$ purely imaginary, in other words $g$ is in 
the center of $N_I$.  

Let $I\subset A$ be an isotropic line. When $\lambda$ is a positive
real number, and $e\in I$, then 
$\exp (\lambda e\otimes e)$ is not unitary, but it will still map 
$\BB$ into itself.  
In fact, the orbits of the ray of positive elements in 
$I\otimes \overline{I}$
are (oriented) geodesic rays in $\BB$ which tend to $[I]\in \partial\BB$. 
Perhaps a more concrete picture is gotten by fixing a generator $e\in I$
so that every point of $\BB$ can be represented in the affine hyperplane
in $V$ defined by $h(z,e)=1$: under the realization of $\BB$ in this
hyperplane, the geodesic ray action becomes simply 
the group of translations over positive multiples of $e$.
We regard the space $\BB (I)$ of these rays as a quotient space of $\BB$
so that we have a fibration by rays $\pi (I): \BB\to \BB (I)$.
The Borel-Serre topology on the disjoint union $\BB\sqcup \BB (I)$
is generated by the open subsets of $\BB$ and the subsets
of the form $U\sqcup \pi (I)(U)$, where $U$ runs over the open subsets
of $\BB$ invariant under $N_I$ and the positive ray in 
$I\otimes \overline{I}$. This adds a partial boundary to $\BB$ so that
it becomes a manifold with boundary.
Let $\BB^+\supset \BB$ be the \emph{Borel-Serre extension} associated 
to $\G$: for every isotropic line $I\subset V$ for which
$\G\cap N_I$ is discrete and cocompact,
we do the above construction. That makes $\BB^+$ a manifold with boundary,
the boundary having in an infinite number of connected components (or being
empty). Notice that the action of $\G$ on this boundary is
properly discrete and cocompact---this is indeed the main justification
for its introduction.

\subsection{Proof of the main theorem}
We now turn to the proof of Theorem \ref{thm:main}. 
Throughout this section the assumptions of that theorem are in force and
we also retain some of the notation introduced in Subsection \ref{subs:finiteholctd},
such as $\Lcal^-(\Hcal),\Lcal^0(\Hcal),\Lcal^+(\Hcal),\cdots $.

We begin with a lemma in which we collect a number of useful properties.

\begin{lemma}\label{lemma:collect}
We have:
\begin{enumerate}
\item[(i)] For any $L\in\Lcal_\irr(\Hcal)$, $h$ induces on 
$(V/L)^\circ$ a flat hermitian form which is positive, semipositive with 
one-dimensional kernel, hyperbolic according to whether $\kappa_L-1$ is negative,
zero, or positive. 
\item[(ii)] The intersection of any two distinct members $L_1,L_2$ of
$\Lcal^0(\Hcal)\cup \Lcal^+(\Hcal)$ is irreducible and (hence) belongs
to $\Lcal^+(\Hcal)$.
\item[(iii)] If $L\in\Lcal^+(\Hcal)$, then the longitudinal Dunkl connection on 
$L^\circ$ has finite holonomy and $L$ satisfies the Schwarz condition (so that
the system satisfies the Schwarz condition).
\end{enumerate}
\end{lemma}
\begin{proof}
The flat hermitian form induces one on the Dunkl system $V/L$. This form is 
nonzero ($L$ cannot be a hyperplane since we assumed that $\kappa$ takes a value 
less than one on these) and so the first statement readily follows from our results
in Section \ref{sect:elliptic}. 

If for $L_1,L_2$ as in the lemma, $L_1\cap L_2$ were reducible, then 
then the flat form on $V/(L_1\cap L_2)$   induced by $h$ would have an 
isotropic plane,  a property which is clearly forbidden by the  signature of $h$.

Let now $L\in\Lcal^+(\Hcal)$. Then the longitudinal holonomy
in $L^\circ$ has a flat positive hermitian form. The desired properties now follow
from Theorem \ref{thm:finitehol2}: in view of the way $\kappa^L$ is defined, 
and part (ii) any one-dimensional member in $\Lcal^+(\Hcal^L)$ is 
in fact a member of $\Lcal^+(\Hcal)$ and so satisfies the Schwarz condition
and any codimension one member in $\Lcal^-(\Hcal^L)$ comes from
a member of $\Lcal^-(\Hcal)$ and hence satisfies the Schwarz condition. 
\end{proof}

\begin{discussion}\label{disc:devmap2}
We introduced in the Discussion \ref{disc:devmap} a blowup $V^+$ under the assumption that
$\Lcal^0(\Hcal)$ is empty and described the behavior of  the projectivized 
developing map on the preimage of the origin of $V$. We generalize this
to the situation where $\Lcal^0(\Hcal)$ is allowed to be nonempty. 

Our $V^+$ will now be obtained by blowing up the members of $\Lcal^+(\Hcal)$
first (in the usual order), and then blowing up each $L\in \Lcal^0(\Hcal)$ in a 
\emph{real-oriented manner}. This is unambiguously defined since 
by Lemma \ref{lemma:collect}-(ii) the intersection
of two such members lies in $\Lcal^+(\Hcal)$ and so their strict 
transforms will not meet. It is clear that $V^+$ is a manifold with smooth
boundary whose manifold interior $V^+ -\partial V^+$ is a quasiprojective 
variety. The latter contains $V^f$ as an open-dense subset
and the complement of $V^f$ in $V^+ -\partial V^+$ is a normal crossing divisor
whose closure in $V^+$  meets the boundary transversally. 

Any $L\in\Lcal^+(\Hcal)$ defines a divisor $E(L)$ in $V^+$
and any $L\in\Lcal^0(\Hcal)$ defines a boundary component
$\partial_LV^+$. These cross normally in an obvious sense so that we get 
a natural stratification of $V^+$.
Let us describe the strata explicitly. For $L\in \Lcal^0(\Hcal)\cup\Lcal^+(\Hcal)$ 
we define $L^-$ as in Discussion \ref{disc:devmap}:
\[
L^-:= L- \cup \{M\; :\; M\in \Lcal^0(\Hcal)\cup \Lcal^+(\Hcal) , L<M\}.
\]
So every $M\in\Lcal_\irr(\Hcal)$ which meets $L^-$ but does not contain $L$
belongs to $\Lcal^-(\Hcal)$.
In particular, $L^-$ is contained in the subset $L^f$ of $L$ defined
by the longitudinal connection. It is clear that $V^-=V^f$.
The preimage of  $L^-$ in $V^+$ is a union
of strata and trivial as a stratified space over $L^-$.
It has a unique open-dense stratum which can be identified with
the product $L^-\times \PP((V/L)^f)$ in case
$L\in\Lcal^+(\Hcal)$. If $L\in\Lcal^0(\Hcal)$, then we must replace 
the factor $\PP((V/L)^f)$ by $\SS(V/L)$, where 
$\SS$ assigns to a (real) vector space the sphere of its real half lines. 
(There is no need to write $(V/L)^f$ here, since the latter equals $V/L-\{ 0\}$.)

An arbitrary stratum is  described inductively: 
the  collection of divisors and boundary walls
defined by a subset of $\Lcal^0(\Hcal)\cup \Lcal^+(\Hcal)$ has a nonempty
intersection if and only if that subset  makes up a flag:
$L_\pt: L_0>L_1>\cdots >L_k>L_{k+1}=V$. 
Their common intersection contains a stratum $S(L_\pt)$ which decomposes as 
\[
S(L_\pt)=L_0^-\times \prod_{i=1}^k\PP((L_i/L_{i-1})^-)\times\PP((V/L_k)^f),
\]
at least, when $L_k\in\Lcal^+(\Hcal)$; if $L_k\in\Lcal^0(\Hcal)$, we must replace 
the last factor by $\SS (V/L_k)$. It is clear that $G .\CC^\times$
naturally acts on $V^+$. The covering $\widetilde{V^f}\to V^f$ extends naturally 
to a ramified covering $\widetilde{V^+}\to V^+$ with $\G\times G$-action. 
Since the holonomy along $S(L_\pt)$ decomposes according to its factors,
a connected component $\tilde{S}(L_\pt)$ of the
preimage of a stratum $S(L_\pt)$ decomposes as a product of coverings
of the factors of $S(L_\pt)$. By Lemma \ref{lemma:collect}, the
covers of these factors are finite except for the
last, which is the holonomy cover of $\PP ((V/L_n)^f)$ or $\SS(V/L_n)$.
 
The preimage $\PP (V^+)$ of the origin of $V$ in $V^+$ is a compact manifold with 
boundary. Let us write $B^+$ for $\PP (V^+)$ and denote its interior 
by $B$. So $B$ is a quasiprojective manifold which contains 
$\PP(V^f)$ as the complement of a normal crossing divisor. The strata in $B^+$ 
are given by the flags $L_\pt$ which begin with $L_0=\{0\}$. We  denote by
$D(L)$ the exceptional divisor in $B^+$ defined by $L\in\Lcal^+(\Hcal)$. (It is easy 
to see that $D(L)=\PP (L^+)\times \PP ((V/L)^+)$.) The group $\G$ acts on $\widetilde{B^+}$
properly discontinuously with compact orbit space $B^+$.
\end{discussion}

\begin{proposition}\label{prop:collect}
The projectivized developing map
extends to this covering as a continuous $\G$-equivariant map 
$\widetilde{B^+}\to \BB^+$ which is constant on the
$G$-orbits. It has the following properties:
\begin{enumerate}
\item[(i)] It maps every boundary component of
$\widetilde{B^+}$ to a Borel-Serre boundary component of $\BB^+$ and the 
restriction $\widetilde{B}\to \BB$ is analytic.
\item[(ii)] Every irreducible component of the preimage in $\widetilde{B}$ 
of an exceptional divisor $D(L)$, $L\in\Lcal^+(\Hcal)$, 
is mapped to an open subset of special subball of $\BB$ of codimension 
$\dim (L)$ and the resulting map from such irreducible components to special subballs 
reverses the inclusion relation.
\item[(iii)] Every connected component of a fiber of the map
$\widetilde{B^+}\to \BB^+$ is compact. If that connected component  is
a singleton, then at the image of this singleton in $G\bs\widetilde{B^+}$, the map
$G\bs\widetilde{B^+}\to \BB^+$ is local isomorphism.
\end{enumerate}
\end{proposition} 
\begin{proof}
The proof amounts to an analysis of the behavior of the projectivized developing map
on $\widetilde{B^+}$. Since we did this already in the case without boundary
components in the proof of Theorem \ref{thm:finitehol2}, we shall now concentrate 
on the case of a boundary stratum.  Such a stratum is given by a flag 
$L_\pt=(\{ 0\}=L_0>L_1>\cdots >L_k>L_{k+1}=V)$, for which
$L_i\in \Lcal^+ (\Hcal)$ for $i<k$ and $L_k\in \Lcal^0 (\Hcal)$:
\[
S(L_\pt)\cong
\PP((L_1/L_0)^-)\times\cdots \times\PP((L_k/L_{k-1})^-)\times\SS (V/L_k) 
\]
Let us write $\partial_k$
for the boundary component of $B^+$ defined by $L_k$.
If we had not blown up the strict transform of $L_k$ in a real-oriented fashion,
but in the conventional manner, then the last factor would be
$\PP (V/L_k)$. On a point over that stratum, the developing map is 
according to Proposition \ref{prop:normalformctd} affine-linearly equivalent to a 
map taking values in  $\CC\times T_1\times\CC\cdots \times T_k\times\CC$
with components
\[
\Big( (t_0^{1-\kappa_0}\cdots t_{i-1}^{1-\kappa_{i-1}}(1, F_i))_{i=1}^{k-1}, 
t_0^{1-\kappa_0}\cdots t_{k-1}^{1-\kappa_{k-1}}(1,F_k,\log t_{k})\Big).
\]
Here $F_i$ is a morphism at a point of this conventional blowup to a linear space $T_i$,
$t_i$ defines the $i$th exceptional divisor and 
$(t_0,F_1,\dots ,F_k,t_k)$ is a chart. However, on the real-oriented blowup,
$\log t_{k}$ is a coordinate: its imaginary part $\arg t_{k}$ helps to
parametrize the ray space $\SS (V/L_k)$ and its real part $\log |t_{k}|$
must be allowed to take the value $-\infty$ (its value on the boundary). 
We denote this coordinate $\tau_{k}$.
On a connected component $\tilde S(L_\pt)$ of the preimage of $S(L_\pt)$ in
$\widetilde{B^+}$, we have defined roots of the normal coordinates: 
$t_i=\tau_i^{q_i}$, $i=0,\dots ,k-1$, so that 
$(F_1,\tau_1\dots ,F_k,\tau_k)$ is a chart for $\widetilde{B^+}$.  In terms of this chart, 
the  projectivized developing map becomes
\[
\Big(\tau_1^{-p_1}\cdots \tau_{k-1}^{-p_{k-1}}(1,F_1), 
\dots , \tau_{k-1}^{-p_{k-1}} (1,F_{k-1}), (1, F_k,\tau_k)
\Big),
\]
where we recall that $-p_i$ is a positive integer and the constant 
component $1$ reminds us of the fact that we are mapping to an affine
chart of a projective space. 
We use this to see that the  projectivized developing map extends to 
$\widetilde{B^+}\to \BB^+$. 
A chart of $\BB^+$ is (implicitly) given by the affine hyperplane $A_1\subset A$ defined by 
$h(-,e)=1$, where $e$ is minus the unit vector corresponding to the slot occupied by 
$\tau_{k}$ (the geodesic action is then given by 
translation over negative multiples of $e$). 
This  normalization is here already in place, for the coordinate in question is 
in the slot with constant 1. So we then have in fact a chart 
of the Borel-Serre compactification, provided that we remember that
$\tau_{k}$ takes its values in $[-\infty,\infty)+\sqrt{-1}\RR$. In particular, we have
the claimed extension $\widetilde{B^+}\to \BB^+$. It sends
the boundary stratum $\tilde S(L_\pt)$ to the Borel-Serre boundary 
(for $\re (\tau_k)$ takes there the value $-\infty$) with image herein
the locus defined by putting all but the last three slots 
equal zero. The fiber passing through $\tilde z$  is locally 
given by putting $\tau_{k-1}=0$ and fixing the values of $F_k$ and $\tau_k$ (with
real part $-\infty$). In particular, this fiber is smooth at $\tilde z$.
This is true everywhere, and hence a connected
component of that fiber is also an irreducible component. Let us denote
the irreducible component passing through $\tilde z$ by $\Phi_{\tilde z}$.
So $\Phi_{\tilde z}$ lies over $\partial_k$.

If $k=1$, then $\Phi_{\tilde z}=\{\tilde z\}$ and the extension is at $\tilde z$ simply
given by $(1,F_1,\tau_1)$ and hence is there a local isomorphism.
If $k>1$, then since $(F_k,\im (\tau_k))$ defines a 
chart for the product $\PP((L_{k-1}/L_k)^+)\times \SS ((V/L_k)^+)$,
$\Phi_{\tilde z}$ is an irreducible component of a fiber of the natural map 
\begin{multline*}
\widetilde{\partial_kD(L_{k-1})}\to \partial_kD(L_{k-1})=
\PP (L_{k-1}^+)\times \PP((L_{k-1}/L_k)^+)\times \SS ((V/L_k)^+)\to\\
\to \PP((L_{k-1}/L_k)^+)\times \SS ((V/L_k)^+).
\end{multline*}
Since $L_{k-1}^\circ$ has finite longitudinal holonomy by Lemma 
\ref{lemma:collect}, the irreducible components of the fibers of this map are 
compact. If $\Phi_{\tilde z}=\{\tilde z\}$, then we must have 
$\dim L_{k-1}=1$. This implies that $k=2$ and that $(\tau_1, F_2,\tau_2)$ is a chart of
$\widetilde{B^+}$ at $\tilde z$ (we have $T_1=\{ 0\}$ in this case).
The  extension at $\tilde z$ is given by $(\tau_1^{-p_1},1, F_2,\tau_2)$.
Since $G_{L_1}$ acts on the first component as multiplication by $|p_1|$th roots 
of unity, we see that the extension  is at $\tilde z$ a local isomorphism 
modulo $G$. The proof of the proposition is now complete.
\end{proof}

\begin{proof}[Proof of Theorem \ref{thm:main}]
According to Proposition \ref{prop:collect}, the map $G\bs 
\widetilde{B^+}\to \BB^+$
has the property that the connected components of its fibers are compact,
that the preimage of the Baily-Borel boundary is in the boundary of the domain
and that where this map is locally finite it is in fact a local 
isomorphism.
So  Lemma \ref{lemma:stein} can be applied (in its entirity) to this situation
and we find that for the  topological Stein factorization  of 
$G\bs \widetilde{B^+}\to \BB^+$, 
\[
\begin{CD}
G\bs\widetilde{B^+} @>>> 
G\bs \widetilde{B^+}_\stein @>>> \BB^+,
\end{CD}
\]
the second map is a local isomorphism over $\BB$.
We first prove that $G\bs \widetilde{B}_\stein\to\BB$ is a 
$\G$-isomorphism.  For this we verify that the hypotheses of Lemma \ref{lemma:cov} 
are verified for the  Stein factor 
$G\bs\widetilde{B^+}_\stein$ with $Y':=\BB$.

We know that $\G$ acts properly discontinuously on
$\widetilde{B^+}$ with compact  fundamental domain. The 
first Stein factor is proper and $\G$-equivariant and so $\G$ acts also 
properly discontinuously on $G\bs\widetilde{B^+}_\stein$. 
Since $\G$ acts on $\BB$ as a group of isometries, Condition (ii)
of \ref{lemma:cov} is fulfilled  as well. The lemma
tells us that $G\bs\widetilde{B}_\stein\to\BB$
is then a covering projection. But $\BB$ is simply connected, 
and so this must be an isomorphism. 
It is easy to see that $G\bs\widetilde{B^+}_\stein\to\BB^+$
is then a $\G$-homeomorphism. Since $\G$ acts on the domain 
discretely and cocompactly, the same is true on its range.
This implies that $\G$ is discrete and of cofinite volume in the
unitary group of $h$.

The irreducible components of the preimages in $\tilde B$ of the 
exceptional divisors $D(L)$ are locally finite in $\widetilde{B}$;
since $\widetilde{B}\to G\bs\widetilde{B}_\stein$ is proper, the image of 
these in $\widetilde{B}_\stein$ are also locally finite. An irreducible
component $\tilde D(L)$ over $D(L)$ gets contracted if $\dim L>1$,
and its image in $\BB$ is the intersection of $\BB$ with a 
special subspace of codimension equal to the dimension of $L$.
The irreducible components of the preimages of the divisors $D(L)$
in $\widetilde{B^+}$ are locally finite. Hence their images
in $\BB$ are locally finite in $\BB$. We get a divisor precisely when
$\dim L=1$. It follows that the collection of special 
hyperplanes is locally finite on $\BB$, and that
$G\bs \PP (V^f)\subset G\bs B_\stein$ maps isomorphically 
onto the complement of the special hyperball arrangement
modulo $\G$, $\G\bs \BB^-$.

Since $G\bs \PP (V^f)\to \G\bs \BB^-$ is an isomorphism, so is
$G\bs V^f\to \G\bs (\LL^\times)^-$.
\end{proof}

\subsection{A presentation for the holonomy group}.
The holonomy group $\G$ is the image of a representation of the 
fundamental group $\pi_1(G\bs V^\circ ,*)$. In case $G$ is a 
Coxeter group and $\Hcal$ is its set of reflection hyperplanes, 
then $\pi_1(G\bs V^\circ ,*)$ is the Artin group of $G$ that we 
encountered in Subsection \ref{subsect:hecke}. But as the Lauricella systems 
show, $\Hcal$ may very well be bigger than 
the set of reflection hyperplanes of $G$. We describe a set of 
generators of the kernel of the holonomy representation and 
thus obtain a presentation of the holonomy group $\G$ in case 
we have one of $\pi_1(G\bs V^\circ ,*)$. 

Let us first note that any $L\in \Lcal(\Hcal)$ unambiguously 
determines a conjugacy class in the fundamental group of $V^\circ$:
blow up $L$ in $V$ and take the conjugacy class of a 
simple loop around the generic point of the exceptional divisor in 
(the preimage of) $V^\circ$. If we pass to the orbit space $G\bs V^\circ$, 
then $L^\circ$ determines a stratum in $G\bs V$. This stratum determines 
in the same way a conjugacy class in $\pi_1(G\bs V^\circ ,*)$.
If $L$ is irreducible and $\alpha_L\in \pi_1(G\bs V^\circ ,*)$ is a 
member of this conjugacy class, then
$\alpha_L^{|G_L|}$ is in the conjugacy class of $\pi_1(V^\circ ,*)$ defined 
above. If  $\kappa_L\not=1$, then the holonomy around this
stratum in $G\bs V^\circ$ has order $q_L$, where $q_L$ is the denominator
of $1-\kappa_L$. So  $\alpha_L^{q_L}$ is then the 
smallest power of $\alpha_L$ which lies in the kernel of the
monodromy representation. 

\begin{theorem}\label{thm:pres}
Suppose that we are in the elliptic, parabolic or hyperbolic case,
that is, in one the cases covered by Theorems \ref{thm:finite2}, \ref{thm:finitehol2}, 
\ref{thm:euclidean} and  \ref{thm:main}. Then $\G$ is obtained from
$\pi_1(G\bs V^\circ ,*)$ by imposing the relations 
$\alpha_L^{q_L}=1$ for
\begin{enumerate}
\item[(1)]  $L\in\Hcal$ and
\item[(2)] $L\in\Lcal_\irr (\Hcal)$ is of dimension $\le 1$ and $\kappa_L>1$.
\end{enumerate}
(Notice that for the complete elliptic and parabolic cases \ref{thm:finite2} and
\ref{thm:euclidean} the relations of the second kind do not occur.)
\end{theorem}
\begin{proof}
We limit ourselves to the hyperbolic case, since the others are easier.
Theorem \ref{thm:main} shows that $G\bs V^\circ$ can be identified with 
an open subset of $\G\bs\LL$.  
Since $\LL$ is a contractible (hence simply connected) complex 
manifold, $\G$ is the orbifold fundamental group of $\G\bs \LL$. 
Hence the quotient $\pi_1(G\bs V^\circ ,*)\to \G$ can be understood as
the map on (orbi)fold fundamental groups of the map 
$G\bs V^\circ\to\G\bs \LL$. It is well-known (and easy to see) that
the kernel of such a map is generated by the powers of the 
conjugacy classes in the fundamental group of $G\bs V^\circ$ defined by 
irreducible components of codimension one of the complement of the image,
$\G\bs \LL-G\bs V^\circ$, the power in question being the order
of local fundamental group at a general point of such an 
irreducible component. These irreducible components are
naturally indexed by the strata of $G\bs V$ of the type described in the 
theorem: the strata of codimension one of $G\bs V$ yield the irreducible
components meeting $G\bs V^f$, the zero dimensional stratum corresponds 
the image of the zero section $\G\BB\subset\G\bs \LL$ and the
strata of dimension one on which $\kappa >1$ correspond to the
remaining irreducible components. The powers are of course as stated in the theorem.
\end{proof}

\begin{remark}\label{rem:coxeterbessis}
Once we seek to apply Theorem \ref{thm:pres} in a concrete case, we 
need of course to have at our disposal a presentation of the fundamental group of 
$G\bs V^\circ$ in which the elements $\alpha_L$ can be identified. 
For $G$ a Coxeter group, this is furnished by the 
Brieskorn-Tits presentation \cite{briesk}, \cite{deligne:tresses};
this produces in the elliptic range the presentations of the associated complex reflection 
groups that are due to Coxeter \cite{coxeter}, Sections $12.1$ and 
$13.4$. For the case of an arbitrary finite complex reflection group, 
one may use a presentation of the fundamental group due to Bessis \cite{bessis}.
\end{remark}

\subsection{Automorphic forms and invariant theory}
According to Theorem \ref{thm:devmap}, the developing map  
$\widetilde{V^f}\to (\LL^\times)^-$ is homogeneous of 
negative degree $p_0$ (recall that $p_0$ is the numerator of the negative 
rational number $\kappa_0 -1$).
We can express this in terms of orbifold line bundles as follows:
if $\Ocal_{\G\bs \BB^-} (-1)$ denotes the $\G$-quotient of the automorphic
line bundle $\Ocal_{\BB^-} (-1)$ over $\G\bs \BB^-$, then 
the pull-back of this bundle over $\PP (V^f)$ is isomorphic to 
$\Ocal_{\PP (V)^f}(-p_0)$. Now $\PP (V)-\PP(V)^f$ is a closed subset of
$\PP(V)$ which is everywhere of codimension $>1$ and so for any $k\ge 0$,
the space of sections of $\Ocal_{\PP (V)^f}(k)$ is the space $\CC [V]_k$ of 
homogeneous polynomials on $V$ of degree $k$. We conclude that we have an
isomorphism of graded algebras
\[
\oplus_{n\ge 0} H^0(\BB^-,\Ocal (-n))^\G\cong 
\oplus_{n\ge 0}\CC [V]_{-np_0}^{G}.
\]
In particular, the lefthand side is finitely generated and its 
$\proj$ can be identified with $G\bs\PP (V)$.
In \cite{looij:ball} a systematic study was made of algebras of 
meromorphic automorphic forms of the type under consideration here.
The upshot is that the $\proj$ of the lefthand side is explicitly 
described as a modification of the Baily-Borel compactification
of $\G\bs \BB$ which leaves $\G\bs \BB^-$ untouched. 

To be more explicit, let us start out with the data consisting of the 
ball $\BB$, the group $\G$ , and the collection of special hyperplanes. 
Let us also make the rather modest assumption that $\dim V\ge 3$, 
so that $\dim\BB\ge 2$. The following lemma verifies the central hypothesis of
Corollary  5.8 of  \cite{looij:ball} (where the hermitian form is given the opposite 
signature).

\begin{lemma}
Every $1$-dimensional intersection of special hyperplanes is positive definite.
\end{lemma}
\begin{proof}
Any $1$-dimensional intersection $K$ of special 
hyperplanes which is negative semidefinite defines a point on the closure of $\BB$. 
If $K$ is negative (which defines an interior point of $\BB$), then $K$ is a 
special subspace and hence
corresponds to a member of $\Lcal^+(\Hcal)$ of codimension one, that is,
a member $H\in \Hcal$. Since $\kappa_H<1$, this is impossible. 
If $K$ is isotropic, then choose a $2$-dimensional intersection $P$ of special 
hyperplanes which contains $K$. Since the projectivization of $P$
meets $\BB$, it is a special subspace and hence corresponds to a member 
$L\in \Lcal^+(\Hcal)$ of codimension $2$. The transversal Dunkl system in $V/L$
has a projectivized developing map taking values in $\BB\cap \PP(P)$.
So $\Hcal_L$ contains a member $H$ with $\kappa_H=1$. But this  we excluded also.
\end{proof}

Although Corollary  5.8 of  \cite{looij:ball} does not apply as its 
stands---$\G$ need not be arithmetic---one can verify that the arguments to prove it
only require $\G$ to be discrete and of cofinite  volume in the relevant 
unitary group. It then tells us something we already know via our main theorem,
namely that the algebra of automorphic forms on $\BB$ with arbitrary poles
along the special hyperplanes is finitely generated with positive degree generators
and that the Proj of this graded algebra defines a certain projective completion
of $\G\bs \BB^-$: in the present situation the latter is just $\PP (G\bs V)$.
But in \cite{looij:ball} the completion is explicitly described as a blowup followed by a 
blowdown of the Baily-Borel  compactification of $\G\bs\BB$. If we go through the
details of this, we find that this intermediate blowup is almost $G\bs B^+$: 
the difference is that we now must blow up the parabolic
$L\in\Lcal^0(\Hcal)$ in the standard manner and not in the real-oriented sense.

\begin{question}
The algebra of $\G$-automorphic forms (of fractional degree) must appear in 
$\CC [V]^{G}$ as a subalgebra. It is in fact the subalgebra of
$G$-invariant polynomials which in degree $n$ vanish on each $L\in \Lcal^+(\Hcal)$ of order 
$\ge n(\kappa_L-1)/(\kappa_0-1)$.  It is only via our main theorem that we can give a geometric 
interpretation of the $\proj$ of this subalgebra as a modification of 
$\PP(G\bs V)$. In the Lauricella case, this can done
directly by means of geometric invariant theory, but is this possible in general?
\end{question}

\section{Classification of orbifolds for reflection arrangements}\label{section:examples}

Our aim is to list the Dunkl systems whose underlying arrangement is
that of a finite reflection group and for which the holonomy is as 
studied in the previous chapters: elliptic, parabolic or hyperbolic with
a discrete holonomy group of cofinite volume. More precisely, 
we classify the cases for which the hypotheses of the Theorems
\ref{thm:finite2}, \ref{thm:euclidean} and  \ref{thm:main} are satisfied. 

In order to display the information in an efficient way, we elaborate a
little on Remark \ref{rem:ab}. 
Given a Dunkl system of type  $A_n$ on $V=\CC^{n+1}/(\text{main 
diagonal})$ with the parameters $\mu_0, \ldots, \mu_n$, then 
for $m=0,\dots ,n$ we have a map 
\[ 
s_m :\CC^n\to V,\quad (u_1,\dots ,u_n)
\mapsto (u_1^2, \dots, u_{m-1}^2 ,0,u_m^2, \ldots, u_n^2).
\]
Remark \ref{rem:ab} tells us that pulling  back the Dunkl 
system along this map  yields a Dunkl system of type $B_n$;
we refer to this way of producing a $B_n$-system 
as \emph{reduction of the $A_n$-system at index $m$}.
Notice that any type $B_k$ subsystem of the $B_n$-system determines 
a $k+1$-element subset $I\subset\{0,\dots ,n\}$ which contains $m$
(and vice versa) with $\kappa$ taking the value 
$-1+2\mu_I$ on its fixed point subspace (where $\mu_I:=\sum_{i\in I} \mu_i$). 
On the other hand, any type $A_k$ subsystem is contained in a unique 
subsystem of type $B_{k+1}$ 
and so determines $(k+1)$-element subset of $J\subset\{0,\dots ,n\}-\{ m\}$; 
$\kappa$ takes then value $\mu_J$ on its fixed point subspace.

If we only wish to consider non-negative weights on arrangements, then
reduction at index $m$ is allowed only if 
$\frac{1}{2} \leq \mu_i + \mu_m < 1$ for all $i \neq m$.
Since the Dunkl system is invariant
under reflection in the short roots, we see that
the Schwarz condition on the weight $\kappa$ for a $B$-type 
intersection becomes: for all $I\ni m$, $1 - \mu_I$ is 
zero or the reciprocal of an integer.  In particular the weights 
on $B_n$ that satisfy the Schwarz conditions are all obtained 
by reduction at an index on $A_n$ that satisfies the Schwarz 
conditions.

The tables below list all the weights for arrangements of type 
$A$ and $B$ that satisfy the Schwarz conditions.  The parameters 
$\mu_i$ are defined by $n_i/d$ where $n_i$ and $d$ appear in the 
table.  If a parameter $n_m$ is typeset in bold then the weight 
obtained by reduction at position $m$ satisfies the Schwarz 
conditions for type $B$.  If additionally $n_i + n_m = d/2$ 
for all $i \neq m$ then the reduced weight can be considered 
as a weight on an arrangement of type $D$.  Note that such a 
weight is then invariant under the Weyl-group of type $D$.  
In the ``remark'' column ``ell'' stands for 
elliptic, ``par'' for parabolic and ``cc'' for co-compact.  
If no remark indicates otherwise, the group will be 
hyperbolic and acts with cofinite volume.

We omit the case $\kappa = 0$ from our tables.  
There is one additional series, corresponding to the full monomial 
groups, that is obtained as follows.  Take integers $n \geq 1$, 
$q \geq 2$ and define a weight on $A_n$ by $\mu_0 = \ldots = \mu_{n-1} = 0$, 
$\mu_{n} = 1-1/q$.  This weight can be reduced at index 
$n$ and satisfies both the Schwarz conditions for type $A$ and type $B$.

\tablecaption{Types $A_*$ and $B_*$}
\tablehead{\hline}
\tabletail{\hline}
\begin{supertabular}{r|rrrrrrrrrrr|c}
$\#$ & $d$ & $n_0$ & $n_1$ & $n_2$ & $n_3$ & $n_4$ & $n_5$ & $n_6$ & $n_7$ & $n_8$ & $n_{9}$ & remark \\
\hline
1      & 3    & \textbf{1} & \textbf{1} & \textbf{1} & \textbf{1} &  &  &  &  &  &  &   \\
2      & 4    & \textbf{1} & \textbf{1} & \textbf{1} & \textbf{1} &  &  &  &  &  &  & par \\
3      & 4    & \textbf{1} & \textbf{1} & \textbf{1} & \textbf{2} &  &  &  &  &  &  &    \\
4      & 5    & \textbf{2} & \textbf{2} & \textbf{2} & \textbf{2} &  &  &  &  &  &  & cc \\
5      & 6    & 1               & 1               & 1               & 1               &  &  &  &  &  &  & ell \\
6      & 6    & 1               & 1               & 1               & \textbf{2} &  &  &  &  &  &  & ell \\
7      & 6    & 1               & 1               & 1               & \textbf{3} &  &  &  &  &  &  & par \\
8      & 6    & 1               & 1               & 1               & \textbf{4} &  &  &  &  &  &  &       \\
9      & 6    & 1               & 1               & \textbf{2} & \textbf{2} &  &  &  &  &  &  & par \\
10    & 6    & 1               & 1               & \textbf{2} & \textbf{3} &  &  &  &  &  &  &\\
11    & 6    & \textbf{1} & \textbf{2} & \textbf{2} & \textbf{2} &  &  &  &  &  &  & \\
12    & 6    & \textbf{1} & \textbf{2} & \textbf{2} & \textbf{3} &  &  &  &  &  &  & \\
13    & 6    & \textbf{2} & \textbf{2} & \textbf{2} & \textbf{3} &  &  &  &  &  &  & \\
14    & 8    & \textbf{1} & \textbf{3} & \textbf{3} & \textbf{3} &  &  &  &  &  &  & cc \\
15    & 8 & \textbf{2} & \textbf{2} & \textbf{2} & \textbf{5} &  &  &  &  &  &  & cc \\
16 & 8 & \textbf{3} & \textbf{3} & \textbf{3} & \textbf{3} &  &  &  &  &  &  & cc \\
17 & 8 & \textbf{3} & \textbf{3} & \textbf{3} & \textbf{4} &  &  &  &  &  &  & cc \\
18 & 9 & \textbf{2} & \textbf{4} & \textbf{4} & \textbf{4} &  &  &  &  &  &  & cc \\
19 & 9 & \textbf{4} & \textbf{4} & \textbf{4} & \textbf{4} &  &  &  &  &  &  & cc \\
20 & 10 & \textbf{1} & \textbf{4} & \textbf{4} & \textbf{4} &  &  &  &  &  &  & cc \\
21 & 10 & \textbf{2} & 3 & 3 & 3 &  &  &  &  &  &  & cc \\
22 & 10 & \textbf{2} & 3 & 3 & \textbf{6} &  &  &  &  &  &  & cc \\
23 & 10 & 3 & 3 & 3 & 3 &  &  &  &  &  &  & cc \\
24 & 10 & 3 & 3 & 3 & \textbf{5} &  &  &  &  &  &  & cc \\
25 & 10 & 3 & 3 & 3 & \textbf{6} &  &  &  &  &  &  & cc \\
26 & 12 & \textbf{1} & \textbf{5} & \textbf{5} & \textbf{5} &  &  &  &  &  &  & cc \\
27 & 12 & 2 & 2 & 2 & \textbf{7} &  &  &  &  &  &  & cc \\
28 & 12 & 2 & 2 & 2 & \textbf{9} &  &  &  &  &  &  & cc \\
29 & 12 & 2 & 2 & \textbf{4} & \textbf{7} &  &  &  &  &  &  &cc \\
30 & 12 & \textbf{2} & \textbf{4} & \textbf{4} & \textbf{7} &  &  &  &  &  &  & cc \\
31 & 12 & \textbf{3} & \textbf{3} & \textbf{3} & \textbf{5} &  &  &  &  &  &  & cc \\
32 & 12 & \textbf{3} & \textbf{3} & \textbf{3} & \textbf{7} &  &  &  &  &  &  & cc \\
33 & 12 & \textbf{3} & \textbf{3} & \textbf{3} & \textbf{8} &  &  &  &  &  &  & cc \\
34 & 12 & \textbf{3} & \textbf{3} & \textbf{5} & \textbf{5} &  &  &  &  &  &  & cc \\
35 & 12 & \textbf{3} & \textbf{3} & \textbf{5} & \textbf{6} &  &  &  &  &  &  &    \\
36 & 12 & \textbf{3} & \textbf{5} & \textbf{5} & \textbf{5} &  &  &  &  &  &  & cc \\
37 & 12 & \textbf{3} & \textbf{5} & \textbf{5} & \textbf{6} &  &  &  &  &  &  & cc \\
38 & 12 & \textbf{4} & \textbf{4} & \textbf{4} & \textbf{5} &  &  &  &  &  &  & \\
39 & 12 & \textbf{4} & \textbf{4} & \textbf{4} & \textbf{7} &  &  &  &  &  &  & \\
40 & 12 & \textbf{4} & \textbf{4} & \textbf{5} & \textbf{5} &  &  &  &  &  &  & cc \\
41 & 12 & \textbf{4} & \textbf{4} & \textbf{5} & \textbf{6} &  &  &  &  &  &  & cc \\
42 & 12 & \textbf{4} & \textbf{5} & \textbf{5} & \textbf{5} &  &  &  &  &  &  & cc \\
43 & 12 & \textbf{4} & \textbf{5} & \textbf{5} & \textbf{6} &  &  &  &  &  &  & cc \\
44 & 12 & \textbf{5} & \textbf{5} & \textbf{5} & \textbf{5} &  &  &  &  &  &  & cc \\
45 & 12 & \textbf{5} & \textbf{5} & \textbf{5} & \textbf{6} &  &  &  &  &  &  & cc \\
46 & 14 & \textbf{2} & 5 & 5 & 5 &   &  &  &  &  &  &cc \\
47 & 14 & 5 & 5 & 5 & 5 &   &  &  &  &  &  &cc \\
48 & 15 & \textbf{4} & \textbf{6} & \textbf{6} & \textbf{6} &  &  &  &  &  &  & cc \\
49 & 15 & \textbf{4} & \textbf{6} & \textbf{6} & \textbf{8} &  &  &  &  &  &  & cc \\
50 & 15 & \textbf{6} & \textbf{6} & \textbf{6} & \textbf{8} &  &  &  &  &  &  & cc \\
51 & 18 & \textbf{1} & \textbf{8} & \textbf{8} & \textbf{8} &  &  &  &  &  &  & cc \\
52 & 18 & \textbf{2} & 7 & 7 & 7 &   &  &  &  &  &  &cc \\
53 & 18 & \textbf{2} & 7 & 7 & \textbf{10} &  &  &  &  &  &  & cc \\
54 & 18 & 3 & 3 & 3 & \textbf{13} &   &  &  &  &  &  &cc \\
55 & 18 & 3 & 3 & 3 & \textbf{14} &   &  &  &  &  &  &cc \\
56 & 18 & \textbf{5} & 7 & 7 & 7 &   &  &  &  &  &  &cc \\
57 & 18 & 7 & 7 & 7 & 7 &   &  &  &  &  &  &cc \\
58 & 18 & 7 & 7 & 7 & \textbf{10} &   &  &  &  &  &  &cc \\
59 & 20 & \textbf{5} & \textbf{5} & \textbf{5} & \textbf{11} &  &  &  &  &  &  & cc \\
60 & 20 & \textbf{5} & \textbf{5} & \textbf{5} & \textbf{14} &   &  &  &  &  &  &cc \\
61 & 20 & 6 & 6 & 6 & \textbf{9} &   &  &  &  &  &  &cc \\
62 & 20 & 6 & 6 & 6 & \textbf{13} &   &  &  &  &  &  &cc \\
63 & 20 & 6 & 6 & \textbf{9} & \textbf{9} &  &  &  &  &  &  &cc \\
64 & 20 & 6 & 6 & \textbf{9} & \textbf{10} &  &  &  &  &  &  & cc \\
65 & 24 & 4 & 4 & 4 & \textbf{17} &  &  &  &  &  &  & cc \\
66 & 24 & 4 & 4 & 4 & \textbf{19} &  &  &  &  &  &  & cc \\
67 & 24 & \textbf{7} & \textbf{9} & \textbf{9} & \textbf{9} &  &  &  &  &  &  & cc \\
68 & 24 & \textbf{7} & \textbf{9} & \textbf{9} & \textbf{14} &  &  &  &  &  &  & cc \\
69 & 24 & \textbf{9} & \textbf{9} & \textbf{9} & \textbf{14} &  &  &  &  &  &  & cc \\
70 & 30 & 5 & 5 & 5 & \textbf{19} &  &  &  &  &  &  & cc \\
71 & 30 & 5 & 5 & 5 & \textbf{22} &  &  &  &  &  &  & cc \\
72 & 30 & 5 & 5 & 5 & \textbf{23} &  &  &  &  &  &  & cc \\
73 & 30 & 9 & 9 & 9 & \textbf{11} &  &  &  &  &  &  & cc \\
74 & 42 & 7 & 7 & 7 & \textbf{29} &  &  &  &  &  &  & cc \\
75 & 42 & 7 & 7 & 7 & \textbf{34} &  &  &  &  &  &  & cc \\
76 & 42 & \textbf{13} & 15 & 15 & 15 &  &  &  &  &  &  & cc \\
77 & 42 & 15 & 15 & 15 & \textbf{26} &  &  &  &  &  &  & cc \\
78 & 3   & \textbf{1} &\textbf{1} & \textbf{1} & \textbf{1} & \textbf{1} &  &  &  &  &  & \\
79 & 4 & \textbf{1} & \textbf{1} & \textbf{1} & \textbf{1} & \textbf{1} &  &  &  &  &   &\\
80 & 4 & \textbf{1} & \textbf{1} & \textbf{1} & \textbf{1} & \textbf{2} &  &  &  &  &   & \\
81 & 6 & 1 & 1 & 1 & 1 & 1 &  &  &  &  &  & ell \\
82 & 6 & 1 & 1 & 1 & 1 & \textbf{2} &  &  &  &  &  & par \\
83 & 6 & 1 & 1 & 1 & 1 & \textbf{3} &  &  &  &  &  & \\
84 & 6 & 1 & 1 & 1 & 1 & \textbf{4} &  &  &  &  &  & \\
85 & 6 & 1 & 1 & 1 & \textbf{2} & \textbf{2} &  &  &  &  &  & \\
86 & 6 & 1 & 1 & 1 & \textbf{2} & \textbf{3} &  &  &  &  &  & \\
87 & 6 & 1 & 1 & \textbf{2} & \textbf{2} & \textbf{2} &  &  &  &  &  & \\
88 & 6 & 1 & 1 & \textbf{2} & \textbf{2} & \textbf{3} &  &  &  &  &  & \\
89 & 6 & \textbf{1} & \textbf{2} & \textbf{2} & \textbf{2} & \textbf{2} &  &  &  &  &  & \\
90 & 6 & \textbf{1} & \textbf{2} & \textbf{2} & \textbf{2} & \textbf{3} &  &  &  &  &  & \\
91 & 6 & \textbf{2} & \textbf{2} & \textbf{2} & \textbf{2} & \textbf{3} &  &  &  &  &  & \\
92 & 8 & \textbf{1} & \textbf{3} & \textbf{3} & \textbf{3} & \textbf{3} &  &  &  &  &  & cc \\
93 & 8 & \textbf{3} & \textbf{3} & \textbf{3} & \textbf{3} & \textbf{3} &  &  &  &  &  & cc \\
94 & 10 & \textbf{2} & 3 & 3 & 3 & 3 &   &  &  &  &  & cc \\
95 & 10 & 3 & 3 & 3 & 3 & 3 &   &  &  &  &  & cc \\
96 & 10 & 3 & 3 & 3 & 3 & \textbf{6} &  &  &  &  &  & cc \\
97 & 12 & 2 & 2 & 2 & 2 & \textbf{7} &  &  &  &  &  & cc \\
98 & 12 & 2 & 2 & 2 & 2 & \textbf{9} &  &  &  &  &  &  cc \\
99 & 12 & 2 & 2 & 2 & \textbf{4} & \textbf{7} &  &  &  &  &  &  cc \\
100 & 12 & \textbf{3} & \textbf{3} & \textbf{3} & \textbf{3} & \textbf{5} &  &  &  &  &  &  \\
101 & 12 & \textbf{3} & \textbf{3} & \textbf{3} & \textbf{3} & \textbf{7} &  &  &  &  &  &  \\
102 & 12 & \textbf{3} & \textbf{3} & \textbf{3} & \textbf{5} & \textbf{5} &  &  &  &  &  &  cc \\
103 & 12 & \textbf{3} & \textbf{3} & \textbf{5} & \textbf{5} & \textbf{5} &  &  &  &  &  &  cc \\
104 & 4   & \textbf{1} & \textbf{1} & \textbf{1} & \textbf{1} & \textbf{1} & \textbf{1} &  &  &  &  &\\
105 & 4 & \textbf{1} & \textbf{1} & \textbf{1} & \textbf{1} & \textbf{1} & \textbf{2} & &  &  &  &\\
106 & 6 & 1 & 1 & 1 & 1 & 1 & 1 &  &  &  &  &par \\
107 & 6 & 1 & 1 & 1 & 1 & 1 & \textbf{2} &  &  &  &  & \\
108 & 6 & 1 & 1 & 1 & 1 & 1 & \textbf{3} &  &  &  &  &  \\
109 & 6 & 1 & 1 & 1 & 1 & 1 & \textbf{4} &  &  &  &  &  \\
110 & 6 & 1 & 1 & 1 & 1 & \textbf{2} & \textbf{2} &  &  &  &  &  \\
111 & 6 & 1 & 1 & 1 & 1 & \textbf{2} & \textbf{3} &  &  &  &  &  \\
112 & 6 & 1 & 1 & 1 & \textbf{2} & \textbf{2} & \textbf{2} &  &  &  &  &  \\
113 & 6 & 1 & 1 & 1 & \textbf{2} & \textbf{2} & \textbf{3} &  &  &  &  &  \\
114 & 6 & 1 & 1 & \textbf{2} & \textbf{2} & \textbf{2} & \textbf{2} &  &  &  &  &  \\
115 & 10 & 3 & 3 & 3 & 3 & 3 & 3 &  &  &  &  &  cc \\
116 & 12 & 2 & 2 & 2 & 2 & 2 & \textbf{7} &  &  &  &  &  cc \\
117 &   4 & \textbf{1} & \textbf{1} & \textbf{1} & \textbf{1} & \textbf{1} & \textbf{1} & \textbf{1} & & & &\\
118 & 6 & 1 & 1 & 1 & 1 & 1 & 1 & 1 & & & &\\
119 & 6 & 1 & 1 & 1 & 1 & 1 & 1 & \textbf{2} & & & &\\
120 & 6 & 1 & 1 & 1 & 1 & 1 & 1 & \textbf{3} & & & &\\
121 & 6 & 1 & 1 & 1 & 1 & 1 & 1 & \textbf{4} & & & &\\
122 & 6 & 1 & 1 & 1 & 1 & 1 & \textbf{2} & \textbf{2} & & & &\\
123 & 6 & 1 & 1 & 1 & 1 & 1 & \textbf{2} & \textbf{3} & & & &\\
124 & 6 & 1 & 1 & 1 & 1 & \textbf{2} & \textbf{2} & \textbf{2} & & & &\\
125 & 6 & 1 & 1 & 1 & 1 & 1 & 1 & 1 & 1 & & & \\
126 & 6 & 1 & 1 & 1 & 1 & 1 & 1 & 1 & \textbf{2}  & & &\\
127 & 6 & 1 & 1 & 1 & 1 & 1 & 1 & 1 & \textbf{3}  & & &\\
128 & 6 & 1 & 1 & 1 & 1 & 1 & 1 & \textbf{2} & \textbf{2}  & & &\\
129 & 6 & 1 & 1 & 1 & 1 & 1 & 1 & 1 & 1 & 1 & & \\
130 & 6 & 1 & 1 & 1 & 1 & 1 & 1 & 1 & 1 & \textbf{2} & & \\
131 & 6 & 1 & 1 & 1 & 1 & 1 & 1 & 1 & 1 & 1 & 1 \\

\end{supertabular}

\smallskip
Tables $2$--$5$ list all remaining cases for the arrangements 
of the exceptional real and complex reflection groups.  
The Shephard groups $G_{25}$, $G_{26}$ and $G_{32}$ 
are omitted because these are already covered by the 
tables for types $A_3$, $B_3$ and $A_4$ respectively.

Only in the $F_4$ case the group has 
more than one orbit in its mirror arrangement.
This number is then two, which means that its discriminant
has two irreducible components; we write $q_1$ and $q_2$
for the ramification indices along these components, 
while we use a single $q$ in all other cases.  
The weight $\kappa$ on the arrangement 
is obtained by setting $\kappa_H = 1 - 2/q_H$ where $q_H$ 
is the ramification index along the image of the mirror $H$ in the 
orbit space.

All listed cases correspond to a hyperbolic reflection 
group except $q_1=2$, $q_2=3$ for type $F_4$ which is of 
parabolic type.  If a number $q$ or $q_i$ is typeset in 
bold then the corresponding group acts co-compactly on a 
hyperbolic ball, otherwise it acts with co-finite volume.  
All the obtained hyperbolic groups for the \emph{real} 
exceptional root systems are arithmetic.

%\newpage
%\suppressfloats

\tablecaption{Types $E_n$}
\tablehead{\hline}
\tabletail{\hline}
\begin{supertabular}{r||l|l|l}
$n$ & 6 & 7 & 8 \\
$q$ & 3, 4 & 3 & 3 \\
\end{supertabular}

\tablecaption{Type $F_4$}
\begin{supertabular}{r||l|l|l|l}
\tablehead{\hline}
\tabletail{\hline}
$q_1$ & 2 & 3 & 4 & 6 \\
$q_2$ & 3, 4, \textbf{5}, 6, \textbf{8}, \textbf{12} & 3, \textbf{4}, 6, \textbf{12} & 4 & 6 \\
\end{supertabular}
\vskip 1ex
The case $q_1 = 2$, $q_2 = 3$ is of parabolic type.

\tablecaption{Types $H_n$}
\tablehead{\hline}
\tabletail{\hline}
\begin{supertabular}{r||l|l}
$n$ & 3 & 4 \\
$q$ & \textbf{3}, \textbf{4}, \textbf{5}, \textbf{10} & \textbf{3}, \textbf{5} \\
\end{supertabular}

\tablecaption{Shephard-Todd groups $G_n$}
\tablehead{\hline}
\tabletail{\hline}
\begin{supertabular}{r||l|l|l|l|l|l}
$n$ & 24 & 27 & 29 & 31 & 33 & 34 \\
$q$ & \textbf{3}, 4, \textbf{5}, 6, \textbf{8}, \textbf{12} & \textbf{3}, 4, \textbf{5} & 3, 4 & 3, \textbf{5} & 3 & 3 \\
\end{supertabular}

%\vfill\eject

%\section*{}
%\newpage
\suppressfloats

\end{document}